\documentclass[reqno,11pt]{amsart}
\usepackage{graphicx}
\usepackage{amsfonts,amsmath, amssymb}
\usepackage{bbm,mathrsfs,mathscinet}
\usepackage[colorlinks,linkcolor=blue,anchorcolor=red,citecolor=blue]{hyperref}
\usepackage{marginnote}
\usepackage{color}
\usepackage{placeins}
\usepackage{verbatim}
\usepackage{tikz}
\usetikzlibrary{decorations,arrows}
\allowdisplaybreaks

\numberwithin{equation}{section}

\newcommand{\R}{\mathbb{R}}

\makeatletter
\newcommand{\rmnum}[1]{\romannumeral #1}
\newcommand{\Rmnum}[1]{\expandafter\@slowromancap\romannumeral #1@}
\makeatother

\newtheorem{theorem}{Theorem}[section]
\newtheorem{corollary}[theorem]{Corollary}
\newtheorem{lemma}[theorem]{Lemma}
\newtheorem{proposition}[theorem]{Proposition}
\newtheorem{remark}{Remark}[section]

\def\pa{\partial}

\def\la{\langle}
\def\ra{\rangle}

\def\var{\varepsilon}
\def\iintxy{\int_{\mathbb{T}}\int_{0}^{1}}

\begin{document}
	\title[Steady Navier--Stokes equations with Navier friction boundary conditions]{Inviscid limit of steady Navier--Stokes equations with Navier friction boundary conditions in a periodic channel}
	
	\author[R. M. Chen]{Robin Ming Chen}
	\address{Department of Mathematics, University of Pittsburgh, Pittsburgh, PA 15260, USA}
	\email{mingchen@pitt.edu}
	
	\author[D. Wang]{Dehua Wang}
	\address{Department of Mathematics, University of Pittsburgh, Pittsburgh, PA 15260, USA}
	\email{dhwang@pitt.edu}

	\author[W. Wang]{Weiqiang Wang}
	\address{Department of Mathematics, University of Pittsburgh, Pittsburgh, PA 15260, USA}
	\email{wew179@pitt.edu}

	\begin{abstract}
    The vanishing viscosity limit of the steady incompressible Navier--Stokes equations is studied in a periodic channel with Navier friction boundary conditions, where the slip length is assumed to be a nonnegative power of the viscosity. Depending on this power, the problem falls into subcritical, critical, and supercritical regimes. In each regime, the Navier--Stokes solutions are constructed via the Prandtl expansion, and the validity of the expansion is rigorously justified. By the Prandtl--Batchelor theory, the leading Euler solution is a Couette flow, and the interaction between the interior and boundary-layer profiles is precisely characterized. Of particular interest is the critical regime, in which the leading order Prandtl profile solves a nonlinear steady Prandtl equation with a Robin-type boundary condition. The well-posedness of this problem is established using the von Mises transformation and the implicit function theorem, and an implicit Batchelor--Wood formula is derived that determines the coefficients of the leading order Couette flow from the boundary data. To the best of our knowledge, this is the first well-posedness result for the steady Prandtl equation with a Robin-type boundary condition. 
	\end{abstract}

	\date{\today}
	
	\maketitle
	
	\setcounter{tocdepth}{2}
	\tableofcontents
	
	\thispagestyle{empty}
	
	\section{Introduction and Main Results}
	
    In this paper, we consider the vanishing viscosity limit for the steady incompressible Navier-Stokes equations in a periodic channel:
    \begin{equation}\label{1.1}
		\left\{\begin{aligned}
			&u^{\varepsilon}\partial_{x}u^{\varepsilon}+v^{\varepsilon}\partial_{y}u^{\varepsilon}+\partial_{x}p^{\varepsilon}=\varepsilon^{2}\big(\partial_{xx}u^{\varepsilon}+\partial_{yy}u^{\varepsilon}\big),\\
			&u^{\varepsilon}\partial_{x}v^{\varepsilon}+v^{\varepsilon}\partial_{y}v^{\varepsilon}+\partial_{y}p^{\varepsilon}=\varepsilon^2\big(\partial_{xx}v^{\varepsilon}+\partial_{yy}v^{\varepsilon}\big),\\
			&\partial_{x}u^{\varepsilon}+\partial_{y}v^{\varepsilon}=0,
		\end{aligned}\right.
	\end{equation}
    with $(x,y)\in \Omega=\mathbb{T}\times (0,1)=:\mathbb{R}\slash (2\pi\mathbb{Z})\times (0,1)$,
	supplemented with Navier friction  boundary conditions
	\begin{equation}\label{1.1-1}
		\begin{aligned}
			&2\delta(D(\mathbf{u}^{\varepsilon})\cdot \vec{n}_{0})\cdot \vec{\tau}_{0}+\mathbf{u}^{\varepsilon}\cdot \vec{\tau}_{0}=f_{0}(x),\quad \mathbf{u}^{\varepsilon}\cdot \vec{n}_{0}=0,\quad \text{on $y=0$},\\
			&2\delta(D(\mathbf{u}^{\varepsilon})\cdot \vec{n}_{1})\cdot \vec{\tau}_{1}+\mathbf{u}^{\varepsilon}\cdot \vec{\tau}_{1}=f_{1}(x),\quad \mathbf{u}^{\varepsilon}\cdot \vec{n}_{1}=0,\quad \text{on $y=1$},
		\end{aligned}
	\end{equation}
    where $\mathbf{u}^{\varepsilon}=(u^{\varepsilon},v^{\varepsilon})$ denotes velocity, $\varepsilon^2$ is the viscosity, $D(\mathbf{u}^{\varepsilon})=\frac{1}{2}(\nabla \mathbf{u}^{\varepsilon}+(\nabla \mathbf{u}^{\varepsilon})^{\top})$ is the deformation tensor. The vectors $\vec{n}_{0}=(0,-1)^{\top}$, $\vec{n}_{1}=(0,1)^{\top}$ are the unit outward normal vectors on the boundaries $\{y=0\}$ and $\{y=1\}$, respectively, while $\vec{\tau}_{0}=\vec{\tau}_{1}=(1,0)^{\top}$ are the corresponding unit tangential vectors. By straightforward calculations, the boundary conditions \eqref{1.1-1} can be rewritten as
\begin{align}\label{1.2}
		\left\{
		\begin{aligned}
			&(u^{\varepsilon}-\delta\partial_{y}u^{\varepsilon},v^{\varepsilon})\vert_{y=0}=(f_{0}(x),0),\quad x\in \mathbb{T},\\
			&(u^{\varepsilon}+\delta\partial_{y}u^{\varepsilon},v^{\varepsilon})\vert_{y=1}=(f_{1}(x),0),\quad x\in \mathbb{T}.
		\end{aligned}
		\right.
	\end{align}
	Here, $\delta>0$ is the scalar coefficient that measures the tendency of a fluid to slip over the boundary (in the literature, it is called the slip length or the inverse of friction coefficient). The boundary condition \eqref{1.1-1} was first proposed
	by Navier \cite{N} and later derived by Maxwell \cite{Maxwell} from the kinetic theory of rarefied gases. It was rigorously justified in \cite{JM-2001} as an effective boundary condition for flows over rough boundaries. From \eqref{1.2}, one recovers the Dirichlet boundary conditions as $\delta\to 0$, and Dirichlet-Neumann-type boundary conditions when $\delta\to+\infty$.

    To obtain a complete expansion of solutions in terms of the viscosity, we assume that the slip length is a power of $\var$, namely, 
    $
    \delta=\varepsilon^{\alpha},
    $
    then the boundary conditions \eqref{1.2} become
	\begin{align}\label{1.3}
		\left\{
		\begin{aligned}
			&(u^{\varepsilon}-\varepsilon^{\alpha}\partial_{y}u^{\varepsilon},v^{\varepsilon})\vert_{y=0}=(f_{0}(x),0),\quad x\in \mathbb{T},\\
			&(u^{\varepsilon}+\varepsilon^{\alpha}\partial_{y}u^{\varepsilon},v^{\varepsilon})\vert_{y=1}=(f_{1}(x),0),\quad x\in \mathbb{T}.
		\end{aligned}
		\right.
	\end{align}
Formally, as $\varepsilon\to 0$, the system \eqref{1.1}--\eqref{1.3} converges to the steady incompressible Euler equations:
\begin{equation}\label{1.1-2}
		\left\{\begin{aligned}
			&u\partial_{x}u+v\partial_{y}u+\partial_{x}p=0,\\
			&u\partial_{x}v+v\partial_{y}v+\partial_{y}p=0,\\
			&\partial_{x}u+\partial_{y}v=0,
		\end{aligned}\right.
	\end{equation}
supplemented with the boundary conditions
\begin{align}\label{1.3-1}
v\vert_{y=0,1}=0.
\end{align}

Due to the mismatch between the tangential boundary conditions in \eqref{1.3} and \eqref{1.3-1}, boundary layers are expected to form. Consequently, justifying the vanishing viscosity limit of Navier-Stokes in the presence of physical boundaries is a fundamental and challenging problem \cite{OS-1999,SG-2018}. For the unsteady Navier-Stokes equations with the no-slip boundary condition, Prandtl \cite{P-1904} proposed in 1905 that, for incompressible fluid in the small viscosity regime, the flow near the boundary is governed by the Prandtl boundary layer equation, whereas
the flow away from the boundary is approximated by the inviscid Euler equations; schematically,
$$
\text{Incompressible NS = Euler + Prandtl layer + $O(\var)$}.
$$
Since then, the well-posedness of the Prandtl equation and the rigorous justification of Prandtl boundary theory have become two core topics in mathematical fluid mechanics. Actually, there has been a lot of important progress on these topics. 
For the well-posedness of the Prandtl equation, see  \cite{AWXY,GM-2015,LMY,MW-2015,Oleinik-1967,PZ-2021,XZ-2004}. The ill-posedness and instability are studied in \cite{G-D,GGT,GT,GT-2011}. For steady and other related boundary-layer problems, see \cite{GM-2019,Iyer-2019,QDW,TW-2002,WZ-2021} and references therein. For results on the validity of unsteady Prandtl expansions and the inviscid limit for Navier-slip boundary conditions, see \cite{Beirao-Crispo,Berselli-Spirito,Chemetov-Cirpriano,CMR,GLNR,IP-2006,IS-2011,LNP,NP,Nguyen-2020,TWZ,WWX-2007,WWX-2011,WYZ-2017,WXZ,XX-2007,Zhong}; for Dirichlet boundary conditions, see \cite{FTZ-2018,FMT-2008,M-2014,MM-2017,SC-1998-1,SC-1998-2,WWZ-2017}.

For the steady case, since the steady Euler equations \eqref{1.1-2} admit a wide variety of solutions, the inviscid limit for the steady Navier--Stokes equations has been studied extensively in settings where the leading Euler flow is a simple shear flow or a small perturbation of shear flow. The first important progress on the validity of steady Prandtl expansion was made by Guo and Nguyen \cite{GT-2017} for the incompressible flow over a moving plate $(0,L)\times\R_{+}$. Later, S. Iyer \cite{Iyer-2017,Iyer-PK} extended that result to a moving sector and a moving half-plane $[1,\infty)\times(0,\infty)$, respectively; see also \cite{Iyer-2019} for the settings in which the leading Euler flow is a perturbation of a shear flow. For no-slip boundary conditions, Guo and Iyer \cite{GI-2021,GI-2023} justified the validity of steady Prandtl expansions on $[0,L]\times(0,\infty)$ for sufficiently small $L$. Later, Gao and Zhang \cite{GZ-2023} extended the analysis to arbitrary finite $L$ under suitable assumptions on the leading Euler shear flow and Prandtl profiles, see also \cite{GZ-2024} for the case when the Euler is a non-shear flow with sufficiently small $L$. We also refer to \cite{GW-arXiv} for steady Prandtl expansions for compressible Navier--Stokes equations.

When the tangential direction is periodic rather than finite, an additional difficulty arises: one should select the possible leading Euler flow and the higher-order interior solutions by the compatibility conditions of the corresponding equations of boundary layer solutions. In fact, the Prandtl--Batchelor theory, which dates back to Prandtl \cite{P-1904} and was rediscovered by Batchelor \cite{B-1956}, predicts that, for a two-dimensional steady flow at small viscosity regime, the vorticity is constant in a region filled with closed streamlines. Therefore, this selects a leading Euler flow with constant vorticity, while its coefficients must be chosen so that the associated leading Prandtl system is solvable. The corresponding selection relations include the Batchelor--Wood formula \cite{Wood} for a disk and the Feynman--Lagerstrom criterion \cite{DIN,FL-1956} for bounded domains with nonconstant curvature. Recently, Gui--Xie--Xu \cite{GXX} removed the closed-streamline assumption in the Prandtl--Batchelor theory. We also refer to \cite{EHSX,GXX-2026,HN-2017,HN-2019,HN-2023} for rigidity results on classification of solutions of the steady Euler equation.

Recent progress on the Prandtl--Batchelor theory includes the works of Fei--Gao--Lin--Tao \cite{FGLT-2023,FGLT-2024} and Pan--Zhao \cite{PZ-2025-2} for a disk, an annulus, and the exterior of a disk, respectively, in regimes where the slip boundary data are small perturbations of constants. See also \cite{GL-2025} for a moving periodic half-space $\mathbb{T}_{\rho}\times(0,\infty)$ with sufficiently small period $\rho$. Total--slip boundary conditions in the periodic channel $\mathbb{T}\times(0,1)$, with or without external forcing, were studied by Song--Tao \cite{ST-2027} and Pan--Zhao \cite{PZ-2025}, respectively. The mixed case with a slip condition at one boundary and a no-slip condition at the other in a periodic channel and an expanding annulus was considered in \cite{FPZ-2025} and \cite{PZZ-2026}, respectively.

In this paper, we investigate the validity of Prandtl expansion and the Prandtl--Batchelor theory of \eqref{1.1} subject to the boundary conditions \eqref{1.3}. To have a clear Prandtl expansion, we assume
\begin{align}\label{alpha}
\alpha\geq 0, \quad  |1-\alpha|\in \mathbb{Z}\text{ or }\frac{1}{|1-\alpha|}\in \mathbb{Z}.
\end{align}
The general $\alpha$ will be discussed in Appendix B. For later use, we denote
\begin{align}\label{talpha}
\tilde{\alpha}=\left\{
\begin{aligned}
&\frac{1}{\min\{|1-\alpha|,1\}},\,\,\quad \text{if $\alpha\neq 1$},\\
&1,\qquad\qquad\qquad\quad\,\,\,\,\text{if $\alpha=1$}.
\end{aligned}
\right.
\end{align}

Throughout the paper, for a function $f=f(x,y)$, we use the notation
$$
\bar{f}(y):=\frac{1}{2\pi}\int_{\mathbb{T}}f(x,y)\,{\rm d}x.
$$
Similarly, for a function $f=f(x)$, we write
$$
\bar{f}:=\frac{1}{2\pi}\int_{\mathbb{T}}f(x){\rm d}x.
$$
We assume that the boundary data $f_0(x)$ and $f_1(x)$ are smooth periodic functions satisfying $\bar{f}_{0}\bar{f}_{1}>0$, which will ensure that the horizontal velocity of the leading Euler flow has no stagnation point. Thus, without loss of generality, we may assume 
\begin{align}\label{bd}
\bar{f}_{0}>0,\quad \bar{f}_{1}>0.
\end{align}
for simplicity of presentation.

Our main theorems are stated as follows:
\begin{theorem}\label{thm1}
For $0\leq \alpha<1$ satisfying \eqref{alpha}, let $f_{0}(x)$ and $f_{1}(x)$ satisfy \eqref{bd}.
Then there exists a constant $\varepsilon_{0}>0$ such that for any $0<\varepsilon<\varepsilon_{0}$,
the Navier–Stokes equations \eqref{1.1} with boundary conditions \eqref{1.3} have a strong solution $(u^{\varepsilon},v^{\varepsilon},p^{\varepsilon})$ satisfying
\begin{align*}
	\|u^{\varepsilon}-(Ay+B)\|_{L^{\infty}}+\|v^{\varepsilon}\|_{L^{\infty}}\leq C\varepsilon^{1-\alpha}.
\end{align*}
Here, the constants $A$ and $B$ in the leading Euler shear flow are determined as follows:
\begin{itemize}
\item [$\bullet$] For $\alpha=0$,
$$
A=\frac{1}{3}\bar{f}_{1}-\frac{1}{3}\bar{f}_{0},\quad B=\frac{2}{3}\bar{f}_{0}+\frac{1}{3}\bar{f}_{1},
$$

	\item [$\bullet$] For $0<\alpha<1$, 
    $$
    A=\bar{f}_{1}-\bar{f}_{0},\quad B=\bar{f}_{0}.
    $$
\end{itemize}
\end{theorem}

\begin{theorem}\label{thm2}
    For $\alpha\geq 1$ satisfying \eqref{alpha}, let $f_{0}(x)$ and $f_{1}(x)$ satisfy \eqref{bd} and be small perturbations around their means:
\begin{align}\label{small-per}
\begin{aligned}
f_{0}(x)=\bar{f}_{0}+\sigma g_{0}(x),\quad f_{1}(x)=\bar{f}_{1}+\sigma g_{1}(x)
\end{aligned}
\end{align}
for some small constant $\sigma$. Then there exist constants $\varepsilon_{0}$ and $\sigma_{0}$ such that for any $0<\varepsilon<\varepsilon_{0}$ and $\sigma\in [0,\sigma_{0})$ the Navier–Stokes equations \eqref{1.1} with boundary conditions \eqref{1.3} have a strong solution $(u^{\varepsilon},v^{\varepsilon},p^{\var})$ satisfying
\begin{align*}
	\Big\| u^{\varepsilon}-(Ay+B)-u_{B}^{0}\Big(x,\frac{y}{\varepsilon}\Big)\chi(y)-\hat{u}_{B}^{0}\Big(x,\frac{y-1}{\varepsilon}\Big)\tilde{\chi}(y)\Big\|_{L^{\infty}}+\|v^{\varepsilon}\|_{L^{\infty}}\leq
\left\{
\begin{aligned}
&C\varepsilon,\qquad \text{for $\alpha=1$},\\
&C\varepsilon^{\frac{1}{\tilde{\alpha}}},\quad\,\text{for $\alpha>1$},
\end{aligned}
  \right.
\end{align*}
Here, $\chi(y)$ is a smooth and non-increasing cutoff function satisfying
\begin{align}\label{cut-off}
\chi(y)=\left\{
\begin{aligned}
&1,\quad \text{for $0<y\leq\frac{1}{4}$},\\
&0,\quad \text{for $\frac{1}{2}\leq y<1$},
\end{aligned}
\right.
\end{align}
and $\tilde{\chi}(y)=1-\chi(y)$. The functions $u_{B}^{0}(x,\frac{y}{\varepsilon})$ and $\hat{u}_{B}^{0}(x,\frac{y-1}{\varepsilon})$ are the leading boundary-layer profiles for the bottom and upper boundaries, respectively. 
Here, the constants $A$ and $B$ in the leading Euler shear flow are determined as follows:
\begin{itemize}
\item [$\bullet$] For $\alpha=1$, $B\geq \bar{f}_{0}$, $A+B\geq \bar{f}_{1}$, and $ B-\bar{f}_{0}$, $(A+B)-\bar{f}_{1}$ are constants of order $O(\sigma^2)$ determined implicitly by
    $$
     (A+B)-\bar{f}_{1}=\frac{\sigma^2}{4\pi(A+B)}\int_{\mathbb{T}}\hat{\phi}^2(\xi,0)\,{\rm d}\xi,\quad B-\bar{f}_{0}=\frac{\sigma^2}{4\pi B}\int_{\mathbb{T}}\phi^2(\xi,0)\,{\rm d}\xi,\quad 
    $$
    where $\phi(\xi,\eta)$, $\hat{\phi}(\xi,\eta)$ are the smooth solutions of \eqref{lowerboundary} and \eqref{upperboundary}, respectively.
	\item [$\bullet$] For $\alpha>1$, $B>0$, $(A+B)>0$ and
	\begin{align}\label{W1}
    (A+B)^2=\overline{f_{1}^2},\quad B^2=\overline{f_{0}^2}.
    \end{align}
\end{itemize}
\end{theorem}

\begin{theorem}\label{thm3}
Let $(u^{\var},v^{\var},p^{\var})$ be the strong solution of the Naver-Stokes equations \eqref{1.1} constructed in Theorem \ref{thm1} or Theorem \ref{thm2}, and $\omega^{\var}=\partial_{y}u^{\var}-\partial_{x}v^{\var}$ be the vorticity. Denote $\Omega_{d}=\mathbb{T}\times (d,1-d)$ for $0<d<\frac{1}{2}$. Then we have
\begin{align}\label{vorticity}
\|\omega^{\var}-A\|_{H^{1}(\Omega_{d})}+\|\omega^{\var}-A\|_{L^{\infty}(\Omega_{d})}\leq C_{d}\varepsilon^{\frac{1}{\tilde{\alpha}}},
\end{align}
where $A$ is the constant determined in Theorem \ref{thm1} or Theorem \ref{thm2}, and $C_{d}>0$ is a constant depending on $d$ but is independent of $\varepsilon$. 
\end{theorem}

\begin{remark}
  By the Prandtl--Batchelor theory, the Euler flow in the inviscid limit of Navier-Stokes equations on $\Omega=\mathbb{T}\times(0,1)$ should be a Couette flow, see \cite{B-1956,FGLT-2023,FGLT-2024,GXX,PZ-2025}. Indeed, it follows from \cite[Remark 1.5]{GXX} that the periodicity of pressure $p^{\var}$ rules out the Poiseuille flow.
\end{remark}

\begin{remark}
    The general $\alpha$ will be briefly discussed in Appendix B. Specifically, the Prandtl expansion also holds for every nonnegative rational $\alpha$.
\end{remark}

\begin{remark}
For the regime $0\leq \alpha<1$, the boundary data do not need to be a small perturbation of their means. That is mainly because the leading boundary-layer profiles vanish.

For the regime $\alpha>1$, the choices of $A$ and $B$ in \eqref{W1} are the well known Batchelor--Wood formulas, see \cite{B-1956,Wood}, which are derived from compatibility conditions for the leading Prandtl equations with the Dirichlet boundary conditions.

For the regime $\alpha=1$, the leading boundary-layer profiles satisfy the Prandtl equations with a Robin-type boundary condition. A corresponding Batchelor--Wood--type relation is required for compatibility, but the relation between the boundary data and the leading Euler flow is implicit. See Section 4.1 for details.
\end{remark}

\begin{remark}
Let 
$$u_{p}^{0}(x,y)=Ay+B+u_{B}^{0}\Big(x,\frac{y}{\varepsilon}\Big)\chi(y)+\hat{u}_{B}^{0}\Big(x,\frac{y-1}{\var}\Big)\tilde{\chi}(y)$$
be the leading Prandtl profile. Then, we will have
\begin{align}\label{p1}
\overline{u_{p}^{0}(x,0)}=B=\bar{f}_{0},\quad \overline{u_{p}^{0}(x,1)}=A+B=\bar{f}_{1}\quad \text{for $0<\alpha<1$},
\end{align}
and
\begin{align}\label{p2}
\overline{u_{p}^{0}(x,0)}=\bar{f}_{0},\quad \overline{u_{p}^{0}(x,1)}=\bar{f}_{1}\quad \text{for $\alpha\geq 1$}.
\end{align}
If we investigate the average of $|u_{p}^{0}|^2$ on the boundaries, we have
\begin{align}\label{p3}
&\overline{|u_{p}^{0}(x,0)|^2}= B^2=|\bar{f}_{0}|^2,\quad \overline{|u_{p}^{0}(x,1)|^2}=(A+B)^2=|\bar{f}_{1}|^2,\quad \text{for $0<\alpha<1$},
\end{align}
\begin{align}\label{p4-0}
&\overline{|u_{p}^{0}(x,0)|^2}=B^2\geq |\bar{f}_{0}|^2,\quad \overline{|u_{p}^{0}(x,1)|^2}=(A+B)^2\geq |\bar{f}_{1}|^2,\quad \text{for $\alpha=1$},
\end{align}
and
\begin{align}\label{p5}
&\overline{|u_{p}^{0}(x,0)|^2}=B^2=\overline{f_{0}^2},\quad \overline{|u_{p}^{0}(x,1)|^2}=(A+B)^2=\overline{f_{1}^2},\quad \text{for $\alpha>1$}.
\end{align}
The corresponding relations for $0<\alpha<1$ can be directly derived from the choices of $A$ and $B$ in Theorem \ref{thm1}. For $\alpha>1$, the corresponding relation \eqref{p2} follows directly from the boundary conditions of $u_{p}^{0}$ (see \eqref{b3} below) and \eqref{p5} follows from the choices of $A$ and $B$ in Theorem \ref{thm2}. For $\alpha=1$, the corresponding relations follow from \eqref{cw4} below. 

Thus, for $\alpha>0$, it follows from \eqref{p1}--\eqref{p2} that the mean of the boundary data, not the mean of leading Euler flow on the boundaries, characterizes the mean of the leading Prandtl profile on the boundaries, and equations \eqref{p3}--\eqref{p5} show that the $L^2$ norm of the leading Euler flow on the boundaries, not the $L^2$ norm of boundary data, characterizes the $L^2$ norm of the leading Prandtl profile on the boundaries.
\end{remark}

Now, we briefly mention the main ingredients in our proof:
\begin{itemize}
    \item[(\rmnum{1})] Both horizontal and vertical velocities in our interior expansions are harmonic by delicate constructions, and the interplay between their coefficients and the corresponding boundary-layer problem is carefully presented.
    \item[(\rmnum{2})] The leading boundary-layer profiles satisfy a nonlinear Prandtl-type equation with a Robin-type boundary condition in the critical regime $\alpha=1$. The well-posedness and the selection principle on the leading Euler flow are established through the von Mises transformation and a delicate implicit function theorem.
    \item[(\rmnum{3})] Two different elegant decompositions of the horizontal velocities are used in the {\it a priori} uniform estimates of the remainder equations to adapt the slip effect in the boundary conditions for $0\leq \alpha<1$ and $\alpha\geq 1$, respectively.
\end{itemize}
We explain the above three ingredients in more details as follows.

(\rmnum{1}) The classical Prandtl expansion shows that the leading interior profile satisfies the Euler equations and, by the Prandtl--Batchelor theory, must be a Couette flow $(Ay+B,0)$. The higher-order interior solutions will satisfy the linearized Euler equation around this Couette flow. Exploiting the structure of these equations, we show that both horizontal and vertical velocities can be chosen to be harmonic; equivalently, all higher-order interior profiles have constant vorticities. Moreover,
once the vertical velocity is determined,
we observe that the horizontal velocity is constrained only by incompressibility and does not affect the construction of the pressure as it appears in its $x-$derivative in the equations. Consequently, an auxiliary linear function $C_{1}y+C_{2}$ may be added to the horizontal velocity without changing the equations for the vertical velocity or the pressure. This freedom is essential to ensure the compatibility conditions of the corresponding boundary layer equations. More precisely, for $0\leq \alpha<1$, the constants $C_{1}$ and $C_{2}$ are determined by the solvability conditions for the corresponding boundary-layer profiles, while for $\alpha\geq1$ they are chosen to absorb the far-field constants of those profiles; see Proposition \ref{prop3} for details.

(\rmnum{2}) For $0\leq \alpha<1$, the leading boundary-layer profiles vanish identically due to homogeneous boundary conditions and the uniqueness of solution to the Prandtl equations. 

For $\alpha>1$, the leading boundary-layer profiles satisfy a nonlinear Prandtl-type equation with Dirichlet boundary data. That is, if we introduce the Prandtl profile $u$ (at the lower boundary for example), it will hold that
\begin{align*}
    \left\{
\begin{aligned}
    &u\partial_{x}u+v\partial_{z}u=\partial_{zz}u,\\
    &\partial_{x}u+\partial_{z}v=0,\\
    &v\vert_{z=0}=0,\quad \lim\limits_{z\to \infty}u=B,
\end{aligned}
    \right.
\end{align*}
with the Dirichlet boundary condition
$u\vert_{z=0}=f_{0}(x)$,
where $z=\frac{y}{\var}$ is the Prandtl scaling at the lower boundary.

Applying the von Mises transformation:
\begin{align*}
\xi=x,\quad \eta=\int_{0}^{z}u(x,s)\,{\rm d}s,\quad w(\xi,\eta)=u^2(x(\xi,\eta),z(\xi,\eta)),
\end{align*}
we obtain
\begin{align}\label{vonMise1-1}
\left\{
\begin{aligned}
&\frac{1}{\sqrt{w}}\partial_{\xi}w=\partial_{\eta\eta}w,\\
&w(\xi,0)=f_{0}^2(\xi),\quad \lim\limits_{\eta\to \infty}w=B^2.
\end{aligned}
\right.
\end{align}
Integrating \eqref{vonMise1-1} over $\xi\in \mathbb{T}$ gives $\bar{w}''=0$, which, together with the boundary condition, yields $\bar{w}=B^2$. In particular, at $\eta=0$, we obtain the Batchelor-Wood formula
\begin{align}\label{Batchelor--Wood}
\overline{f_{0}^2}=B^2.
\end{align}
Under the compatibility condition \eqref{Batchelor--Wood} and the assumption \eqref{small-per}, the well-posedness of \eqref{vonMise1-1} was established in \cite{FGLT-2023,FGLT-2024,GL-2025,GT-2017}.

In the critical case $\alpha=1$, however, the leading Prandtl profile satisfies the Robin boundary condition
\begin{align*}
(u-\partial_{z}u)\vert_{z=0}=f_{0}(x).
\end{align*}
For this boundary condition, an analogue of the compatibility condition \eqref{Batchelor--Wood} is not easy to identify in terms of $w=u^2$. We therefore work instead with
\begin{align*}
\mathfrak{w}(\xi,\eta)=u(x(\xi,\eta),z(\xi,\eta)),
\end{align*}
which satisfies
\begin{align*}
\left\{
\begin{aligned}
&\partial_{\xi}\mathfrak{w}=(\partial_{\eta}\mathfrak{w})^2+\mathfrak{w}\partial_{\eta\eta}\mathfrak{w},\\
&(\mathfrak{w}-\mathfrak{w}\partial_{\eta}\mathfrak{w})\vert_{\eta=0}=f_{0}(\xi),\quad \lim\limits_{\eta\to\infty}\mathfrak{w}=B.
\end{aligned}
\right.
\end{align*}
By a straightforward calculation, we find
\begin{align*}
\overline{\mathfrak{w}}(0)=\bar{f}_0, \qquad \overline{\mathfrak{w}^2}(0)=B^2.
\end{align*}
Unlike \eqref{Batchelor--Wood}, these identities do not yield an explicit relation between $f_{0}$ and $B$. This is the main difficulty in establishing well-posedness of Prandtl equation in the critical regime.

To overcome this difficulty, we assume \eqref{small-per}, and seek a solution of the form
$$
\mathfrak{w}=B+\sigma\phi,
$$
then $\phi$ satisfies
\begin{align}\label{per-equation}
\left\{
\begin{aligned}
  &\partial_{\xi}\phi=B\partial_{\eta\eta}\phi+\frac{\sigma}{2}\partial_{\eta\eta}(\phi^2),\\
  &(\phi-B\partial_{\eta}\phi-\frac{\sigma}{2}\partial_{\eta}(\phi^2))\vert_{\eta=0}=\frac{1}{\sigma}(\bar{f}_{0}-B)+g_{0}(\xi),\quad \lim\limits_{\eta\to \infty}\phi=0.
\end{aligned}
\right.
\end{align}
Denote the nonlinear operator
$$
\mathcal{F}(u;B,\sigma,q)=\mathcal{L}_{B}u-\sigma\mathcal{K}u-[0,q]^{\top},
$$
where
$$
\mathcal{L}_{B}u=[\partial_{\xi}u-B\partial_{\eta\eta}u,\Pi_{\perp}(u-B\partial_{\eta}u)\vert_{\eta=0}]^{\top},\quad \mathcal{K}u=\Big[\frac{1}{2}\partial_{\eta\eta}(u^2),\Pi_{\perp}\big(\frac{1}{2}\partial_{\eta}(u^2)\vert_{\eta=0}\big)\Big]^{\top},
$$
the notation $[a,b]^{\top}$ denotes the transpose of the vector $[a,b]$, and $\Pi_{\perp}f=f-\bar{f}(y)$. Now, we solve the abstract equation
\begin{align}\label{ae}
    \mathcal{F}(u;B,\sigma,q)=[0,0]^{\top}.
\end{align}
First, the linear operator $\mathcal{L}_{B}$ is shown to be an isomorphism using the classical parabolic theory. Next, it is direct to check that \eqref{ae} has a solution
$\phi^{0}:=\mathcal{L}_{\bar{f}_{0}}^{-1}[0,g_{0}]^{\top}$ for $(B,\sigma,q)=(\bar{f}_{0},0,g_{0})$. Then, applying the implicit function theorem to \eqref{ae}, we obtain a $C^{\infty}$ solution map
$$
(B,\sigma,q)\longmapsto U(B,\sigma,q)
$$
in a small neighborhood of $(\bar{f}_{0},0,g_{0})$, which is unique in a fixed small neighborhood of $\phi^{0}$. Setting $\Phi(B,\sigma)=U(B,\sigma,g_{0})$, one has that $\Phi(B,\sigma)$ satisfies
\begin{align}\label{per-equation2}
\left\{
\begin{aligned}
  &\partial_{\xi}\Phi-B\partial_{\eta\eta}\Phi-\frac{\sigma}{2}\partial_{\eta\eta}(\Phi^2)=0,\\
  &\Pi_{\perp}(\Phi-B\partial_{\eta}\Phi-\frac{\sigma}{2}\partial_{\eta}(\Phi^2))\vert_{\eta=0}=g_{0}(\xi), &\lim\limits_{\eta\to \infty}\Phi=0.
\end{aligned}
\right.
\end{align}
Comparing \eqref{per-equation} and \eqref{per-equation2}, we need to require
\begin{align}\label{B-sequence}
\frac{1}{\sigma}(\bar{f}_{0}- B)&=\frac{1}{2\pi}\int_{\mathbb{T}}(\Phi-B\partial_{\eta}\Phi-\frac{\sigma}{2}\partial_{\eta}(\Phi^2))(\xi,0)\,{\rm d}\xi\nonumber\\
&=\frac{1}{2\pi}\int_{\mathbb{T}}\Phi(\xi,0)\,{\rm d}\xi=-\frac{\sigma}{4\pi B}\int_{\mathbb{T}}\Phi^2(\xi,0)\,{\rm d}\xi.
\end{align}
The last two equalities come from the compatibility condition of \eqref{per-equation2}. However, since $\Phi$ depends on $B$, we must show that there exists indeed a constant $B$ that satisfies \eqref{B-sequence} in the small neighborhood of $\bar{f}_{0}$. 
To do that, we construct the following iteration
\begin{align*}
B_{k}=\bar{f}_{0}+\frac{\sigma^2}{4\pi B_{k-1}}\int_{\mathbb{T}}|\Phi(B_{k-1},\sigma)(\xi,0)|^2\,{\rm d}\xi,
\end{align*}
with $B_{0}=\bar{f}_{0}$. Then, based on delicate estimates on $\Phi(B,\sigma)$, we obtain that $\{B_{k}\}$ is a Cauchy sequence provided that $\sigma$ is sufficiently small, which yields that there exists a constant $B$ in the small neighborhood of $\bar{f}_{0}$ such that \eqref{B-sequence} holds. Finally, set $\phi^{\sigma}=\Phi(B,\sigma)$, then $\phi^{\sigma}$ is a solution of \eqref{per-equation}. It is worth pointing out that, in the above process, we need to ensure that the smallness of $\sigma$ can be determined independent of both the order of regularity of $\Phi(B,\sigma)$ and each iteration of $B_{k}$; see Lemma \ref{lem1} for details.

(\rmnum{3}) Once the interior and boundary-layer profiles have been constructed, it remains to control the remainder $(u_{R},v_{R},p_{R})$, which satisfies
\begin{align}\label{remainder}
\left\{
\begin{aligned}
 &u^{a}\partial_{x}u_{R}+v^{a}\partial_{y}u_{R}+u_{R}\partial_{x}u^{a}+v_{R}\partial_{y}u^{a}+\partial_{x}p_{R}=\var^2\Delta u_{R}+R_{u}^{a}-(u_{R}\partial_{x}u_{R}+v_{R}\partial_{y}u_{R}),\\
 &u^{a}\partial_{x}v_{R}+v^{a}\partial_{y}v_{R}+u_{R}\partial_{x}v^{a}+v_{R}\partial_{y}v^{a}+\partial_{y}p_{R}=\varepsilon^2\Delta v_{R}+R_{v}^{a}-(u_{R}\partial_{x}v_{R}+v_{R}\partial_{y}v_{R}),\\
 &\partial_{x}u_{R}+\partial_{y}v_{R}=0,\\
 &(u_{R}-\varepsilon^{\alpha}\partial_{y}u_{R},v_{R})\vert_{y=0}=(0,0),\quad (u_{R}+\varepsilon^{\alpha}\partial_{y}u_{R},v_{R})\vert_{y=1}=(0,0),
\end{aligned}
\right.
\end{align}
where $(u^{a},v^{a})$ is the approximate solution constructed by the Prandtl expansion, and
\begin{align*}
&R_{u}^{a}=-(u^{a}\partial_{x}u^{a}+v^{a}\partial_{y}u^{a}+\partial_{x}p^{a})+\varepsilon^2 \Delta u^{a},\\
&R_{v}^{a}=-(u^{a}\partial_{x}v^{a}+v^{a}\partial_{y}v^{a}+\partial_{y}p^{a})+\varepsilon^2\Delta v^{a},
\end{align*}
are the approximation error terms. As in \cite{FGLT-2023,FGLT-2024}, we use the Banach fixed point theory, together with basic energy estimates, positivity estimates and Stokes' estimates, to show the existence of \eqref{remainder}.

For the basic energy estimates, we  multiply $\eqref{remainder}_{1}$ by $u_{R}$, and $\eqref{remainder}_{2}$ by $v_{R}$, and then add the resultant equations together to obtain the following good terms:
$$
\varepsilon^2\|(\nabla u_{R},\nabla v_{R})\|_{L^2}^2+\varepsilon^{2-\alpha}\int_{\mathbb{T}}\big(|u_{R}(x,0)|^2+|u_{R}(x,1)|^2\big)\,{\rm d}x.
$$
For the bad terms, the most singular term is
$$
\int_{\mathbb{T}}\int_{0}^{1}u_{R}v_{R}\partial_{y}u^{a}\,{\rm d}y{\rm d}x,
$$
because, $\partial_{y}u^{a}$ behaves like $\varepsilon^{-\alpha}$ for $0\leq \alpha<1$, $\varepsilon^{-1}$ for $\alpha\geq 1$ when $\partial_{y}$ is applied to the leading boundary-layer profile $u_{B}^{a}$. Moreover, because $u_R$ does not vanish on the boundary, the Hardy inequality cannot be applied directly to $u_R$ to gain a factor of $\varepsilon$. We therefore use different decompositions in the regimes $\alpha<1$ and $\alpha\geq 1$, respectively to overcome this difficulty.

For $0\leq \alpha<1$, we rewrite $u_{R}$ as
\begin{align}\label{decom1}
u_{R}(x,y)&=\bar{u}_{R}(0)+(\bar{u}_{R}(y)-\bar{u}_{R}(0))+\tilde{u}_{R}(x,y)\nonumber\\
&=:\bar{u}_{R}(0)+\hat{u}_{R}(y)+\tilde{u}_{R}(x,y),
\end{align}
near the lower boundary $y=0$. Near the upper boundary $y=1$, we replace $\bar{u}_{R}(0)$ by $\bar{u}_{R}(1)$. By this decomposition, one has $\hat{u}_{R}(0)=0$ when \eqref{decom1} applied near $y=0$ and $\hat{u}_{R}(1)=0$ when \eqref{decom1} applied near $y=1$, and
$$
\int_{\mathbb{T}}\tilde{u}_{R}(x,y)\,{\rm d}x=0.
$$
For term involving $\tilde{u}_{R}$, we can apply the Poinc{a}r\'{e} inequality to $\tilde{u}_{R}$ and apply the Hardy inequality to $v_{R}$ to obtain
\begin{align*}
\int_{\mathbb{T}}\int_{0}^{1}\tilde{u}_{R}v_{R}\partial_{y}u^{a}\,{\rm d}y{\rm d}x\leq C\big(\|\partial_{x}u_{R}\|_{L^2}^2+\|\partial_{x}v_{R}\|_{L^2}^2\big).
\end{align*}
The estimates of tangential derivatives $\|(\partial_{x}u_{R},\partial_{x}v_{R})\|_{L^2}$ will be closed by the positivity estimates. For the term involving $\hat{u}_{R}(y)$, we can apply the Hardy inequality to obtain
\begin{align*}
\int_{\mathbb{T}}\int_{0}^{1}\hat{u}_{R}v_{R}\partial_{y}u_{B}^{a}\,{\rm d}y{\rm d}x\leq C\varepsilon^{4-2\alpha}\|\partial_{y}u_{R}\|_{L^2}^2+C\|\partial_{y}v_{R}\|_{L^2}^2=C\varepsilon^{4-2\alpha}\|\partial_{y}u_{R}\|_{L^2}^2+C\|\partial_{x}u_{R}\|_{L^2}^2.
\end{align*}
Here, we would like to point out that $4-2\alpha>2$ due to $0\leq \alpha<1$. The term involving the boundary-average terms $\bar{u}_{R}(0)$ and $\bar{u}_{R}(1)$ is bounded by
\begin{align*}
&C\varepsilon^{1-\alpha+\lambda}
\Big(\int_{\mathbb{T}}\int_{0}^{1}\frac{\bar{u}_{R}^2(0)}{y^{2\lambda}}\,{\rm d}y{\rm d}x\Big)^{\frac{1}{2}}\Big(\int_{\mathbb{T}}\int_{0}^{1}\frac{v_{R}^2}{y^2}\,{\rm d}y{\rm d}x\Big)^{\frac{1}{2}}\nonumber\\
&\quad +C\varepsilon^{1-\alpha+\lambda}\Big(\int_{\mathbb{T}}\int_{0}^{1}\frac{\bar{u}_{R}^2(1)}{|y-1|^{2\lambda}}\,{\rm d}y{\rm d}x\Big)^{\frac{1}{2}}\Big(\int_{\mathbb{T}}\int_{0}^{1}\frac{v_{R}^2}{(y-1)^2}\,{\rm d}y{\rm d}x\Big)^{\frac{1}{2}}\\
&\leq C\|\partial_{x}u_{R}\|_{L^2}^2+C\varepsilon^{2-2\alpha+2\lambda}\Big(\int_{\mathbb{T}}|u_{R}(x,0)|^2+|u_{R}(x,1)|^2\,{\rm d}x\Big).
\end{align*}
Here, we should require $2\lambda<1$ and $2-2\alpha+2\lambda>2-\alpha$, which leads to $\frac{\alpha}{2}<\lambda<\frac{1}{2}$.

At $\alpha=1$, however, the above arguments fail to work because no admissible choice of $\lambda$ is available. To overcome this difficulty, for $\alpha\geq 1$, we decompose $u_{R}$ as follows:
\begin{align}\label{decom2}
u_{R}=(u_{R}-\varepsilon^{\alpha} \partial_{y}u_{R})+\varepsilon^{\alpha} \partial_{y}u_{R}\,\, \text{near $y=0$},\quad u_{R}=(u_{R}+\varepsilon^{\alpha} \partial_{y}u_{R})-\varepsilon^{\alpha} \partial_{y}u_{R}\,\,\text{near $y=1$},
\end{align}
which allows us to apply the Hardy inequality to $u_{R}-\varepsilon^{\alpha}\partial_{y}u_R$ and $u_{R}+\varepsilon^{\alpha}\partial_{y}u_R$  by virtue of the Robin-type boundary condition. Specifically, the corresponding terms can be bounded by
\begin{align*}
    &C\varepsilon^2(\varepsilon^{\frac{1}{\tilde{\alpha}}}+\sigma)\|\partial_{y}(u_{R}\pm \varepsilon^{\alpha}\partial_{y}u_{R})\|_{L^2}^2+C\|\partial_{y}v_{R}\|_{L^2}^2\nonumber\\
    &\leq  C\varepsilon^2(\varepsilon^{\frac{1}{\tilde{\alpha}}}+\sigma)\|\partial_{y}u_{R}\|_{L^2}^2+C\varepsilon^{2+2\alpha}(\varepsilon^{\frac{1}{\tilde{\alpha}}}+\sigma)\|\partial_{yy}u_{R}\|_{L^2}^2+C\|\partial_{x}u_{R}\|_{L^2}^2.
\end{align*}
and the remaining terms are bounded as
\begin{align*}
    \int_{\mathbb{T}}\int_{0}^{1}\varepsilon^{\alpha}\partial_{y}u_{R}v_{R}\partial_{y}u_{B}^{a}\,{\rm d}y{\rm d}x\leq C\varepsilon^{2\alpha}(\varepsilon^{\frac{1}{\tilde{\alpha}}}+\sigma)\|\partial_{y}u_{R}\|_{L^2}^2+C\|\partial_{x}u_{R}\|_{L^2}^2.
\end{align*}
To close the estimate of $\|\partial_{yy}u_{R}\|_{L^2}^2$, we establish Stokes' estimates for the Navier--Stokes equations with a Robin-type boundary conditions. The decomposition \eqref{decom2} fails to work for the subcritical regime $0\leq \alpha<1$, since it will bring a term
\begin{align*}
\int_{\mathbb{T}}\int_{0}^{1}\varepsilon^{\alpha}\partial_{y}u_{R}v_{R}\partial_{y}u_{B}^{a}\,{\rm d}y{\rm d}x\leq C\varepsilon^2\|\partial_{y}u_{R}\|_{L^2}^2+C\|\partial_{x}u_{R}\|_{L^2}^2,
\end{align*}
where there is no smallness in front of $\varepsilon^2\|\partial_{y}u_{R}\|_{L^2}^2$ when $0\leq \alpha<1$.

The decompositions \eqref{decom1} and \eqref{decom2} are also used in the positivity estimates, where we test $\eqref{remainder}_{1}$ by $\partial_{x}u_{R}$ and $\eqref{remainder}_{2}$ by $\partial_{x}v_{R}$ to control $\|(\partial_{x}u_{R},\partial_{x}v_{R})\|_{L^2}^2$; see Section 6 for details.

\medskip

The rest of the paper is organized as follows. In Section 2, we construct the classical Prandtl expansions, including the scales induced by the Navier-friction boundary condition, and derive the equations for the interior and boundary-layer profiles at each order. The general well-posedness theory of the equation of interior and boundary-layer profiles will be presented in Sections 3 and 4, respectively. Then we apply these results to construct the approximate solutions and explain how compatibility conditions determine the free constants in the interior profiles in Section 5. In Section 6, we establish the existence and uniform estimates for remainders via the basic energy estimates, positivity estimates, Stokes' estimates and the Banach fixed point theorem. The proof of Theorems \ref{thm1}--\ref{thm3} will be presented at the end of Section 6. In Appendix A, we present a construction of the auxiliary correction terms used in the Prandtl expansions that enforce both incompressibility and the desired boundary conditions. A sketch of proof of the case where $\alpha$ is a general nonnegative rational number is presented in Appendix B.

\medskip

\textbf{Notation}: We denote $L^p(\Omega)$ ($1\leq p\leq \infty)$ by the standard Lebesgue space over $\Omega$ with the norm $\|\cdot\|_{L^{p}}$, and $H^{k}(\Omega)$ by the standard Sobolev space $W^{k,2}(\Omega)$ over $\Omega$ with the norm $\|\cdot\|_{H^{k}}$. Also, we denote $|\cdot|_{L^{p}(\mathbb{T})}$ and $|\cdot|_{H^{k}(\mathbb{T})}$ by the $L^{p}-$norm and $H^{k}-$norm over $\mathbb{T}$, respectively. We set $\langle k\rangle=(1+k^2)^{\frac{1}{2}}$, we denote by $C$ a generic constant independent of $\varepsilon$ and $\sigma$, whose value may change from line to line. 

\section{Prandtl Expansions of the Solution} 

It is well known from the Prandtl boundary layer theory that the boundary-layer thickness is of order $\varepsilon$,
which leads to a scale $\varepsilon^{-1}$ in the normal derivative $\partial_{y}$ near the boundaries. On the other hand, the slip length $\delta=\varepsilon^{\alpha}$ introduces the additional scale $\varepsilon^{\alpha-1}$ in the boundary conditions. So, to have clear expansions in powers of $\varepsilon$, we assume that $\alpha$ satisfies \eqref{alpha}, and let $\tilde{\alpha}$ be the one defined in \eqref{talpha}.
It is clear that $\tilde{\alpha}\geq 1$ and $\tilde{\alpha}\in \mathbb{Z}_{+}$. The case where $\alpha$ is a general nonnegative rational number is discussed in Appendix B.

\subsection{Interior expansions} We construct the interior expansions as follows
\begin{align}\label{1.5}
	u_{I}(x,y)\sim \sum\limits_{i=0}^{\infty}\var^{\frac{i}{\tilde{\alpha}}}u_{I}^{i}(x,y),\quad v_{I}(x,y)\sim \sum\limits_{i=0}^{\infty}\var^{\frac{i}{\tilde{\alpha}}}v_{I}^{i}(x,y),\quad p_{I}(x,y)\sim \sum\limits_{i=0}^{\infty}\var^{\frac{i}{\tilde{\alpha}}}p_{I}^{i}(x,y).
\end{align}
Then we substitute the expansions \eqref{1.5} into \eqref{1.1} and collect equal orders of $\var$.

\textit{The order of $\varepsilon^{0}$: the equations for $(u_{I}^{0},v_{I}^{0},p_{I}^{0})$}:
\begin{align}\label{1.4}
	\left\{
	\begin{aligned}
		&u_{I}^{0}\partial_{x}u_{I}^{0}+v_{I}^{0}\partial_{y}u_{I}^{0}+\partial_{x}p_{I}^{0}=0,\\
		&u_{I}^{0}\partial_{x}v_{I}^{0}+v_{I}^{0}\partial_{y}v_{I}^{0}+\partial_{y}p_{I}^{0}=0,\\
		&\partial_{x}u_{I}^{0}+\partial_{y}v_{I}^{0}=0.
	\end{aligned}
	\right.
\end{align}
For the incompressible Euler equation \eqref{1.4}, we impose following boundary conditions:
\begin{align}\label{1.6}
	v_{I}^{0}\vert_{y=0}=0,\quad v_{I}^{0}\vert_{y=1}=0.
\end{align}

\textit{The order of $\varepsilon^{\frac{k}{\tilde{\alpha}}}$: the equation for $(u_{I}^{k},v_{I}^{k},p_{I}^{k})$, $k\geq 1$}:
\begin{align}\label{1.9}
	\left\{
	\begin{aligned}
		&u_{I}^{0}\partial_{x}u_{I}^{k}+u_{I}^{k}\partial_{x}u_{I}^{0}+v_{I}^{0}\partial_{y}u_{I}^{k}+v_{I}^{k}\partial_{y}u_{I}^{0}+\partial_{x}p_{I}^{k}=F_{I}^{k},\\
		&u_{I}^{0}\partial_{x}v_{I}^{k}+u_{I}^{k}\partial_{x}v_{I}^{0}+v_{I}^{0}\partial_{y}v_{I}^{k}+v_{I}^{k}\partial_{y}v_{I}^{0}+\partial_{y}p_{I}^{k}=G_{I}^{k},\\
		&\partial_{x}u_{I}^{k}+\partial_{y}v_{I}^{k}=0.
	\end{aligned}
	\right.
\end{align}
Here we have denoted $(u_{I}^{i},v_{I}^{i})\equiv 0$ for $i<0$,
\begin{align}\label{1.9-1}
\begin{aligned}
&F_{I}^{k}=-\sum\limits_{j=1}^{k-1}\big[u_{I}^{j}\partial_{x}u_{I}^{k-j}+v_{I}^{j}\partial_{y}u_{I}^{k-j}\big]+\Delta u_{I}^{k-2\tilde{\alpha}},\\
&G_{I}^{k}=-\sum\limits_{j=1}^{k-1}\big[u_{I}^{j}\partial_{x}v_{I}^{k-j}+v_{I}^{j}\partial_{y}v_{I}^{k-j}\big]+\Delta v_{I}^{k-2\tilde{\alpha}}.
\end{aligned}
\end{align}
Here, the source terms $F_{I}^{k}$ and $G_{I}^{k}$ depend only on the lower-order interior profiles $(u_{I}^{j},v_{I}^{j})$ with $1\leq j\leq k-1$.

\subsection{Boundary layer expansions} To introduce the boundary layer expansions, we denote $z=\frac{y}{\varepsilon}$ and $\zeta=\frac{y-1}{\varepsilon}$ by the scaled coordinate at lower and upper boundaries, respectively.

\underline{\it Lower boundary layer expansions near $y=0$}:
\begin{align}\label{1.8-2}
	u_{B}\sim \sum\limits_{i=0}^{\infty}(\varepsilon^{\frac{i}{\tilde{\alpha}}})u_{B}^{i}(x,z),\quad v_{B}\sim\sum\limits_{i=0}^{\infty}(\varepsilon^{\frac{i}{\tilde{\alpha}}})v_{B}^{i}(x,z),\quad p_{B}\sim \sum\limits_{i=0}^{\infty}(\varepsilon^{\frac{i}{\tilde{\alpha}}})p_{B}^{i}(x,z).
\end{align}

\underline{\it Upper boundary layer expansions near $y=1$}:
\begin{align}\label{1.8-3}
	\hat{u}_{B}\sim \sum\limits_{i=0}^{\infty}(\varepsilon^{\frac{i}{\tilde{\alpha}}})\hat{u}_{B}^{i}(x,\zeta),\quad \hat{v}_{B}\sim\sum\limits_{i=0}^{\infty}(\varepsilon^{\frac{i}{\tilde{\alpha}}})\hat{v}_{B}^{i}(x,\zeta),\quad \hat{p}_{B}\sim \sum\limits_{i=0}^{\infty}(\varepsilon^{\frac{i}{\tilde{\alpha}}})\hat{p}_{B}^{i}(x,\zeta).
\end{align}
The equations for the upper boundary-layer profiles are analogous to those for the lower boundary layer, with $z$ replaced by $\zeta$ and the lower-boundary traces of the interior profiles replaced by their upper-boundary counterparts. We therefore derive only the lower boundary-layer equations for simplicity of presentation. Substituting $(u^{\varepsilon},v^{\varepsilon},p^{\varepsilon})\sim (u_{I}+u_{B},v_{I}+v_{B},p_{I}+p_{B})$ into \eqref{1.1}, and collecting equal orders of $\varepsilon$ gives:

{\it The order of $\varepsilon^{-1+\frac{k}{\tilde{\alpha}}}$: the equations of $(v_{B}^{k},p_{B}^{k})$ with $0\leq k\leq \tilde{\alpha}-1$}:
\begin{align}\label{1.11}
	\partial_{z}v_{B}^{k}=0,\quad \partial_{z}p_{B}^{k}=-\sum\limits_{i=0}^{k}\big[v_{I}^{k-i,-}+v_{B}^{k-i}\big]\partial_{z}v_{B}^{i},
\end{align}	
where we have denoted by $f^{-}=f(x,0)$ the trace of $f$ on the lower boundary $y=0$, and used the following Taylor expansions for interior solutions near the boundary $y=0$:	
\begin{align}
	&u_{I}^{k}(x,y)=u_{I}^{k}(x,0)+y\partial_{y}u_{I}^{k}(x,0)+\cdots+\frac{y^{\mathfrak{b}}\partial_{y}^{\mathfrak{b}}u_{I}^{k}(x,\xi_{u})}{\mathfrak{b}!},\label{1.11-1}\\
	&v_{I}^{k}(x,y)=v_{I}^{k}(x,0)+y\partial_{y}v_{I}^{k}(x,0)+\cdots+\frac{y^{\mathfrak{b}}\partial_{y}^{\mathfrak{b}}v_{I}^{k}(x,\xi_{v})}{\mathfrak{b}!},\label{1.11-2}
\end{align}
for some $\xi_{u},\xi_{v}\in (0,y)$. Here $\mathfrak{b}$ is a large enough integer which will be chosen in Section 6. Since the boundary-layer profiles should vanish away from the boundary $z=0$, i.e.,
$$
\lim\limits_{z\to \infty}v_{B}^{k}(x,z)=0,\quad \lim\limits_{z\to \infty}p_{B}^{k}(x,z)=0,
$$
we obtain from \eqref{1.11} that
\begin{align}\label{1.12}
	v_{B}^{k}\equiv 0,\quad p_{B}^{k}\equiv 0 \quad \text{for $0\leq k\leq \tilde{\alpha}-1$}.
\end{align}

{\it The order of $\varepsilon^{0}$: the equations of $(u_{B}^{0},v_{B}^{\tilde{\alpha}},p_{B}^{\tilde{\alpha}})$}:
\begin{align}\label{1.13}
	\left\{
	\begin{aligned}
		&(u_{I}^{0,-}+u_{B}^{0})\partial_{x}u_{B}^{0}+u_{B}^{0}\partial_{x}u_{I}^{0,-}+(v_{I}^{\tilde{\alpha},-}+v_{B}^{\tilde{\alpha}}+z\partial_{y}v_{I}^{0,-})\partial_{z}u_{B}^{0}=\partial_{zz}u_{B}^{0}\\
		&\partial_{z}p_{B}^{\tilde{\alpha}}=G_{B}^{\tilde{\alpha}},\\
		&\partial_{x}u_{B}^{0}+\partial_{z}v_{B}^{\tilde{\alpha}}=0,
	\end{aligned}
	\right.
\end{align}
where
\begin{align*}
	G_{B}^{\tilde{\alpha}}&=-(u_{I}^{0,-}+u_{B}^{0})\partial_{x}v_{B}^{0}-u_{B}^{0}\partial_{x}v_{I}^{0,-}-v_{B}^{0}\partial_{y}v_{I}^{0,-}+\partial_{zz}v_{B}^{0}\\
	&\quad-(v_{I}^{\tilde{\alpha},-}+v_{B}^{\tilde{\alpha}}+z\partial_{y}v_{I}^{0,-})\partial_{z}v_{B}^{0}-\sum\limits_{i=1}^{\tilde{\alpha}}(v_{I}^{\tilde{\alpha}-i,-}+v_{B}^{\tilde{\alpha}-i})\partial_{z}v_{B}^{i}.
\end{align*}
Using the boundary conditions \eqref{1.6} for $v_{I}^{0}$ and \eqref{1.12}, one has $G_{B}^{\tilde{\alpha}}\equiv 0$, which yields
\begin{align}\label{1.14}
	p_{B}^{\tilde{\alpha}}\equiv 0.
\end{align}

{\it The order of $\varepsilon^{\frac{k}{\tilde{\alpha}}}$: the equations of $(u_{B}^{k},v_{B}^{\tilde{\alpha}+k},p_{B}^{\tilde{\alpha}+k})$ with $k\geq 1$}:
\begin{align}\label{1.15}
	\left\{
	\begin{aligned}
		&(u_{I}^{0,-}+u_{B}^{0})\partial_{x}u_{B}^{k}+u_{B}^{k}\partial_{x}(u_{I}^{0,-}+u_{B}^{0})+(v_{I}^{\tilde{\alpha}+k,-}+v_{B}^{\tilde{\alpha}+k})\partial_{z}u_{B}^{0}\\
		&\quad +(v_{I}^{\tilde{\alpha},-}+v_{B}^{\tilde{\alpha}}+z\partial_{y}v_{I}^{0,-})\partial_{z}u_{B}^{k}=\partial_{zz}u_{B}^{k}+F_{B}^{k}\\
		&\partial_{z}p_{B}^{\tilde{\alpha}+k}=G_{B}^{\tilde{\alpha}+k},\\
		&\partial_{x}u_{B}^{k}+\partial_{z}v_{B}^{\tilde{\alpha}+k}=0,
	\end{aligned}
	\right.
\end{align}
where
\begin{align}
	F_{B}^{k}&=-\sum\limits_{i=1}^{k-1}u_{B}^{i}\partial_{x}u_{B}^{k-i}-\sum\limits_{i=1}^{k}\Big(\sum\limits_{0\leq j\leq \frac{i}{\tilde{\alpha}}}\frac{z^{j}\partial_{y}^{j}u_{I}^{i-\tilde{\alpha}j,-}}{j!}\Big)\partial_{x}u_{B}^{k-i}-\sum\limits_{i=1}^{k}u_{B}^{k-i}\sum\limits_{0\leq j\leq \frac{i}{\tilde{\alpha}}}\frac{z^{j}\partial_{x}\partial_{y}^{j}u_{I}^{i-\tilde{\alpha}j,-}}{j!}
	\nonumber\\
	&\quad 
	-\sum\limits_{i=0}^{k}v_{B}^{k-i}\sum\limits_{0\leq j\leq \frac{i}{\tilde{\alpha}}	}\frac{z^j\partial_{y}^{j+1}u_{I}^{i-\tilde{\alpha}j,-}}{j!}-\sum\limits_{i=1}^{k-1}\Big(v_{B}^{\tilde{\alpha}+i}+\sum\limits_{0\leq j\leq 1+\frac{i}{\tilde{\alpha}}	}\frac{z^{j}\partial_{y}^{j}v_{I}^{i-(j-1)\tilde{\alpha},-}}{j!}\Big)\partial_{z}u_{B}^{k-i}\nonumber\\
	&\quad-\sum\limits_{1\leq j\leq 1+\frac{k}{\tilde{\alpha}}}\frac{z^j\partial_{y}^{j}v_{I}^{k-(j-1)\tilde{\alpha},-}}{j!}\partial_{z}u_{B}^{0} -\partial_{x}p_{B}^{k}+\partial_{xx}u_{B}^{k-2\tilde{\alpha}},\label{1.15-4}\\
	G_{B}^{\tilde{\alpha}+k}&=-\sum\limits_{i=0}^{k}\Big[\Big(u_{B}^{i}+\sum\limits_{0\leq j\leq \frac{i}{\tilde{\alpha}}	}\frac{z^j\partial_{y}^{j}u_{I}^{i-\tilde{\alpha}j,-}}{j!}\Big)\partial_{x}v_{B}^{k-i}+u_{B}^{k-i}\sum\limits_{0\leq j\leq \frac{i}{\tilde{\alpha}}}\frac{z^j\partial_{x}\partial_{y}^{j}v_{I}^{i-\tilde{\alpha}j,-}}{j!}\Big]\nonumber\\
	&\quad -\sum\limits_{i=0}^{\tilde{\alpha}+k}\Big(v_{B}^{i}+\sum\limits_{0\leq j\leq \frac{i}{\tilde{\alpha}}}\frac{z^j\partial_{y}^{j}v_{I}^{i-\tilde{\alpha}j,-}}{j!}\Big)\partial_{z}v_{B}^{\tilde{\alpha}+k-i}-\sum\limits_{i=0}^{k}v_{B}^{k-i}\sum\limits_{0\leq j\leq \frac{i}{\tilde{\alpha}}}\frac{z^j\partial_{y}^{j+1}v_{I}^{i-\tilde{\alpha}j,-}}{j!}\nonumber\\
	&\quad +\partial_{xx}v_{B}^{k-2\tilde{\alpha}}+\partial_{zz}v_{B}^{k}.\label{1.15-5}
\end{align}

\subsection{Boundary conditions} Substituting 
$(u^{\var},v^{\var})\sim (u_{I}+\chi(y)u_{B}+\tilde{\chi}(y)\hat{u}_{B},v_{I}+\chi(y)v_{B}+\tilde{\chi}(y)\hat{v}_{B})$ into \eqref{1.3} and collecting equal orders of $\varepsilon$, we have

\textbf{The boundary conditions for interior profiles}:
\begin{align}\label{1.10}
	v_{I}^{k}\vert_{y=0}=-v_{B}^{k}\vert_{z=0},\quad v_{I}^{k}\vert_{y=1}=-\hat{v}_{B}^{k}\vert_{\zeta=0},\quad \text{for $k\geq 1$}.
\end{align}

\textbf{The boundary conditions for boundary-layer profiles}:
\begin{itemize}
\item[(C1)] for $\alpha=0$,
\begin{align}\label{b00}
    \begin{aligned}
&\partial_{z}u_{B}^{0}(x,0)=\partial_{\zeta}\hat{u}_{B}^{0}(x,0)=0,\\
&\partial_{z}u_{B}^{1}(x,0)=-f_{0}(x)+u_{I}^{0}(x,0)+u_{B}^{0}(x,0)-\partial_{y}u_{I}^{0}(x,0),\\
&\partial_{\zeta}\hat{u}_{B}^{(1-\alpha)\tilde{\alpha}}(x,0)=f_{1}(x)-u_{I}^{0}(x,1)-\hat{u}_{B}^{0}(x,0)-\partial_{y}u_{I}^{0}(x,1),\\
&
\begin{aligned}
&\partial_{z}u_{B}^{(1-\alpha)\tilde{\alpha}+k}(x,0)=u_{I}^{k}(x,0)+u_{B}^{k}(x,0)-\partial_{y}u_{I}^{k}(x,0),\\ 
&\partial_{\zeta}\hat{u}_{B}^{(1-\alpha)\tilde{\alpha}+k}(x,0)=-u_{I}^{k}(x,1)-\hat{u}_{B}^{k}(x,0)-\partial_{y}u_{I}^{k}(x,1),
\end{aligned}\quad \text{for $k\geq 1$.}
\end{aligned}
\end{align}
    \item[(C2)] for $0<\alpha<1$,
    \begin{align}\label{b}
    \begin{aligned}
&\partial_{z}u_{B}^{i}(x,0)=\partial_{\zeta}\hat{u}_{B}^{i}(x,0)=0,\quad \text{for $0\leq i\leq (1-\alpha)\tilde{\alpha}-1$}\\
&\partial_{z}u_{B}^{(1-\alpha)\tilde{\alpha}}(x,0)=-f_{0}(x)+u_{I}^{0}(x,0)+u_{B}^{0}(x,0),\\
&\partial_{\zeta}\hat{u}_{B}^{(1-\alpha)\tilde{\alpha}}(x,0)=f_{1}(x)-u_{I}^{0}(x,1)-\hat{u}_{B}^{0}(x,0),\\
&
\begin{aligned}
&\partial_{z}u_{B}^{(1-\alpha)\tilde{\alpha}+k}(x,0)=u_{I}^{k}(x,0)+u_{B}^{k}(x,0)-\partial_{y}u_{I}^{k-\alpha\tilde{\alpha}}(x,0),\\ 
&\partial_{\zeta}\hat{u}_{B}^{(1-\alpha)\tilde{\alpha}+k}(x,0)=-u_{I}^{k}(x,1)-\hat{u}_{B}^{k}(x,0)-\partial_{y}u_{I}^{k-\alpha\tilde{\alpha}}(x,1),
\end{aligned}\quad \text{for $k\geq 1$.}
\end{aligned}
\end{align}
\item[(C3)] for $\alpha=1$,
\begin{align}\label{b2}
    \begin{aligned}
&u_{B}^{0}(x,0)-\partial_{z}u_{B}^{0}(x,0)=f_{0}(x)-u_{I}^{0}(x,0),\quad \hat{u}_{B}^{0}(x,0)+\partial_{\zeta}\hat{u}_{B}^{0}(x,0)=f_{1}(x)-u_{I}^{0}(x,1),\\
&\begin{aligned}
&u_{B}^{k}(x,0)-\partial_{z}u_{B}^{k}(x,0)=-u_{I}^{k}(x,0)+\partial_{y}u_{I}^{k-1}(x,0),\\ &\hat{u}_{B}^{k}(x,0)+\partial_{\zeta}\hat{u}_{B}^{k}(x,0)=-u_{I}^{k}(x,1)-\partial_{y}u_{I}^{k-1}(x,1),
\end{aligned}\quad \text{for $k\geq 1$.}
\end{aligned}
\end{align}
\item[(C4)] for $\alpha>1$,
\begin{align}\label{b3}
    \begin{aligned}
&u_{B}^{0}(x,0)=f_{0}(x)-u_{I}^{0}(x,0),\quad \hat{u}_{B}^{0}(x,0)=f_{1}(x)-u_{I}^{0}(x,1),\\
&\begin{aligned}
&u_{B}^{k}(x,0)=-u_{I}^{k}(x,0)+\partial_{z}u_{B}^{k-(\alpha-1)\tilde{\alpha}}(x,0)+\partial_{y}u_{I}^{k-\alpha\tilde{\alpha}}(x,0),\\
&\hat{u}_{B}^{k}(x,0)=-u_{I}^{k}(x,1)-\partial_{\zeta}\hat{u}_{B}^{k-(\alpha-1)\tilde{\alpha}}(x,0)-\partial_{y}u_{I}^{k-\alpha\tilde{\alpha}}(x,1),
\end{aligned}\quad \text{for $k\geq 1$.}
\end{aligned}
\end{align}
\end{itemize}
Here, we have used $u_{B}^{i}=\hat{u}_{B}^{i}\equiv 0$ for $i<0$. It is worth pointing out that for the case $0\leq \alpha<1$, the leading boundary-layer profiles $(u_{B}^{0},v_{B}^{\tilde{\alpha}})$ and $(\hat{u}_{B}^{0},\hat{v}_{B}^{\tilde{\alpha}})$ satisfy the Prandtl equations with homogeneous boundary conditions. It will be shown in Lemma \ref{lem2} that $(u_{B}^{0},v_{B}^{\tilde{\alpha}})=(\hat{u}_{B}^{0},\hat{v}_{B}^{\tilde{\alpha}})\equiv (0,0)$ from, so the leading boundary layer profiles are vanishing. This is why $0\leq \alpha<1$ is called subcritical case. When $\alpha=1$, the leading boundary-layer profile $(u_{B}^{0},v_{B}^{\tilde{\alpha}})$ and $(\hat{u}_{B}^{0},\hat{v}_{B}^{\tilde{\alpha}})$ satisfy the Prandtl equations with a Robin boundary condition, whereas for $\alpha>1$, it satisfies a Dirichlet boundary condition. Therefore, these regimes are referred to as the critical and supercritical cases, respectively.

\section{Well-Posedness of the Interior Solutions} In this section, we establish the well-posedness theory of the interior solutions. Recall \eqref{1.4} and \eqref{1.9}, for every $\alpha\in \R$, the leading interior profile $(u_{I}^{0},v_{I}^{0},p_{I}^{0})$ satisfies the incompressible Euler equations, the higher-order interior profile $(u_{I}^{k},v_{I}^{k},p_{I}^{k})$ satisfies the linearized Euler equation around $(u_{I}^{0},v_{I}^{0},p_{I}^{0})$.

\subsection{Well-posedness of zeroth-order interior profiles}

According to the Prandtl--Batchelor theory (cf. \cite{B-1956,FGLT-2023,FGLT-2024,GXX,PZ-2025}), the leading Euler flow $(u_{I}^{0},v_{I}^{0},p_{I}^{0})$ must be the Couette flow:
\begin{align}\label{2.0}
	u_{I}^{0}=Ay+B,\quad v_{I}^{0}=0,\quad p_{I}^{0}=0
\end{align}
for some constants $A$ and $B$. However, the values of $A$ and $B$ are determined by the  solvability conditions for the corresponding boundary-layer problems, that is, $(u_{B}^{(1-\alpha)\tilde{\alpha}},\hat{u}_{B}^{(1-\alpha)\tilde{\alpha}})$ for subcritical case $0\leq \alpha<1$, and $(u_{B}^{0},\hat{u}_{B}^{0})$ for critical and supercritical cases $\alpha\geq 1$. Such choices are given in the Proposition \ref{prop3}.

\subsection{Well-posedness of higher-order interior profiles}

Recalling \eqref{1.9}, \eqref{1.9-1} and \eqref{1.10}, the higher-order interior profiles $(u_{I}^{k},v_{I}^{k},p_{I}^{k})$ satisfy the following linearized incompressible Euler equations with external forces around Couette flow $(u_{I}^{0},v_{I}^{0})=(Ay+B,0)$:
\begin{align}\label{2.1}
	\left\{
	\begin{aligned}
		&u_{I}^{0}\partial_{x}u_{I}^{k}+v_{I}^{k}\partial_{y}u_{I}^{0}+\partial_{x}p_{I}^{k}=F_{I}^{k},\\
		&u_{I}^{0}\partial_{x}v_{I}^{k}+\partial_{y}p_{I}^{k}=G_{I}^{k},\\
		&\partial_{x}u_{I}^{k}+\partial_{y}v_{I}^{k}=0,\\
		&v_{I}^{k}\vert_{y=0}=-v_{B}^{k}\vert_{z=0},\quad v_{I}^{k}\vert_{y=1}=-\hat{v}_{B}^{k}\vert_{\zeta=0},
	\end{aligned}
	\right.
\end{align}
where $F_{I}^{k}$ and $G_{I}^{k}$ are the ones given in \eqref{1.9-1}.

\begin{proposition}\label{interior}
  For any $k\geq 1$, there exists a solution $(u_{I}^{k},v_{I}^{k},p_{I}^{k})$ to \eqref{2.1} satisfying
    \begin{align}\label{harmonic}
    \Delta u_{I}^{k}=\Delta v_{I}^{k}=0,
    \end{align}
    and
    \begin{align}\label{v}
    \|v_{I}^{k}\|_{H^m}\leq C_{m}\big(|v_{B}^{k}(\cdot,0)|_{H^{m-\frac{1}{2}}(\mathbb{T})}+|\hat{v}_{B}^{k}(\cdot,0)|_{H^{m-\frac{1}{2}}(\mathbb{T})}\big)\quad \text{for any $m\geq 1$},
    \end{align}
    where $C_{m}$ is a positive constant depending only on $m$.
    Furthermore, $u_{I}^{k}$ can be chosen as
    \begin{align}\label{u}
        u_{I}^{k}=-\int_{0}^{x}\partial_{y}v_{I}^{k}(x',y)\,{\rm d}x'+
	\int_{0}^{y}\partial_{x}v_{I}^{k}(0,y')\,{\rm d}y'+\mathfrak{C}_{k}y+\tilde{\mathfrak{C}}_{k},
    \end{align}
    where $\tilde{\mathfrak{C}}_{k}$ and $\mathfrak{C}_{k}$ are two constants determined uniquely by the compatibility conditions for the corresponding boundary-layer equations. And $p_{I}^{k}$ is constructed by
    \begin{align}\label{p}
    p_{I}^{k}(x,y)&=\big[Ay+B]\int_{0}^{x}\partial_{y}v_{I}^{k}(x',y)\,{\rm d}x'-A\int_{0}^{x}v_{I}^{k}(x',y)\,{\rm d}x'+\int_{0}^{x}F_{I}^{k}(x',y)\,{\rm d}x'\nonumber\\
    &\quad +\int_{0}^{y}\Big[-(Ay'+B)\partial_{x}v_{I}^{k}(0,y')-\int_{0}^{x}\partial_{y}F_{I}^{k}(x',y')\,{\rm d}x'+G_{I}^{k}(x,y')\Big]\,{\rm d}y'.
    \end{align} 
\end{proposition}

\noindent\textbf{Proof}. We apply $\partial_{y}$ to $\eqref{2.1}_{1}$ and $\partial_{x}$ to $\eqref{2.1}_{2}$ to eliminate the pressure $p_{I}^{k}$:
\begin{align}\label{2.3}
	-u_{I}^{0}\Delta v_{I}^{k}&=\partial_{y}F_{I}^{k}-\partial_{x}G_{I}^{k}=\sum\limits_{j=1}^{k-1}\big(u_{I}^{j}\Delta v_{I}^{k-j}-v_{I}^{j}\Delta u_{I}^{k-j}\big)+(\partial_{y}\Delta u_{I}^{k-2\tilde{\alpha}}-\partial_{x}\Delta v_{I}^{k-2\tilde{\alpha}}).
\end{align}

\textbf{Step 1.} For $k=1$, it follows from \eqref{1.9-1} and \eqref{2.0} that $F_{I}^{1}=G_{I}^{1}\equiv 0$. Then \eqref{2.3} leads to
\begin{align}\label{2.4}
	\Delta v_{I}^{1}=0,\quad v_{I}^{1}\vert_{y=0}=-v_{B}^{1}\vert_{z=0},\quad v_{I}^{1}\vert_{y=1}=-\hat{v}_{B}^{1}\vert_{\zeta=0}.
\end{align}
It follows from the elliptic theory that there exists a unique classical solution to \eqref{2.4} satisfying \eqref{v} for $k=1$. Moreover, noting from the incompressible condition that
$$
\int_{\mathbb{T}}v_{B}^{1}(x,z)\,{\rm d}x=\int_{\mathbb{T}}\hat{v}_{B}^{1}(x,\zeta)\,{\rm d}x=0,
$$
one obtains from $\eqref{2.1}_{3}$ and \eqref{2.4} that
\begin{align}\label{vzeromean}
\int_{\mathbb{T}}\partial_{y}v_{I}^{1}(x,y)\,{\rm d}x=\int_{\mathbb{T}}v_{I}^{1}(x,y)\,{\rm d}x=0.
\end{align}

To solve $u_{I}^{1}$, we observe that it only needs to satisfy the incompressible condition $\eqref{2.1}_{3}$. Actually, once $v_{I}^{1}$ is known, the choice of $u_{I}^{1}$ doesn't affect the solvability of the pressure $p_{I}^{1}$ because $u_{I}^{1}$ enters the equations \eqref{2.1} through $\partial_{x}u_{I}^{1}$, which can be replaced by $-\partial_{y}v_{I}^{1}$ via the incompressible condition. We therefore define
\begin{align}\label{2.11}
	u_{I}^{1}=-\int_{0}^{x}\partial_{y}v_{I}^{1}(x',y)\,{\rm d}x'+\int_{0}^{y}\partial_{x}v_{I}^{1}(0,y')\,{\rm d}y'+\mathfrak{C}_{1}y+\tilde{\mathfrak{C}}_{1}.
\end{align}
where $\mathfrak{C}_{1}$ and $\tilde{\mathfrak{C}}_{1}$ are determined uniquely by the compatibility conditions of the corresponding boundary-layer profiles. Since  $\Delta v_{I}^{1}=0$,
\begin{align*}
	\Delta u_{I}^{1}&=-\partial_{xy}v_{I}^{1}(x,y)-\int_{0}^{x}\partial_{yyy}v_{I}^{1}(x',y)\,{\rm d}x'+\partial_{xy}v_{I}^{1}(0,y)\\
	&=-\partial_{xy}v_{I}^{1}(x,y)+\int_{0}^{x}\partial_{xxy}v_{I}^{1}(x',y)\,{\rm d}x'+\partial_{xy}v_{I}^{1}(0,y)=0.
\end{align*}
Therefore, $u_{I}^{1}$ is also harmonic by our construction. Furthermore, it follows from \eqref{vzeromean} and \eqref{2.11} that $u_{I}^{1}$ is periodic with respect to $x$.

Now, we construct $p_{I}^{1}$. It follows from $\eqref{2.1}_{1}$ that
\begin{align}\label{2.8}
	p_{I}^{1}=[Ay+B]\int_{0}^{x}\partial_{y}v_{I}^{1}(x',y)\,{\rm d}x'-A\int_{0}^{x}v_{I}^{1}(x',y)\,{\rm d}x'+\phi_{1}(y),
\end{align}
where we have used $F_{I}^{1}\equiv 0$, and $\phi_{1}(y)$ is a function depending only on $y$. Substituting \eqref{2.8} into $\eqref{2.1}_{2}$ gives
\begin{align}\label{2.9}
	\phi_{1}'(y)&=-(Ay+B)\partial_{x}v_{I}^{1}-(Ay+B)\int_{0}^{x}\partial_{yy}v_{I}^{1}(x',y)\,{\rm d}x'\nonumber\\
	&\quad -A\int_{0}^{x}\partial_{y}v_{I}^{1}(x',y)\,{\rm d}x'+A\int_{0}^{x}\partial_{y}v_{I}^{1}(x',y)\,{\rm d}x'\nonumber\\
	&=-(Ay+B)\partial_{x}v_{I}^{1}(0,y),
\end{align}
where we have used $G_{I}^{1}\equiv 0$ and $\partial_{yy}v_{I}^{1}=-\partial_{xx}v_{I}^{1}$ due to $\Delta v_{I}^{1}=0$. Combining \eqref{2.8}--\eqref{2.9} together, we obtain \eqref{p} for $k=1$. Moreover, it follows from \eqref{vzeromean} and \eqref{2.8} that $p_{I}^{1}$ is periodic with respect to $x$.

\textbf{Step 2}. Now, we assume
\begin{align}\label{H1}
\Delta v_{I}^{j}=\Delta u_{I}^{j}=0\quad \text{for $1\leq j\leq k-1$}.
\end{align}
Then, we obtain from \eqref{2.3} that
\begin{align*}
	-u_{I}^{0}\Delta v_{I}^{k}=\partial_{y}F_{I}^{k}-\partial_{x}G_{I}^{k}\equiv 0,
\end{align*}
which yields that $v_{I}^{k}$ satisfies
\begin{align}\label{2.13}
	\Delta v_{I}^{k}=0,\quad v_{I}^{k}\vert_{y=0}=-v_{B}^{k}\vert_{z=0},\quad v_{I}^{k}\vert_{y=1}=-\hat{v}_{B}^{k}\vert_{\zeta=0}.
\end{align}
It follows from the elliptic theory that there exists a unique classical solution $v_{I}^{k}$ to \eqref{2.13} satisfying \eqref{v}. Moreover, noting from the incompressible condition that 
$$
\int_{\mathbb{T}}v_{B}^{k}(x,z)\,{\rm d}x=\int_{\mathbb{T}}\hat{v}_{B}^{k}(x,\zeta)\,{\rm d}x=0,
$$
we obtain from $\eqref{2.1}_{3}$ and \eqref{2.13} that
\begin{align}\label{vzeromean2}
\int_{\mathbb{T}}\partial_{y}v_{I}^{k}(x,y)\,{\rm d}x=\int_{\mathbb{T}}v_{I}^{k}(x,y)\,{\rm d}x=0.
\end{align}

For $u_{I}^{k}$, similar to the construction of $u_{I}^{1}$, we only need to construct $u_{I}^{k}$ such that it satisfies the incompressible condition and the harmonic equation. Therefore, we may construct $u_{I}^{k}$ as follows
\begin{align}\label{2.15}
	u_{I}^{k}=-\int_{0}^{x}\partial_{y}v_{I}^{k}(x',y)\,{\rm d}x'+
	\int_{0}^{y}\partial_{x}v_{I}^{k}(0,y')\,{\rm d}y'+\mathfrak{C}_{k}y+\tilde{\mathfrak{C}}_{k},
\end{align}
where $\mathfrak{C}_{k}$ and $\tilde{\mathfrak{C}}_{k}$ are uniquely determined by the compatibility conditions for the corresponding boundary-layer profiles. A direct calculation shows that
\begin{align}\label{2.15-1}
	\Delta u_{I}^{k}&=-\partial_{xy}v_{I}^{k}(x,y)-\int_{0}^{x}\partial_{yyy}v_{I}^{k}(x',y)\,{\rm d}x'+\partial_{xy}v_{I}^{k}(0,y)\nonumber\\
	&=-\partial_{xy}v_{I}^{k}(x,y)+\int_{0}^{x}\partial_{xxy}v_{I}^{k}(x',y)\,{\rm d}x'+\partial_{xy}v_{I}^{k}(0,y)=0,
\end{align}
where we have used $\Delta v_{I}^{k}=0$. Combining \eqref{2.13} and \eqref{2.15-1} closes the induction assumption \eqref{H1}. Therefore \eqref{harmonic} and \eqref{u} hold for all $k\geq 1$. Moreover, it follows from \eqref{vzeromean2} and \eqref{2.15} that $u_{I}^{k}$ is periodic with respect to $x$.

As for $p_{I}^{k}$, we first obtain from $\eqref{2.1}_{1}$ that
\begin{align}\label{2.16}
	p_{I}^{k}=\big(Ay+B)\int_{0}^{x}\partial_{y}v_{I}^{k}(x',y)\,{\rm d}x'-A\int_{0}^{x}v_{I}^{k}(x',y)\,{\rm d}x'+\int_{0}^{x}F_{I}^{k}(x',y)\,{\rm d}x'+\phi_{k}(y),
\end{align}
where $\phi_{k}(y)$ is a function depending only on $y$. To determine $\phi_{k}(y)$, we substitute \eqref{2.16} into $\eqref{2.1}_{2}$ and obtain
\begin{align}\label{2.17}
	\phi_{k}'(y)&=-(Ay+B)\partial_{x}v_{I}^{k}-(Ay+B)\int_{0}^{x}\partial_{yy}v_{I}^{k}(x',y)\,{\rm d}x'\nonumber\\
	&\quad -A\int_{0}^{x}\partial_{y}v_{I}^{k}(x',y)\,{\rm d}x'+A\int_{0}^{x}\partial_{y}v_{I}^{k}(x',y)\,{\rm d}x'-\int_{0}^{x}\partial_{y}F_{I}^{k}(x',y)\,{\rm d}x'+G_{I}^{k}\nonumber\\
	&=-(Ay+B)\partial_{x}v_{I}^{k}(0,y)-\int_{0}^{x}\partial_{y}F_{I}^{k}(x',y)\,{\rm d}x'+G_{I}^{k}.
\end{align}
It remains to show $-\int_{0}^{x}\partial_{y}F_{I}^{k}(x',y)\,{\rm d}x'+G_{I}^{k}$ depends only on $y$. In fact, it follows from a straightforward computation and \eqref{2.3} that
\begin{align*}
	\partial_{x}\Big(-\int_{0}^{x}\partial_{y}F_{I}^{k}(x',y)\,{\rm d}x'+G_{I}^{k}\Big)&=-\partial_{y}F_{I}^{k}+\partial_{x}G_{I}^{k}\nonumber\\
	&=-\sum\limits_{j=1}^{k-1}\big(u_{I}^{j}\Delta v_{I}^{k-j}-v_{I}^{j}\Delta u_{I}^{k-j}\big)-(\partial_{y}\Delta u_{I}^{k-2\tilde{\alpha}}-\partial_{x}\Delta v_{I}^{k-2\tilde{\alpha}})\nonumber\\
	&=0,
\end{align*} 
where we used \eqref{harmonic}. This indicates that $-\int_{0}^{x}\partial_{y}F_{I}^{k}(x',y)\,{\rm d}x'+G_{I}^{k}$ depends only on $y$, and \eqref{2.17} yields
\begin{align*}
\phi_{k}(y)=-\int_{0}^{y}(As+B)\partial_{x}v_{I}^{k}(0,s)\,{\rm d}s+\int_{0}^{y}\Big(-\int_{0}^{x}\partial_{y}F_{I}^{k}(x',s)\,{\rm d}x'+G_{I}^{k}(x,s)\Big)\,{\rm d}s,
\end{align*}
which completes the proof of \eqref{p}. 

Noting $\eqref{1.9-1}_{1}$ and using \eqref{harmonic}, a direct calculation shows that
\begin{align}\label{FIK}
    F_{I}^{k}=-\frac{1}{2}\sum\limits_{j=1}^{k-1}\partial_{x}\big(u_{I}^{j}u_{I}^{k-j}+v_{I}^{j}v_{I}^{k-j}\big)-\sum\limits_{j=1}^{k-1}\big(\partial_{y}u_{I}^{k-j}-\partial_{x}v_{I}^{k-j}\big)v_{I}^{j}.
\end{align}
Denoting $\omega^{k-j}:=\partial_{y}u_{I}^{k-j}-\partial_{x}v_{I}^{k-j}$ the $(k-j)$th order vorticity, one has from \eqref{harmonic} that $\partial_{x}\omega^{k-j}=-\Delta v_{I}^{k-j}=0$ and $\partial_{y}\omega^{k-j}=\Delta u_{I}^{k-j}=0$. So, $\omega^{k-j}$ is a constant. Thus, using \eqref{vzeromean2}, we obtain from \eqref{FIK} that
\begin{align}\label{Fzeromean}
    \int_{\mathbb{T}}F_{I}^{k}(x,y)\,{\rm d}x=0.
\end{align}
It follows from \eqref{p}, \eqref{vzeromean2}, \eqref{Fzeromean} and the periodicity of $F_{I}^{k}$ that $p_{I}^{k}$ is periodic with respect to $x$. Therefore, the proof of Proposition \ref{interior} is complete. $\hfill\square$

\section{Well-Posedness of the Boundary-Layer Expansions} 
In this section, we shall focus on the well-posedness of boundary-layer equations.

\subsection{Well-posedness of zeroth-order boundary-layer profiles} Recall \eqref{1.13} and \eqref{b00}--\eqref{b3}, $(u_{B}^{0}+u_{I}^{0}(x,0),v_{B}^{\tilde{\alpha}}+v_{I}^{\tilde{\alpha}}(x,0))$ and $(\hat{u}_{B}^{0}+u_{I}^{0}(x,1),\hat{v}_{B}^{\tilde{\alpha}}+v_{I}^{\tilde{\alpha}}(x,1))$ will satisfy the Prandtl equation with homogeneous boundary conditions, Robin boundary conditions and Dirichlet boundary conditions for $0\leq \alpha<1$, $\alpha=1$ and $\alpha>1$, respectively. We first deal with the case $\alpha\geq 1$ and then case $0\leq \alpha<1$.

Consider the Prandtl equation:
\begin{align}\label{p0}
\left\{
\begin{aligned}
&u\partial_{x}u+v\partial_{z}u=\partial_{zz}u,\\
&\partial_{x}u+\partial_{z}v=0,
\end{aligned}
\right.
\end{align}
supplemented with the far field condition
\begin{align}\label{f0}
    \lim\limits_{z\to\infty}u=\mathfrak{B},
\end{align}
and the Dirichlet boundary conditions:
\begin{align}\label{D1}
v(x,0)=0,\quad u(x,0)=\mathfrak{h}_{1}(x),
\end{align}
or the Robin-type boundary conditions
\begin{align}\label{R1}
v(x,0)=0,\quad (u-\partial_{z}u)(x,0)=\mathfrak{h}_{2}(x),
\end{align}
Here, $\bar{\mathfrak{h}}_{1}>0$ and $\bar{\mathfrak{h}}_{2}>0$.

\subsubsection{Supercritical: $\alpha>1$.} For the well-posedness of Prandtl equation \eqref{p0} with far-field condition \eqref{f0} and the boundary condition \eqref{D1}, we use the von Mises transformation:
\begin{align*}
\xi=x,\quad \eta=\int_{0}^{z}u(x,s)\,{\rm d}s,
\end{align*}
and denote
\begin{align}\label{w0}
w=u^2(\xi,\eta).
\end{align}
Then a direct calculation shows that
\begin{align*}
\partial_{\xi}w=2u\partial_{x}u+2u\partial_{z}u\big(\frac{v}{u}\big)=2(u\partial_{x}u+v\partial_{z}u),
\end{align*}
and
\begin{align*}
\partial_{\eta}w=2\partial_{z}u,\quad \partial_{\eta\eta}w=\frac{2}{u}\partial_{zz}u,
\end{align*}
which implies that $w$ satisfies
\begin{align}\label{w}
\left\{
\begin{aligned}
&\frac{1}{\sqrt{w}}\partial_{\xi}w=\partial_{\eta\eta}w,\\
&w(x,0)=\mathfrak{h}_{1}^2,\quad \lim\limits_{\eta\to \infty}w=\mathfrak{B}^2.
\end{aligned}
\right.
\end{align}
Then, integrating \eqref{w} over $\mathbb{T}$ gives $\bar{w}''=0$. Using the far-field condition $\eqref{w}_{2}$, we obtain $\bar{w}=\mathfrak{B}^2$. In particular, it holds that
\begin{align*}
\int_{\mathbb{T}}w(\xi,0)\,{\rm d}\xi=\lim\limits_{\eta\to \infty}\int_{\mathbb{T}}w(\xi,\eta)\,{\rm d}\xi,
\end{align*}
which yields
\begin{align}\label{bw}
\mathfrak{B}^2=\frac{1}{2\pi}\int_{\mathbb{T}}\mathfrak{h}_{1}^2(\xi)\,{\rm d}\xi.
\end{align}
The formula \eqref{bw} is well known as the Batchelor--Wood formula \cite{B-1956,Wood}, and is a compatibility condition for the smooth solution of \eqref{p0} with \eqref{f0}--\eqref{D1}.

Under the compatibility condition \eqref{bw}, we further assume $\mathfrak{h}_{1}$ is a small perturbation around its average over $\mathbb{T}$, i.e., 
\begin{align}\label{h1}
\mathfrak{h}_{1}(\xi)=\frac{1}{2\pi}\int_{\mathbb{T}}\mathfrak{h}_{1}(\xi)\,{\rm d}\xi+\sigma\tilde{\mathfrak{h}}_{1}(\xi)=:\bar{\mathfrak{h}}_{1}+\sigma \tilde{\mathfrak{h}}_{1}(\xi),
\end{align}
where $\int_{\mathbb{T}}\tilde{\mathfrak{h}}_{1}(\xi)\,{\rm d}\xi=0$.
Then, using similar arguments as in \cite{FGLT-2023,FGLT-2024,GL-2025,GT-2017}, we have the following proposition.
\begin{proposition}
   Let $\mathfrak{B}> 0$ and smooth function $\mathfrak{h}_{1}(\xi)$ satisfy $\bar{\mathfrak{h}}_{1}>0$ and \eqref{bw}--\eqref{h1}. Then, there exists a constant $\sigma_{0}>0$ such that for any $\sigma\in (0,\sigma_{0})$, there exists a unique solution $w(\xi,\eta)$ to \eqref{w} satisfying $\frac{\mathfrak{B}^2}{2}\leq w\leq \frac{3\mathfrak{B}^2}{2}$ and
   \begin{align*}
\sum\limits_{j+k\leq m}\int_{\mathbb{T}}\int_{0}^{\infty}|\partial_{\xi}^{j}\partial_{\eta}^{k}(w-\mathfrak{B}^2)|^2\langle \eta\rangle^{2\ell}\,{\rm d}\eta{\rm d}\xi\leq C(m,\ell)|(\mathfrak{h}_{1}^2-\mathfrak{B}^2)|_{H^{m+1}(\mathbb{T})}^2\leq C(m,\ell)\sigma^2,\quad \forall m,\ell\geq 0.
   \end{align*}
Here, $\langle \eta\rangle=(1+\eta^2)^{\frac{1}{2}}$.
\end{proposition}

Returning to the original variable $(x,z)$ and the original solution $u(x,z)$, we obtain the following corollary.
\begin{corollary}\label{cor1}
 Let $\mathfrak{B}>0$ and smooth function $\mathfrak{h}_{1}(x)$ satisfy $\bar{\mathfrak{h}}_{1}>0$ and \eqref{bw}--\eqref{h1}. Then, there exists a constant $\sigma_{0}>0$ such that for any $\sigma\in (0,\sigma_{0})$, there exists a unique solution $(u,v)(x,z)$ to \eqref{p0} with \eqref{f0}--\eqref{D1} satisfying
  \begin{align*}
\sum\limits_{j+k\leq m}\int_{\mathbb{T}}\int_{0}^{\infty}\big[|\partial_{x}^{j}\partial_{z}^{k}(u-\mathfrak{B})|^2+|\partial_{x}^{j}\partial_{z}^{k}(v-v_{\infty})|^2\big]\langle z\rangle^{2\ell}\,{\rm d}z{\rm d}x\leq C(m,\ell)\sigma^2,\quad \forall m,\ell\geq 0.
   \end{align*}
\end{corollary}
Here $v_{\infty}=\lim\limits_{z\to \infty}v(x,z)=-\int_{0}^{\infty}u_{x}(x,s)\,{\rm d}s$ is the far-field state of $v$.

\begin{remark}
We would like to point out that it is essential that $\mathfrak{B}$ should have the same sign with $\bar{\mathfrak{h}}_{1}$. Indeed, it follows from \eqref{D1} that
\begin{align*}
\int_{\mathbb{T}}u(x,0)\,{\rm d}x=\int_{\mathbb{T}}\mathfrak{h}_{1}(x)\,{\rm d}x=2\pi\bar{\mathfrak{h}}_{1},
\end{align*}
which, together with the fact that $u(x,z)$ is around the constant $\mathfrak{B}$, implies that $\mathfrak{B}$ should have the same sign with $\bar{\mathfrak{h}}_{1}$.
\end{remark}

\medskip

\subsubsection{Critical regime: $\alpha=1$.} For the well-posedness of Prandtl equation \eqref{p0} with far-field condition \eqref{f0} and the boundary condition \eqref{R1}, there is also an analogous Batchelor--Wood formula between $\mathfrak{B}$ and boundary data $\mathfrak{h}_{2}$. Unlike the explicit formula \eqref{bw}, however, their relation is difficult to derive explicitly. Moreover, due to the nonlinear relation \eqref{w0} between $w(\xi,\eta)$ and $u(x,\eta)$, such relation is not transparent in \eqref{w}. To address this issue, we denote
\begin{align*}
\mathfrak{w}(\xi,\eta)=u(\xi,\eta).
\end{align*}
Then a similar calculation shows that $\mathfrak{w}$ satisfies
\begin{align}\label{w11}
\left\{
\begin{aligned}
&\partial_{\xi}\mathfrak{w}=(\partial_{\eta}\mathfrak{w})^2+\mathfrak{w}\partial_{\eta\eta}\mathfrak{w},\\
&(\mathfrak{w}-\mathfrak{w}\partial_{\eta}\mathfrak{w})\vert_{\eta=0}=\mathfrak{h}_{2}(\xi),\quad \lim\limits_{\eta\to \infty}\mathfrak{w}=\mathfrak{B}.
\end{aligned}
\right.
\end{align}

\begin{proposition}\label{prop2}
 If $\mathfrak{w}$ is a smooth periodic solution in $\xi$ of \eqref{w11}, then it must hold that
  \begin{align}\label{cw4}
\frac{1}{2\pi}\int_{\mathbb{T}}\mathfrak{w}(\xi,0)\,{\rm d}\xi=\bar{\mathfrak{h}}_{2},\quad \frac{1}{2\pi}\int_{\mathbb{T}}\mathfrak{w}^2(\xi,0)\,{\rm d}\xi=\mathfrak{B}^2.
  \end{align}
Furthermore, if $\bar{\mathfrak{h}}_{2}=\mathfrak{B}$, \eqref{w11} has a smooth solution satisfying that $\mathfrak{w}$ doesn't change sign if and only if $\mathfrak{h}_{2}(\xi)=\mathfrak{B}$.
\end{proposition}

\noindent\textbf{Proof}. Integrating \eqref{w11} over $\mathbb{T}\times (0,\infty)$ gives
\begin{align*}
0&=\int_{\mathbb{T}}\int_{0}^{\infty}(\partial_{\eta}\mathfrak{w})^2\,{\rm d}\eta{\rm d}\xi+\int_{\mathbb{T}}\int_{0}^{\infty}\mathfrak{w}\partial_{\eta\eta}\mathfrak{w}\,{\rm d}\eta{\rm d}\xi\nonumber\\
&=-\int_{\mathbb{T}}\mathfrak{w}\partial_{\eta}\mathfrak{w}(\xi,0)\,{\rm d}\xi=\int_{\mathbb{T}}\mathfrak{h}_{2}(\xi)\,{\rm d}\xi-\int_{\mathbb{T}}\mathfrak{w}(\xi,0)\,{\rm d}\xi,
\end{align*}
which yields
\begin{align}\label{cw0}
\int_{\mathbb{T}}\mathfrak{w}(\xi,0)\,{\rm d}\xi=\int_{\mathbb{T}}\mathfrak{h}_{2}(\xi)\,{\rm d}\xi.
\end{align}
Moreover, denote
$$
\psi=\mathfrak{w}-\mathfrak{B},
$$
we obtain from \eqref{w11} that $\psi$ satisfies
\begin{align}\label{w12}
\left\{
\begin{aligned}
&\partial_{\xi}\psi=\mathfrak{B}\partial_{\eta\eta}\psi+(\partial_{\eta}\psi)^2+\psi\partial_{\eta\eta}\psi,\\
&(\psi-\mathfrak{B}\partial_{\eta}\psi)\vert_{\eta=0}=\mathfrak{h}_{2}(\xi)-\mathfrak{B}+\frac{1}{2}\partial_{\eta}(\psi^2)\vert_{\eta=0},\quad \lim\limits_{\eta\to \infty}\psi=0.
\end{aligned}
\right.
\end{align}
Integrating \eqref{w12} with respect to $\xi$ over $\mathbb{T}$ gives
\begin{align}\label{cw1}
\mathfrak{B}\partial_{\eta\eta}\int_{\mathbb{T}}\psi(\xi,\eta)\,{\rm d}\xi&=-\int_{\mathbb{T}}[(\partial_{\eta}\psi)^2+\psi\partial_{\eta\eta}\psi]\,{\rm d}\xi=-\frac{1}{2}\partial_{\eta\eta}\int_{\mathbb{T}}\psi^2\,{\rm d}\xi.
\end{align}
Denoting $H(\eta):=\mathfrak{B}\bar{\psi}(\eta)+\frac{1}{2}\bar{\psi^2}(\eta)$. The equation \eqref{cw1} together with $\eqref{w12}_{2}$ give
$$
H''(\eta)=0,\quad \text{and }\lim\limits_{\eta\to \infty}H(\eta)=0.
$$
Thus, $H(\eta)=H'(\eta)=0$, that is,
\begin{align}\label{cw2}
\mathfrak{B}\partial_{\eta}\int_{\mathbb{T}}\psi(\xi,\eta)\,{\rm d}\xi=-\frac{1}{2}\partial_{\eta}\int_{\mathbb{T}}\psi^2(\xi,\eta)\,{\rm d}\xi,\quad \mathfrak{B}\int_{\mathbb{T}}\psi(\xi,\eta)\,{\rm d}\xi=-\frac{1}{2}\int_{\mathbb{T}}\psi^2(\xi,\eta)\,{\rm d}\xi.
\end{align}
Evaluating \eqref{cw2} at $\eta=0$, one has
\begin{align*}
\mathfrak{B}\int_{\mathbb{T}}\psi(\xi,0)\,{\rm d}\xi=-\frac{1}{2}\int_{\mathbb{T}}\psi^2(\xi,0)\,{\rm d}\xi,
\end{align*}
which, together with \eqref{cw0}, yields
\begin{align}\label{cw3}
\mathfrak{B}\int_{\mathbb{T}}\psi(\xi,0)\,{\rm d}\xi=-\frac{1}{2}\int_{\mathbb{T}}\psi^2(\xi,0)\,{\rm d}\xi=\mathfrak{B}\int_{\mathbb{T}}[\mathfrak{h}_{2}(\xi)-\mathfrak{B}]\,{\rm d}\xi.
\end{align}
Rewriting $\psi=\mathfrak{w}-\mathfrak{B}$, \eqref{cw4} follows directly from \eqref{cw3}.

If $\mathfrak{h}_{2}(\xi)=\mathfrak{B}$, it is clear that \eqref{w11} has a trivial solution $\mathfrak{w}=\mathfrak{B}$. If $\mathfrak{h}_{2}(\xi)$ satisfies $\int_{\mathbb{T}}\mathfrak{h}_{2}(\xi)\,{\rm d}\xi=2\pi\mathfrak{B}$, i.e., $\int_{\mathbb{T}}(\mathfrak{h}_{2}-\mathfrak{B})\,{\rm d}\xi=0$, it follows from \eqref{cw3} that
\begin{align*}
\psi(\xi,0)=0. 
\end{align*}
Then we multiply \eqref{w12} by $\psi$ and integrate $(\xi,\eta)$ over $\mathbb{T}\times (0,\infty)$ to have
\begin{align*}
\int_{\mathbb{T}}\int_{0}^{\infty}\mathfrak{w}|\partial_{\eta}\psi|^2\,{\rm d}\eta{\rm d}\xi=0,
\end{align*}
which, together with the fact that $\mathfrak{w}$ does not change the sign, directly yields
\begin{align*}
\partial_{\eta}\psi=0\Rightarrow \mathfrak{w}=\mathfrak{B}.
\end{align*}
This, in turn, indicates that $\mathfrak{h}_{2}\equiv \mathfrak{B}$. Therefore, the proof of proposition \ref{prop2} is complete. $\hfill\square$

It follows from Proposition \ref{prop2} that the boundary data and the far-field state cannot be arbitrarily imposed. In other words, the far-field state must be appropriately chosen such that \eqref{w11} admits a nontrivial solution with desired properties. To this end, we assume that $\mathfrak{h}_{2}(\xi)$ is a small perturbation of its average over $\mathbb{T}$, i.e.,
\begin{align}\label{h2}
\mathfrak{h}_{2}(\xi)
=\frac{1}{2\pi}\int_{\mathbb{T}}\mathfrak{h}_{2}(\xi)\,{\rm d}\xi+\sigma\tilde{\mathfrak{h}}_{2}(\xi)=:\bar{\mathfrak{h}}_{2}+\sigma\tilde{\mathfrak{h}}_{2}(\xi),
\end{align}
where $\bar{\mathfrak{h}}_{2}>0$, $\int_{\mathbb{T}}\tilde{\mathfrak{h}}_{2}(\xi)\,{\rm d}\xi=0$ and $\sigma\geq 0$ is a small constant. For $\sigma=0$, it follows from Proposition \ref{prop2} that $\mathfrak{B}=\bar{\mathfrak{h}}_{2}$ and $\mathfrak{w}=\mathfrak{B}$. For $\sigma\neq 0$, we may expect $\mathfrak{w}$ in the form of
\begin{align*}
\mathfrak{w}=\mathfrak{B}+\sigma \phi,
\end{align*}
where $\phi(\xi,\eta)$ satisfies 
\begin{align}\label{P4}
	\left\{
	\begin{aligned}
		&\pa_{\xi}\phi=\mathfrak{B}\pa_{\eta\eta}^2\phi+\sigma[(\partial_{\eta}\phi)^2+\phi\partial_{\eta\eta}\phi]=\mathfrak{B}\partial_{\eta\eta}\phi+\frac{\sigma}{2}\partial_{\eta\eta}(\phi^2)\\
		&\phi(\xi,0)-\mathfrak{B}\partial_{\eta}\phi(\xi,0)=\frac{1}{\sigma}(\bar{\mathfrak{h}}_{2}-\mathfrak{B})+\tilde{\mathfrak{h}}_{2}(\xi)+\sigma(\phi\partial_{\eta}\phi)(\xi,0),\\
		&\lim\limits_{\eta\to \infty}\phi=0.
	\end{aligned}
	\right.
\end{align}
Moreover, for $\sigma\neq 0$, the compatibility condition \eqref{cw4} reads
\begin{align}\label{cw8}
\frac{1}{2\pi}\int_{\mathbb{T}}\phi(\xi,0)\,{\rm d}\xi=-\frac{\sigma}{4\pi\mathfrak{B}}\int_{\mathbb{T}}\phi^2(\xi,0)\,{\rm d}\xi=\frac{1}{\sigma}(\bar{\mathfrak{h}}_{2}-\mathfrak{B}).
\end{align}
The formula \eqref{cw8} implies that $\bar{\mathfrak{h}}_{2}-\mathfrak{B}$ should be of order $O(\sigma^2)$.

\begin{lemma}\label{lem1}
Let $\mathfrak{h}_{2}(\xi)$ be a smooth function satisfying \eqref{h2} and $\bar{\mathfrak{h}}_{2}> 0$. There exists $\sigma_0>0$ such that, for every $0\leq \sigma<\sigma_0$, there exist a constant $\mathfrak{B}=\mathfrak{B}(\sigma)$ and a smooth solution $\phi$ to \eqref{P4} satisfying \eqref{cw8} and
$$
\mathfrak{B}\geq \bar{\mathfrak{h}}_{2},\quad \mathfrak{B}+\sigma\phi\geq \frac{1}{2}\bar{\mathfrak{h}}_{2}.
$$
Moreover, $\mathfrak{B}=\bar{\mathfrak{h}}_{2}$ if and only if $\sigma\tilde{\mathfrak{h}}_{2}=0$, and for every $0<\beta<(2|\bar{\mathfrak{h}}_{2}|)^{-\frac{1}{2}}$, there is $\sigma_\beta\in(0,\sigma_0]$, independent of $m$, such that
\begin{align}\label{P4-1}
	\sum\limits_{i+j\leq m}\int_{\mathbb{T}}\int_{0}^{\infty}e^{2\beta\eta}|\partial_{\xi}^{i}\partial_{\eta}^{j}\phi|^2\,{\rm d}\eta{\rm d}\xi\leq C(m,\beta)<\infty\quad \text{for any $m\geq 0$ and $0\leq \sigma<\sigma_\beta$.}
\end{align}
\end{lemma}

\noindent\textbf{Proof}.
Let $\gamma$ be a positive constant with $\gamma<(3\bar{\mathfrak{h}}_{2})^{-\frac{1}{2}}$ and set the interval $\mathcal{I}=[\frac{\bar{\mathfrak{h}}_{2}}{2},\frac{3\bar{\mathfrak{h}}_{2}}{2}]$. It is clear that
$$
(2\mathfrak{B})^{-\frac{1}{2}}-\gamma\geq (3\bar{\mathfrak{h}}_{2})^{-\frac{1}{2}}-\gamma=:c_{\gamma}>0\quad \text{for any $\mathfrak{B}\in \mathcal{I}$}.
$$
Now, we introduce following function spaces:
$$
L_{\gamma}^2=\{f\in L^2(0,\infty)\,|\,\int_{0}^{\infty}e^{2\gamma\eta}|f(\eta)|^2\,{\rm d}\eta<\infty\},\quad H_{\gamma}^m=\{f\in H^m(0,\infty)\,|\,e^{\gamma\eta}f\in H^m(0,\infty)\},
$$
and
$$
X_{\gamma}^{m}=H^{m}(\mathbb{T};L_{\gamma}^2)\cap H^{m-\frac{1}{2}}(\mathbb{T};H_{\gamma}^1)\cap H^{m-1}(\mathbb{T};H_{\gamma}^2),
$$
equipped with norm 
$$
\|u\|_{X_{\gamma}^{m}}^2=\sum\limits_{k\in\mathbb{Z}}\Big(\la k\ra^{2m}\|u_{k}\|_{L_{\gamma}^2}^2+\la k\ra^{2m-1}\|u_{k}'\|_{L_{\gamma}^2}^2+\la k\ra^{2m-2}\|u_{k}''\|_{L_{\gamma}^2}^2\Big)
$$
for $u(\xi,\eta)=\sum\limits_{j\in \mathbb{Z}}u_{j}(\eta)e^{ij\xi}\in X_{\gamma}^{m}$. Here $\la k\ra=(1+k^2)^{\frac{1}{2}}$. We also denote
$$
Y_{\gamma}^{m-1}=H^{m-1}(\mathbb{T};L_{\gamma}^2),\quad \mathring{H}^{m}(\mathbb{T})=\{f\in H^{m}(\mathbb{T})\,|\,\int_{\mathbb{T}}f(\xi)\,{\rm d}\xi=0\}.
$$
We divide the rest of the proof into five steps.

\textbf{Step 1}. We claim that for $u\in X_{\gamma}^{m}$ with $m\geq 2$,
\begin{align}\label{P6-1}
	\sup_{\eta\geq 0}e^{\gamma\eta}\|u(\cdot,\eta)\|_{H^{m-\frac{1}{4}}(\mathbb{T})}+\sup_{\eta\geq 0}e^{\gamma\eta}\|u_{\eta}(\cdot,\eta)\|_{H^{m-\frac{3}{4}}(\mathbb{T})}\leq C(m,\gamma)\|u\|_{X_{\gamma}^{m}}.
\end{align}
In particular, it holds that $u,u_{\eta}\in L^{\infty}(\mathbb{T}\times (0,\infty))$ and
$$
u(\cdot,0)\in H^{m-\frac{1}{4}}(\mathbb{T}),\quad u_{\eta}(\cdot,0)\in H^{m-\frac{3}{4}}(\mathbb{T}),
$$
and
$$
\|u(\cdot,0)\|_{H^{m-\frac{1}{4}}}+\|u_{\eta}(\cdot,0)\|_{H^{m-\frac{3}{4}}}\leq C(m,\gamma)\|u\|_{X_{\gamma}^{m}}.
$$
In fact, let $U_{k}=e^{\gamma\eta}u_{k}(\eta)$. Noting $U_{k}\in H^{1}(0,\infty)$, we have from the Sobolev embedding that
\begin{align}\label{P6-2}
	\|U_{k}\|_{L^{\infty}}^2\leq 2\|U_{k}\|_{L^2}\|U_{k}'\|_{L^2}\leq C\|u_{k}\|_{L_{\gamma}^2}\big(\|u_{k}\|_{L_{\gamma}^2}+\|u_{k}'\|_{L_{\gamma}^2}\big).
\end{align}
Multiplying \eqref{P6-2} by $\la k\ra^{2m-\frac{1}{2}}$, and summing in $k$ gives
\begin{align}\label{P6-3}
	\sup_{\eta\geq 0}e^{2\gamma\eta}\|u(\cdot,\eta)\|_{H^{m-\frac{1}{4}}}^2\leq C\|u\|_{X_{\gamma}^{m}}^2.
\end{align}
Applying similar arguments to $u_{k}'$, we also have
\begin{align*}
	\sup_{\eta\geq 0}e^{2\gamma\eta}\|u_{\eta}(\cdot,\eta)\|_{H^{m-\frac{3}{4}}}^2\leq C\|u\|_{X_{\gamma}^{m}}^2,
\end{align*}
which together with \eqref{P6-3} yields \eqref{P6-1}. Moreover, by the Sobolev embedding, we have directly $u,u_{\eta}\in L^{\infty}(\mathbb{T}\times (0,\infty))$. This concludes the claim.

For later use, let $X=X_{\gamma}^2$ and $Z=Y_{\gamma}^{1}\times \mathring{H}^{\frac{5}{4}}(\mathbb{T})$. Consider the nonlinear map:
$$
\mathcal{K}u=\Big[\frac{1}{2}\partial_{\eta\eta}(u^2),\Pi_{\perp}\big(\frac{1}{2}\partial_{\eta}(u^2)\vert_{\eta=0}\big)\Big]^{\top},
$$
where $\Pi_{\perp}$ is the projection operator defined by
$
\Pi_{\perp}f:=f-\bar{f}(\eta).
$
It follows directly from the above claim that $\mathcal{K}$ is well defined from $X$ to $Z$. Moreover, a direct calculation shows that $\mathcal{K}$ is a $C^{\infty}$ map, and for any $u,v\in X$
\begin{align*}
	\|\mathcal{K}(u)\|_{Z}\leq C\|u\|_{X}^2,\quad \|\mathcal{K}(u)-\mathcal{K}(v)\|_{Z}\leq C(\|u\|_{X}+\|v\|_{X})(\|u-v\|_{X}).
\end{align*}

\textbf{Step 2}. Now, for $\mathfrak{B}\in \mathcal{I}$, we consider the linear operator:
\begin{align*}
	\mathcal{L}_{\mathfrak{B}}: u\mapsto [\partial_{\xi}u-\mathfrak{B}\partial_{\eta\eta}u,\Pi_{\perp}\big((u-\mathfrak{B}\partial_{\eta}u)\vert_{\eta=0}\big)]^{\top},
\end{align*}
It is clear that $\mathcal{L}_{\mathfrak{B}}$ is a $C^{\infty}$ map from $X_{\gamma}^{m}$ to $Y_{\gamma}^{m-1}\times \mathring{H}^{m-\frac{3}{4}}(\mathbb{T})$. We claim that $\mathcal{L}_{\mathfrak{B}}$ is an isomorphism. In fact, for any $[F,h]\in Y_{\gamma}^{m-1}\times \mathring{H}^{m-\frac{3}{4}}(\mathbb{T})$, set $\mathcal{L}(u)=[F,h]^{\top}$, that is, $u$ satisfies the following equation
\begin{align*}
	\left\{
	\begin{aligned}
		&u_{\xi}-\mathfrak{B}u_{\eta\eta}=F,\\
		&\Pi_{\perp}\big((u-\mathfrak{B}\partial_{\eta}u)\vert_{\eta=0}\big)=h.
	\end{aligned}
	\right.
\end{align*}
Using the Fourier transform, one has
$$
u(\xi,\eta)=\sum\limits_{k\in \mathbb{Z}}u_{k}(\eta)e^{ik\xi},\quad F(\xi,\eta)=\sum\limits_{k\in \mathbb{Z}}F_{k}(\eta)e^{ik\xi},\quad h(\xi)=\sum\limits_{k\neq 0}h_{k}e^{ik\xi}.
$$
Then for $k\neq 0$, $u_{k}(\eta)$ satisfies
\begin{align*}
	ik u_{k}-\mathfrak{B}u_{k}''=F_{k}(\eta),\quad (u_{k}-\mathfrak{B}u_{k}')\vert_{\eta=0}=h_{k},\quad \lim\limits_{\eta\to \infty}u_{k}=0,
\end{align*}
which implies
\begin{align}\label{P6-5}
	u_{k}(\eta)=\int_0^{\infty}G_{k}(\eta,s)F_{k}(s)\,{\rm d}s+\frac{h_{k}}{1+\mathfrak{B}\lambda_{k}}e^{-\lambda_{k}\eta},
\end{align}
where $\lambda_{k}=\sqrt{\frac{ik}{\mathfrak{B}}}$, and the Green's kernel
$$
G_{k}(\eta,s)=\frac{1}{2\mathfrak{B}\lambda_{k}}\big(e^{-\lambda_{k}|\eta-s|}+\frac{\mathfrak{B}\lambda_{k}-1}{\mathfrak{B}\lambda_{k}+1}e^{-\lambda_{k}(\eta+s)}\big).
$$
And for $k=0$, $u_{0}(\eta)$ satisfies
$$
\mathfrak{B}u_{0}''=-F_{0}(\eta),\quad \lim\limits_{\eta\to \infty}u_{0}(\eta)=0,
$$
which implies that
\begin{align}\label{P6-6}
	u_{0}(\eta)=-\frac{1}{\mathfrak{B}}\int_{\eta}^{\infty}\int_{s}^{\infty}F_{0}(t)\,{\rm d}t{\rm d}s.
\end{align}
It is worth pointing out that there is no boundary condition at $\eta=0$ for $u_{0}(\eta)$ but the solution $u_{0}(\eta)$ is uniquely determined by \eqref{P6-6} in the space $X_{\gamma}^{m}$ due to the far field decay condition.

A direct calculation shows that
$$
\lambda_{k}=\sqrt{\frac{ik}{\mathfrak{B}}}=\sqrt{\frac{|k|}{2\mathfrak{B}}}+i\operatorname{sgn}\big(k\big)\sqrt{\frac{|k|}{2\mathfrak{B}}}\Rightarrow \mathfrak{R}(\lambda_{k})=\sqrt{\frac{|k|}{2\mathfrak{B}}}\geq \sqrt{\frac{1}{2\mathfrak{B}}}\quad \text{for }k\in \mathbb{Z}\backslash \{0\}.
$$
Noting that $\mathfrak{R}(\lambda_{k})-\gamma>\tilde{c}\mathfrak{R}(\lambda_{k})$ uniformly for $\mathfrak{B}\in \mathcal{I}$, we have from \eqref{P6-5} that
\begin{align}\label{P6-7-0}
\langle k\rangle\|u_{k}\|_{L_{\gamma}^{2}}+\langle k\rangle^{\frac{1}{2}}\|u_{k}'\|_{L_{\gamma}^2}+\|u_{k}''\|_{L_{\gamma}^2}\leq C\|F_{k}\|_{L_{\gamma}^2}+C\langle k\rangle^{\frac{1}{4}}|h_{k}|.
\end{align}
Multiplying the \eqref{P6-7-0} by $\langle k\rangle^{m-1}$ and summing over $k\in\mathbb{Z}\backslash\{0\}$ gives
that
\begin{align}\label{P6-7}
\|\Pi_{\perp}u\|_{X_{\gamma}^{m}}\leq C(m,\gamma)\big(\|F\|_{Y_{\gamma}^{m-1}}+\|h\|_{H^{m-\frac{3}{4}}}\big),
\end{align}
where $C$ is a uniform constant for $\mathfrak{B}\in \mathcal{I}$.

For $u_{0}(\eta)$, by \eqref{P6-6}, a direct calculation shows that
\begin{align}\label{P6-8}
\gamma^2\|u_{0}\|_{L_{\gamma}^2}+\gamma\|u_{0}'\|_{L_{\gamma}^2}+\|u_{0}''\|_{L_{\gamma}^2}\leq \frac{C}{\mathfrak{B}}\|F_{0}\|_{L_{\gamma}^2}.
\end{align}
Combining \eqref{P6-7} and \eqref{P6-8}, we have
\begin{align}\label{P6-9}
	\|u\|_{X_{\gamma}^{m}}\leq C\big(\|F\|_{Y_{\gamma}^{m-1}}+\|h\|_{H^{m-\frac{3}{4}}}\big),
\end{align}
which implies that $\mathcal{L}_{\mathfrak{B}}$ is surjective from $X_{\gamma}^{m}$ to $Y_{\gamma}^{m-1}\times \mathring{H}^{m-\frac{3}{4}}$. Here, $C$ is a uniform constant for $\mathfrak{B}\in \mathcal{I}$.

To show $\mathcal{L}_{\mathfrak{B}}$ is injective, we set $\mathcal{L}_{\mathfrak{B}}u=[0,0]$, then we obtain from \eqref{P6-5} and \eqref{P6-6} that $u_{k}=0$ for all $k\in \mathbb{Z}$. Therefore, $u\equiv 0$, which concludes that $\mathcal{L}_{\mathfrak{B}}$ is injective. Hence, $\mathcal{L}_{\mathfrak{B}}$ is an isomorphism. Moreover, \eqref{P6-9} implies that
\begin{align*}
\|\mathcal{L}_{\mathfrak{B}}^{-1}\|_{Y_{\gamma}^{m-1}\times \mathring{H}^{m-\frac{3}{4}}\longmapsto X_{\gamma}^{m}}\leq C.
\end{align*}
Here, $C$ is a uniform constant for $\mathfrak{B}\in \mathcal{I}$.

\textbf{Step 3}. For any $q\in \mathring{H}^{\frac{5}{4}}(\mathbb{T})$, we consider the operator:
\begin{align*}
	\mathcal{F}(u;\mathfrak{B},\sigma,q)
    &:=\mathcal{L}_{\mathfrak{B}}u-\sigma \mathcal{K}u-[0,q]^{\top}.
\end{align*}
Then $\mathcal{F}$ is a $C^{\infty}$ map into $Z$. When $(\mathfrak{B},\sigma,q)=(\bar{\mathfrak{h}}_{2},0,\tilde{\mathfrak{h}}_{2})$, we obtain from $\mathcal{F}(u;\mathfrak{B},\sigma,q)=[0,0]^{\top}$ that
$$
u=\phi^{0}:=\mathcal{L}_{\bar{\mathfrak{h}}_{2}}^{-1}[0,\tilde{\mathfrak{h}}_{2}]^{\top}\Rightarrow \phi_{k}^{0}=\frac{\tilde{\mathfrak{h}}_{2,k}}{1+\bar{\mathfrak{h}}_{2}\sqrt{ik\slash\bar{\mathfrak{h}}_{2}}}e^{-\sqrt{ik\slash \bar{\mathfrak{h}}_{2}}\eta}\,\,(k\neq 0),\quad \phi_{0}^{0}=0,
$$
where we have used \eqref{P6-5}--\eqref{P6-6}. Moreover, a direct calculation shows that the partial derivative
$$
D_{u}\mathcal{F}\vert_{(\phi^{0},\bar{\mathfrak{h}}_{2},0,\tilde{\mathfrak{h}}_{2})}=\mathcal{L}_{\bar{\mathfrak{h}}_{2}}.
$$
Thus, $D_{\phi}\mathcal{F}$ is an isomorphism at $(\phi^{0},\bar{\mathfrak{h}}_{2},0,\tilde{\mathfrak{h}}_{2})$. The implicit function theorem gives a $C^{\infty}$ map
\begin{align*}
(\mathfrak{B},\sigma,q)\longmapsto U(\mathfrak{B},\sigma,q)\in X
\end{align*}
in a neighborhood of $(\bar{\mathfrak{h}}_{2},0,\tilde{\mathfrak{h}}_{2})$, which is unique in a fixed and small neighborhood of $\phi^{0}$. In particular, there exist $\delta>0$ and $s>0$ so that $\mathcal{J}=[\bar{\mathfrak{h}}_{2}-\delta,\bar{\mathfrak{h}}_{2}+\delta]\subset \mathcal{I}$, and $\mathcal{J}\times [-s,s]\times \{\tilde{\mathfrak{h}}_{2}\}$ lies strictly inside such neighborhood of $(\bar{\mathfrak{h}}_{2},0,\tilde{\mathfrak{h}}_{2})$. Let
$$
\Phi(\mathfrak{B},\sigma)=U(\mathfrak{B},\sigma,\tilde{\mathfrak{h}}_{2})\quad \text{for $(\mathfrak{B},\sigma)\in \mathcal{J}\times [-s,s]$}.
$$
Then $\Phi(\mathfrak{B},\sigma)\in C^{\infty}(\mathcal{J}\times [-s,s];X)$ satisfies
\begin{align}\label{phisigma0}
\left\{
\begin{aligned}
&\partial_{\xi}\Phi=\mathfrak{B}\partial_{\eta\eta}\Phi+\frac{\sigma}{2}\partial_{\eta\eta}(\Phi^2),\\
&\Pi_{\perp}\big[\big(\Phi(\xi,0)-\mathfrak{B}\partial_{\eta}\Phi-\frac{\sigma}{2}\partial_{\eta}((\Phi)^2)\big)\vert_{\eta=0}\big]=\tilde{\mathfrak{h}}_{2}(\xi),
\end{aligned}
\right.
\end{align}
and there exists a uniform constant $C$ such that
\begin{align}\label{Phi}
\|\Phi(\mathfrak{B},\sigma)\|_{X}\leq C,\,\, \|\Phi(\mathfrak{B},\sigma)-\Phi(\mathfrak{B}',\sigma)\|_{X}\leq C|\mathfrak{B}-\mathfrak{B}'|\quad \text{for $\mathfrak{B},\mathfrak{B}'\in \mathcal{J}$ and $\sigma\in [-s,s]$}.
\end{align}
Furthermore, by reducing $\delta$ and $s$ if necessary, one has from \eqref{Phi} that
\begin{align}\label{coer}
\mathfrak{B}+\sigma \Phi(\mathfrak{B},\sigma)\geq \frac{\bar{\mathfrak{h}}_{2}}{2}.
\end{align}

\textbf{Step 4}. Now, for any $\sigma\in [-s,s]$, we shall show the uniform-in-$\sigma$ higher regularity of the derivatives of $\Phi$ with respect to $\xi$ and $\eta$. Let 
$$
T_{t}f(\xi,\eta):=f(\xi+t,\eta)\text{ and }q_{t}:=T_{t}\tilde{\mathfrak{h}}_{2}=\tilde{\mathfrak{h}}_{2}(\xi+t).
$$
Noting that $\mathcal{L}_{\mathfrak{B}}$ and $\mathcal{K}$ are both translation invariant, one has
$$
\mathcal{F}(T_{t}u;\mathfrak{B},\sigma,q_{t})=T_{t}\mathcal{F}(u,\mathfrak{B},\sigma,\tilde{\mathfrak{h}}_{2}).
$$
Since $\mathcal{J}\times [-s,s]$ is compact and $\Phi$ is smooth over $(\mathfrak{B},\sigma)\in \mathcal{J}\times [-s,s]$, the set
$$
\mathcal{D}=\{\Phi(\mathfrak{B},\sigma)\,|\,(\mathfrak{B},\sigma)\in \mathcal{J}\times [-s,s]\}
$$
is compact in $X$. Then for any $\epsilon>0$, there exists a finite $\epsilon-$net $\{u_{1},u_{2},\cdots,u_{N}\}$ for $\mathcal{D}$. Since $T_{t}u\to u$ in $X$ as $t\to 0$, there exists $t_{0}>0$ such that
\begin{align*}
    \|T_{t}u_{j}-u_{j}\|_{X}<\epsilon\quad \text{for every $j$ whenever $|t|<t_0$}.
\end{align*}
For $u\in \mathcal{D}$, choosing $u_{j}$ such that $\|u_{j}-u\|_{X}<\epsilon$, one has from the fact $T_{t}$ is an isometry on $X$ that
\begin{align*}
    \|T_{t}u-u\|_{X}\leq \|T_{t}(u-u_{j})\|_{X}+\|T_{t}u_{j}-u_{j}\|_{X}+\|u_{j}-u\|_{X}\leq 3\epsilon\quad \text{whenever $|t|<t_{0}$},
\end{align*}
which implies
\begin{align}\label{translation}
    \sup\limits_{u\in\mathcal{D}}\|T_{t}u-u\|_{X}\to 0\quad \text{as $t\to 0$}.
\end{align}

Next, it follows from $q_{t}\to \tilde{\mathfrak{h}}_{2}$ in $\mathring{H}^{\frac{5}{4}}$ that there exists a sufficiently small $t_{0}$ such that $(\mathfrak{B},\sigma,q_{t})$ lies inside the small neighborhood of $(\bar{\mathfrak{h}}_{2},0,\tilde{\mathfrak{h}}_{2})$ for all $\mathfrak{B}\in \mathcal{J}$, $|\sigma|\leq s$ and $|t|<t_{0}$. Then, using the strong continuity \eqref{translation}, one has
$$
\mathcal{F}(T_{t}\Phi(\mathfrak{B},\sigma);\mathfrak{B},\sigma,q_{t})=0,\quad \text{for $|t|<t_{0}$}.
$$
So, the local uniqueness of the implicit function theorem gives 
\begin{align}\label{TPhi}
T_{t}\Phi(\mathfrak{B},\sigma)=U(\mathfrak{B},\sigma,q_{t}),\quad \text{for $\mathfrak{B}\in \mathcal{J}$, $|\sigma|\leq s$ and $|t|<t_{0}$}.
\end{align}
Here, the smallness of $t_{0}$ is uniform $(\mathfrak{B},\sigma)\in \mathcal{J}\times [-s,s]$. 

Since $U(\mathfrak{B},\sigma,q_{t})$ is $C^{\infty}$, it follows from \eqref{TPhi} that 
\begin{align*}
    \partial_{\xi}^{n}\Phi(\mathfrak{B},\sigma)=\frac{d^{n}}{dt^{n}}T_{t}\Phi(\mathfrak{B},\sigma)\Big\vert_{t=0}=\frac{d^{n}}{dt^{n}}U(\mathfrak{B},\sigma,q_{t})\Big\vert_{t=0}\in X,
\end{align*}
for any integer $n\geq 0$, and 
\begin{align}\label{tangential}
    \sup_{\mathfrak{B}\in \mathcal{J},|\sigma|\leq s}\|\partial_{\xi}^{n}\Phi(\mathfrak{B},\sigma)\|_{X}=\sup_{\mathfrak{B}\in \mathcal{J},|\sigma|\leq s}\Big\|\frac{d^{n}}{dt^{n}}U(\mathfrak{B},\sigma,q_{t})\Big\vert_{t=0}\Big\|\leq C(n,\gamma)<\infty.
\end{align}

To control the normal derivatives, we rewrite $\eqref{phisigma0}_{1}$ as
\begin{align}\label{phisigma2}
\partial_{\eta\eta}\Phi=(\mathfrak{B}+\sigma\Phi)^{-1}(\Phi_{\xi}-\sigma(\partial_{\eta}\Phi)^2).
\end{align}
Then applying $\partial_{\xi}^{i}\partial_{\eta}^{j}$ to \eqref{phisigma2}, and using \eqref{Phi}, \eqref{coer} and \eqref{tangential}, we obtain from the induction method that
\begin{align}\label{fullder}
\sum\limits_{i+j\leq m}\|\partial_{\xi}^{i}\partial_{\eta}^{j}\Phi(\mathfrak{B},\sigma)\|_{L_{\gamma}^2}\leq C(m,\gamma)<\infty,
\end{align}
for all $(\mathfrak{B},\sigma)\in \mathcal{J}\times [-s,s]$.

\textbf{Step 5}. Then, for $\sigma\neq 0$, we consider the following equation about $\mathfrak{B}$:
\begin{align}\label{frakB}
\frac{1}{\sigma}(\bar{\mathfrak{h}}_{2}-\mathfrak{B})&=\frac{1}{2\pi}\int_{\mathbb{T}}(\Phi-\mathfrak{B}\partial_{\eta}\Phi-\sigma\Phi\partial_{\eta}\Phi)\vert_{\eta=0}\,{\rm d}\xi\nonumber\\
&=\frac{1}{2\pi}\int_{\mathbb{T}}\Phi(\xi,0)\,{\rm d}\xi=-\frac{\sigma}{4\pi\mathfrak{B}}\int_{\mathbb{T}}|\Phi(\xi,0)|^2\,{\rm d}\xi.
\end{align}
The last two equalities follow directly from the compatibility condition of \eqref{phisigma0}. Since, $\eqref{phisigma0}$ involves the parameter $\mathfrak{B}\in \mathcal{J}$, we must show the equation \eqref{frakB} can admit a solution $\mathfrak{B}$ satisfying $\mathfrak{B}\geq \bar{\mathfrak{h}}_{2}$. 

Let $\mathcal{J}_{+}=[\bar{\mathfrak{h}}_{2},\bar{\mathfrak{h}}_{2}+\delta]\subset \mathcal{J}$. We construct the iterative sequence $\{\mathfrak{B}_{k}\}$ through
\begin{align*}
\mathfrak{B}_{k+1}=\bar{\mathfrak{h}}_{2}+\frac{\sigma^2}{4\pi\mathfrak{B}_{k}}\int_{\mathbb{T}}|\Phi(\mathfrak{B}_{k},\sigma)(\xi,0)|^2\,{\rm d}\xi,
\end{align*}
with $\mathfrak{B}_{0}=\bar{\mathfrak{h}}_{2}>0$. It is clear that $\mathfrak{B}_{k}\geq \bar{\mathfrak{h}}_{2}>0$ for any $k\geq 0$ by the induction. Moreover, it follows from \eqref{Phi} that
\begin{align*}
\mathfrak{B}_{k+1}\leq \bar{\mathfrak{h}}_{2}+C\sigma^2.
\end{align*}
Choosing $\sigma_{0}\in (0,s)$ so that $C\sigma_{0}^2\leq \delta$, we obtain that $\mathfrak{B}_{k}\in \mathcal{J}_{+}$ for any $k\geq 0$. Next, using \eqref{Phi}, a direct calculation shows that
\begin{align*}
|\mathfrak{B}_{k+1}-\mathfrak{B}_{k}|\leq C\sigma^2|\mathfrak{B}_{k}-\mathfrak{B}_{k-1}|\leq \frac{1}{2}|\mathfrak{B}_{k}-\mathfrak{B}_{k-1}|
\end{align*}
provided that $C\sigma_{0}^2\leq \frac{1}{2}$ and $\sigma\in (-\sigma_{0},\sigma_{0})$. Therefore, $\{\mathfrak{B}_{k}\}$ is a Cauchy sequence and there exists a constant $\mathfrak{B}\in \mathcal{J}_{+}$ such that \eqref{frakB} holds. 

Set
$$
\phi^{\sigma}=\Phi(\mathfrak{B}(\sigma),\sigma).
$$
Combining \eqref{phisigma0} and \eqref{frakB}, $\phi^{\sigma}$ will satisfy \eqref{P4}.
It follows from the Proposition \ref{prop2} that $\mathfrak{B}=\bar{\mathfrak{h}}_{2}$ if and only if $\sigma\tilde{\mathfrak{h}}_{2}=0$. 

Finally, let $0<\beta<(2\bar{\mathfrak{h}}_{2})^{-\frac{1}{2}}$. Noting from \eqref{Phi} and \eqref{frakB} that $\mathfrak{B}-\bar{\mathfrak{h}}_{2}=O(\sigma^2)$, there exists a positive gap between $\beta$ and $(2\mathfrak{B})^{-\frac{1}{2}}$. Then, by similar calculations in Steps 3--4, we find a constant $\sigma_{\beta}$ with $0<\sigma_{\beta}\leq \sigma_{0}$ so that  the estimate \eqref{P4-1} follows from \eqref{fullder}. Therefore, the proof of Lemma \ref{lem1} is complete.  $\hfill\square$

Returning to the original
coordinates $(x,z)$ and using the fact $\eta=\mathfrak{B}z+O(\sigma)$, we obtain the following corollary.

\begin{corollary}\label{cor2}
Let $\mathfrak{h}_{2}(x)$ be a smooth function satisfying \eqref{h2} and $\bar{\mathfrak{h}}_{2}>0$. There exists $\sigma_0>0$ such that, for every $0\leq \sigma<\sigma_0$, there exist a constant $\mathfrak{B}=\mathfrak{B}(\sigma)$ satisfying
$$
\mathfrak{B}\geq \bar{\mathfrak{h}}_{2},\quad \bar{\mathfrak{h}}_{2}-\mathfrak{B}=\frac{\sigma}{2\pi}\int_{\mathbb{T}}\phi(\xi,0)\,{\rm d}\xi=-\frac{\sigma^2}{4\pi\mathfrak{B}}\int_{\mathbb{T}}\phi^2(\xi,0)\,{\rm d}\xi,
$$
where $\phi$ is a smooth solution of \eqref{phisigma0}, and a smooth solution $(u,v)(x,z)$ of \eqref{p0} with \eqref{f0} and \eqref{R1}. 
Moreover, $\mathfrak{B}=\bar{\mathfrak{h}}_{2}$ if and only if $\sigma\tilde{\mathfrak{h}}_{2}=0$, and 
for every $0<\beta<\sqrt{|\bar{\mathfrak{h}}_{2}|/2}$, there is $\sigma_\beta\in(0,\sigma_0]$, independent of $m$, such that
$$
\sum\limits_{i+j\leq m}\int_{\mathbb{T}}\int_{0}^{\infty} e^{2\beta z}[|\partial_{x}^{i}\partial_{z}^{j}(u-\mathfrak{B})|^2+|\partial_{x}^{i}\partial_{z}^{j}(v-v_{\infty})|^2]\,{\rm d}z{\rm d}x\leq C(m,\beta)\sigma^2
$$
for $m\geq 0$ and $0\leq \sigma<\sigma_{\beta}$. Here, $v_{\infty}=\lim\limits_{z\to \infty}v(x,z)=-\int_{0}^{\infty}\partial_{x}u(x,s)\,{\rm d}s$ is the far-field of $v$.
\end{corollary}

\subsubsection{Subcritical regime: $0\leq \alpha<1$} Now, we consider the well-posedness of \eqref{p0} with homogeneous Neumann-Dirichlet boundary conditions:
\begin{align}\label{N1}
\partial_{z}u(x,0)=0,\quad v(x,0)=0,
\end{align}
and the far-field condition \eqref{f0}.

\begin{lemma}\label{lem2}
The smooth solution $(u,v)$ satisfying $u-\mathfrak{B}\in H^1(\mathbb{T}\times (0,\infty))$ to the equations \eqref{p0} with \eqref{f0} and \eqref{N1} must be $(u,v)=(\mathfrak{B},0)$.
\end{lemma}

\noindent\textbf{Proof}. We divide the proof into three steps.

\noindent\textbf{Step 1}. Denote $w=u-\mathfrak{B}$. For any $N>0$, we multiply \eqref{p0} by $w$ and integrate $(x,z)$ over $\mathbb{T}\times (0,N)$ to get
\begin{align}\label{N2}
\int_{0}^{N}\int_{\mathbb{T}}|w_{z}|^2\,{\rm d}x{\rm d}z=\int_{\mathbb{T}}ww_{z}(x,N)\,{\rm d}x-\frac{1}{2}\int_{\mathbb{T}}vw^{2}(x,N)\,{\rm d}x.
\end{align}
Let
$$
\psi(x,N)=\int_{0}^{N}w(x,s)\,{\rm d}s,
$$
then $v(x,N)=-\psi_{x}(x,N)$ and
\begin{align}\label{N3}
\int_{\mathbb{T}}vw^{2}(x,N)\,{\rm d}x=-\int_{\mathbb{T}}(\psi_{x}w^{2})(x,N)\,{\rm d}x=2\int_{\mathbb{T}}(\psi ww_{x})(x,N)\,{\rm d}x.
\end{align}
Substituting \eqref{N3} into \eqref{N2} gives
\begin{align}\label{N4}
   \int_{0}^{N}\int_{\mathbb{T}}|w_{z}|^2\,{\rm d}x{\rm d}z=\int_{\mathbb{T}}(ww_{z})(x,N)\,{\rm d}x-\int_{\mathbb{T}}(\psi ww_{x})(x,N)\,{\rm d}x.
\end{align}

\noindent\textbf{Step 2}. We claim that there exists a sequence $N_{k}$ satisfying $\lim\limits_{k\to \infty}N_{k}= \infty$ such that
\begin{align}\label{N5}
    \lim\limits_{k\to \infty}N_{k}\big(\|w_{x}(x,N_{k})\|_{L^2}^2+\|w_{z}(x,N_{k})\|_{L^2}^2\big)=0\quad \text{ and }\quad \sup_{k\geq 1}\|w(\cdot,N_{k})\|_{L^{\infty}}<\infty.
\end{align}
In fact, denoting $D=\mathbb{T}\times (0,\infty)$, since $w\in H^1(D)$, one has from the Sobolev embedding that
\begin{align}\label{N5-1}
\sup_{z\geq 0}\|w(\cdot,z)\|_{L^2}\leq C\|w\|_{H^1(D)}.
\end{align}
For almost every $z\in (0,\infty)$, we define
\begin{align*}
h(z):=\|w_{x}(\cdot,z)\|_{L^2}^2+\|w_{z}(\cdot,z)\|_{L^2}^2.
\end{align*}
Then, $h(z)\in L^{1}(0,\infty)$, and
\begin{align*}
    a_{k}:=\int_{k}^{2k}h(z)\,{\rm d}z\to 0\quad \text{as $k\to \infty$}.
\end{align*}
Let
$$
G_{k}=\big\{z\in[k,2k]:h(z)\leq \frac{2a_{k}}{k}\big\}.
$$
If $a_{k}=0$, then $G_{k}=[k,2k]$. If $a_{k}>0$, then, for any $z\not\in G_{k}$, one has $h(z)>\frac{2a_{k}}{k}$, and 
$$
a_{k}\geq \int_{G_{k}^{c}}h(z)\,{\rm d}z\geq \frac{2a_{k}}{k}|G_{k}^{c}|\Rightarrow |G_{k}^{c}|\leq \frac{k}{2},
$$
where $|G_{k}^{c}|$ denotes the Lebesgue measure of the complement set of $G_{k}$ relative to $[k,2k]$. Hence, $|G_{k}|>0$. 

On the other hand, by Fubini's theorem, it follows from $w\in H^{1}(D)$ that there is set $E\subset (0,\infty)$ with $|E^{c}|=0$ such that $w(\cdot,z)\in H^{1}(\mathbb{T})$ with its weak derivative given by $w_{x}(x,z)$ for $z\in E$. We therefore can choose
$$
N_{k}\in G_{k}\cap E
$$
such that $k\leq N_{k}\leq 2k$ and
\begin{align*}
    N_{k}h(N_{k})\leq 4a_{k}\to 0\quad \text{as $k\to \infty$}.
\end{align*}
In particular,
\begin{align}\label{N6}
\|w_{x}(\cdot,N_{k})\|_{L^2(\mathbb{T})}+\|w_{z}(\cdot,N_{k})\|_{L^2(\mathbb{T})}\to 0\quad \text{as $k\to \infty$}.
\end{align}
Moreover, since $w(\cdot,N_{k})\in H^{1}(\mathbb{T})$, one has from the Sobolev embedding that
\begin{align*}
\|w(\cdot,N_{k})\|_{L^{\infty}}\leq C\|w(\cdot,N_{k})\|_{H^{1}(\mathbb{T})},
\end{align*}
which, together with \eqref{N5-1} and \eqref{N6}, yields that
\begin{align*}
    \sup_{k\geq 1}\|w(\cdot,N_{k})\|_{L^{\infty}}<\infty.
\end{align*}
This concludes \eqref{N5}. 

\noindent\textbf{Step 3}. Noting that $\|\psi(\cdot,N_{k})\|_{L^2}\leq \sqrt{N_{k}}\|w\|_{L^2(D)}$, and using \eqref{N5}--\eqref{N5-1}, a direct calculation shows that
\begin{align}\label{N7}
    \Big\vert\int_{\mathbb{T}}(\psi w w_{x})(x,N_{k})\,{\rm d}x\Big\vert&\leq 2\|\psi(\cdot,N_{k})\|_{L^2}\|w(\cdot,N_{k})\|_{L^{\infty}}\|w_{x}(\cdot,N_{k})\|_{L^2}\nonumber\\
    &\leq C\sqrt{N_{k}h(N_{k})}\to 0\quad \text{as $k\to \infty$}.
\end{align}
and
\begin{align}\label{N8}
    \Big\vert\int_{\mathbb{T}}w(x,N_{k})w_{z}(x,N_{k})\,{\rm d}x\Big\vert\leq C\|w(\cdot,N_{k})\|_{L^{\infty}}\|w_{z}(\cdot,N_{k})\|_{L^2}\leq C\sqrt{h(N_{k})}\to 0\quad \text{as $k\to \infty$}.
\end{align}
Choosing $N=N_{k}$ in \eqref{N4}, and then taking limit $k\to\infty$, we obtain from \eqref{N7}--\eqref{N8} that
\begin{align*}
\lim\limits_{k\to \infty}\int_{0}^{N_{k}}\int_{\mathbb{T}}|w_{z}|^2\,{\rm d}x{\rm d}z=\int_{0}^{\infty}\int_{\mathbb{T}}|w_{z}|^2\,{\rm d}x{\rm d}z=0,
\end{align*}
which implies that $w_{z}\equiv 0$. Using the far-field condition \eqref{f0}, we further obtain $w\equiv 0$. This gives $u\equiv \mathfrak{B}$ and $v\equiv 0$. Therefore, the proof of Lemma \ref{lem2} is complete. $\hfill\square$

\subsection{Well-posedness of higher-order boundary-layer profiles} Recall \eqref{1.15} with \eqref{b00}--\eqref{b3}, for any $\alpha\geq 0$, the higher-order boundary-layer profile $(u_{k},v_{k+\tilde{\alpha}})$ will satisfy the linearized Prandtl equation around the zeroth-order boundary-layer profile $(u,v)$ constructed in Corollary \ref{cor1} for $\alpha>1$, Corollary \ref{cor2} for $\alpha=1$ and Lemma \ref{lem2} for $0\leq \alpha<1$:
\begin{align}\label{lp0}
\left\{
\begin{aligned}
&u\partial_{x}u_{k}+u_{k}\partial_{x}u+(v_{k+\tilde{\alpha}}-v_{k+\tilde{\alpha}}(x,0))\partial_{z}u+v\partial_{z}u_{k}=\partial_{zz}u_{k}+\mathfrak{f}_{k},\\
&\partial_{x}u_{k}+\partial_{z}v_{k+\tilde{\alpha}}=0,
\end{aligned}
\right.
\end{align}
supplemented with the Robin boundary condition
\begin{align}\label{b2-1}
u_{k}(x,0)-\partial_{z}u_{k}(x,0)=\mathfrak{R}_{1}(x)\quad \text{if $\alpha=1$},
\end{align}
or the Dirichlet boundary condition
\begin{align}\label{b3-1}
u_{k}(x,0)=\mathfrak{R}_{2}(x)\quad \text{if $\alpha>1$},
\end{align}
or the Neumann boundary conditions
\begin{align}\label{b1-1}
\partial_{z}u_{k}(x,0)=\mathfrak{R}_{3}(x)\quad \text{if $0\leq \alpha<1$}.
\end{align}
We again deal with the case $\alpha\geq 1$, and then deal with $0\leq \alpha<1$.

\medskip

\subsubsection{Critical and supercritical regimes: $\alpha\geq 1$.} For the system \eqref{lp0} with \eqref{b2-1} or \eqref{b3-1}, as mentioned in \cite{FGLT-2023,FGLT-2024,GL-2025}, there is an implicit compatibility condition between the far-field condition and the boundary data such that we are not able to expect
$$
\lim\limits_{z\to \infty}u_{k}(x,z)=0.
$$
Motivated by \cite{FGLT-2023,FGLT-2024,GL-2025}, we impose the following far-field conditions:
\begin{align}\label{4-2}
\lim\limits_{z\to \infty}(\partial_{z}u_{k},v_{k+\tilde{\alpha}})=(0,0).
\end{align}
For later use, we denote
$$
\mathfrak{R}(x)=
\left\{
\begin{aligned}
&\mathfrak{R}_{1}(x),\quad \text{in case of \eqref{b2-1}},\\
&\mathfrak{R}_{2}(x),\quad \text{in case of \eqref{b3-1}}.
\end{aligned}
\right.
$$

\begin{lemma}\label{lem3}
Let $\mathfrak{R}$ and $\mathfrak{f}_{k}$ be smooth, periodic in $x$ satisfying
\begin{align*}
\sum\limits_{i+j\leq m}\big[\|\langle z\rangle^{\ell}\partial_{x}^{i}\partial_{z}^{j}\mathfrak{f}_{k}\|_{L^2(D)}+\|\partial_{x}^{i}\mathfrak{R}\|_{L^2(\mathbb{T})}\big]<\infty,\quad \text{for any integers $m,\ell\geq 0$}.
\end{align*}
Here $D=\mathbb{T}\times (0,\infty)$. Then, there exists a unique smooth solution $(u_{k},v_{k+\tilde{\alpha}})$ to \eqref{lp0} with \eqref{b2-1} ${\rm (}$or \eqref{b3-1}${\rm )}$ and \eqref{4-2} satisfying
\begin{align}\label{4-30}
&\sum\limits_{i+j\leq m}\int_{\mathbb{T}}\int_{0}^{\infty}\big[|\partial_{x}^{i}\partial_{z}^{j}(u_{k}-\mathfrak{A}_{k,\infty})|^2+|\partial_{x}^{i}\partial_{z}^{j}v_{k+\tilde{\alpha}}|^2\big]\langle z\rangle^{2\ell}\,{\rm d}z{\rm d}x\nonumber\\
&\leq C(m,\ell)\Big[\sum\limits_{i+j\leq m}\|\langle z\rangle^{\ell+3}\partial_{x}^{i}\partial_{z}^{j}\mathfrak{f}_{k}\|_{L^2(D)}^2+\sum\limits_{i\leq m+1}|\partial_{x}^{i}\mathfrak{R}|_{L^2(\mathbb{T})}^2\Big]
\end{align}
for any integers $m,\ell\geq 0$, where $\mathfrak{A}_{k,\infty}=\lim\limits_{z\to \infty}u_{k}(x,z)$ is a constant.
\end{lemma}

\noindent\textbf{Proof}. For simplicity, we focus on the case of the Robin boundary condition \eqref{b2-1} in \eqref{lp0}. The case of Dirichlet boundary condition \eqref{b3-1} can be dealt with in essentially the same way. We also refer to \cite{FGLT-2023,FGLT-2024} for the Dirichlet case. We divide the proof into four steps.

\textbf{Step 1}. Let $\kappa(z)\in C_{c}^{\infty}[0,\infty)$ satisfy
$$
\kappa(0)=1,\quad \kappa'(0)=0,\quad \int_{0}^{\infty}\kappa(z)\,{\rm d}z=0,
$$
we denote
\begin{align*}
\mathfrak{u}=u_{k}-\mathfrak{R}(x)\kappa(z),\quad \mathfrak{v}=v_{k+\tilde{\alpha}}-v_{k+\tilde{\alpha}}(x,0)+\mathfrak{R}'(x)\int_{0}^{z}\kappa(s)\,{\rm d}s.
\end{align*}

Then, it follows from \eqref{lp0}, \eqref{b2-1} and \eqref{4-2} that $(\mathfrak{u},\mathfrak{v})$ satisfies 
\begin{align}\label{4-3}
	\left\{
	\begin{aligned}
		&u\partial_{x}\mathfrak{u}+\mathfrak{u}\partial_{x}u+\mathfrak{v}\partial_{z}u+v\partial_{z}\mathfrak{u}=\partial_{zz}^{2}\mathfrak{u}+\mathfrak{f},\\
		&\partial_{x}\mathfrak{u}+\partial_{z}\mathfrak{v}=0,\\
		&(\mathfrak{u}-\partial_{z}\mathfrak{u})\vert_{z=0}=0,\quad \mathfrak{v}\vert_{z=0}=0,\quad \lim\limits_{z\to \infty}\partial_{z}\mathfrak{u}=0,
	\end{aligned}
	\right.
\end{align}
where
$$
\mathfrak{f}=\mathfrak{f}_{k}-\kappa(z)\big(u\mathfrak{R}'(x)+\mathfrak{R}(x)\partial_{x}u\big)-v\kappa'(z)\mathfrak{R}+\mathfrak{R}'(x)\Big(\int_{0}^{z}\kappa(s)\,{\rm d}s\Big)\partial_{z}u+\mathfrak{R}(x)\kappa''(z).
$$

\textbf{Step 2}. We consider the following approximate system:
\begin{align}\label{4-4}
	\left\{
	\begin{aligned}
		&u\partial_{x}\mathfrak{u}^{\lambda}+\mathfrak{u}^{\lambda}\partial_{x}u+v\partial_{z}\mathfrak{u}^{\lambda}+\mathfrak{v}^{\lambda}\partial_{z}u=\partial_{zz}^{2}\mathfrak{u}^{\lambda}+\lambda\partial_{xx}\mathfrak{u}^{\lambda}+\mathfrak{f},\\
		&\partial_{x}\mathfrak{u}^{\lambda}+\partial_{z}\mathfrak{v}^{\lambda}=0,\\
		&(\mathfrak{u}^{\lambda}-\partial_{z}\mathfrak{u}^{\lambda})\vert_{z=0}=0,\quad \mathfrak{v}^{\lambda}\vert_{z=0}=0,\quad \lim\limits_{z\to \infty}\partial_{z}\mathfrak{u}^{\lambda}=0.
	\end{aligned}
	\right.
\end{align}
For each fixed $\lambda>0$, existence of a smooth solution to
\eqref{4-4} follows by standard truncation and approximation arguments. We omit
the details and focus on the uniform estimates in $\lambda$.

\textbf{Step 3}. Multiplying $\eqref{4-4}_{1}$ by $\mathfrak{u}^{\lambda}$ and integrating over $(x,z)\in D$, we obtain
\begin{align}\label{4-5}
&\int_{\mathbb{T}}\int_{0}^{\infty}\Big[u\mathfrak{u}^{\lambda}\partial_{x}\mathfrak{u}^{\lambda}+v\mathfrak{u}^{\lambda}\partial_{z}\mathfrak{u}^{\lambda}\Big]\,{\rm d}z{\rm d}x+\int_{\mathbb{T}}\int_{0}^{\infty}\Big[-\Big(\int_{0}^{z}\partial_{x}\mathfrak{u}^{\lambda}(x,s)\,{\rm d}s\Big)\mathfrak{u}^{\lambda}\partial_{z}u+(\mathfrak{u}^{\lambda})^2\partial_{x}u\Big]\,{\rm d}z{\rm d}x\nonumber\\
&\quad =\int_{\mathbb{T}}\int_{0}^{\infty}[\partial_{zz}\mathfrak{u}^{\lambda}+\lambda\partial_{xx}\mathfrak{u}^{\lambda}]\mathfrak{u}^{\lambda}\,{\rm d}z{\rm d}x+\int_{\mathbb{T}}\int_{0}^{\infty}\mathfrak{f}\mathfrak{u}^{\lambda}\,{\rm d}z{\rm d}x.
\end{align}
It follows from the boundary conditions $\eqref{4-4}_{3}$ and integration by parts that
\begin{align}\label{4-6}
	-\int_{\mathbb{T}}\int_{0}^{\infty}[\partial_{zz}\mathfrak{u}^{\lambda}+\lambda\partial_{xx}\mathfrak{u}^{\lambda}]\mathfrak{u}^{\lambda}\,{\rm d}z{\rm d}x=\|\partial_{z}\mathfrak{u}^{\lambda}\|_{L^2}^2+\lambda\|\partial_{x}\mathfrak{u}^{\lambda}\|_{L^2}^2+\int_{\mathbb{T}}|\mathfrak{u}^{\lambda}(x,0)|^2\,{\rm d}x.
\end{align}
We decompose $\mathfrak{u}^{\lambda}$ as follows:
\begin{align*}
	\mathfrak{u}^{\lambda}(x,z)=\mathfrak{u}^{\lambda}(x,0)+\big[\mathfrak{u}^{\lambda}(x,z)-\mathfrak{u}^{\lambda}(x,0)\big]=:\mathfrak{u}^{\lambda}(x,0)+\bar{\mathfrak{u}}^{\lambda}(x,z),
\end{align*}
then $\bar{\mathfrak{u}}^{\lambda}(x,0)=0$, and it follows from the Hardy inequality that
\begin{align}\label{4-6-4}
	\int_{\mathbb{T}}\int_{0}^{\infty}\frac{|\mathfrak{u}^{\lambda}(x,z)|^2}{1+z^2}\,{\rm d}z{\rm d}x&\leq C\int_{\mathbb{T}}\int_{0}^{\infty}\frac{|\mathfrak{u}^{\lambda}(x,0)|^2}{1+z^2}\,{\rm d}z{\rm d}x+C\int_{\mathbb{T}}\int_{0}^{\infty}\frac{|\bar{\mathfrak{u}}^{\lambda}(x,z)|^2}{z^2}\,{\rm d}z{\rm d}x\nonumber\\
&\leq C\Big(\int_{\mathbb{T}}|\mathfrak{u}^{\lambda}(x,0)|^2\,{\rm d}x+\int_{\mathbb{T}}\int_{0}^{\infty}|\partial_{z}\mathfrak{u}^{\lambda}(x,z)|^2\,{\rm d}z{\rm d}x\Big).
\end{align}
Moreover, it follows from Corollary \ref{cor2} and the Sobolev embedding that for $j,k\in \{0,1\}$
\begin{align*}
	&|\partial_{x}^{j}\partial_{z}^{k}(u-\mathfrak{B})\langle z\rangle^{2}|\leq C(j,k)\sigma,\quad |v|\leq C\sigma.
\end{align*}
Then we obtain from \eqref{4-6-4} and the H\"{o}lder inequality that
\begin{align}\label{4-8}
	&-\int_{\mathbb{T}}\int_{0}^{\infty}\Big[u\partial_{x}\mathfrak{u}^{\lambda}+v\partial_{z}\mathfrak{u}^{\lambda}-\Big(\int_{0}^{z}\partial_{x}\mathfrak{u}^{\lambda}(x,s)\,{\rm d}s\Big)\partial_{z}u+\mathfrak{u}^{\lambda}\partial_{x}u\Big]\mathfrak{u}^{\lambda}\,{\rm d}z{\rm d}x+\int_{\mathbb{T}}\int_{0}^{\infty}\mathfrak{f}\mathfrak{u}^{\lambda}\,{\rm d}z{\rm d}x\nonumber\\
	&\leq \int_{\mathbb{T}}\int_{0}^{\infty}\frac{1}{2}\big[\partial_{x}u+\partial_{z}v](\mathfrak{u}^{\lambda})^2\,{\rm d}z{\rm d}x+\|\langle z\rangle^2\partial_{z}u\|_{\infty}\Big\|\frac{\int_{0}^{z}\partial_{x}\mathfrak{u}^{\lambda}(x,s)\,{\rm d}s}{\langle z\rangle}\Big\|_{L^2}\Big\|\frac{\mathfrak{u}^{\lambda}}{\langle z\rangle}\Big\|_{L^2}\nonumber\\
	&\quad +\|\langle z\rangle^2\partial_{x}u\|_{\infty}\Big\|\frac{\mathfrak{u}^{\lambda}}{\langle z\rangle}\Big\|_{L^2}^2+\|\langle z\rangle \mathfrak{f}\|_{L^2}\Big\|\frac{\mathfrak{u}^{\lambda}}{\langle z\rangle}\Big\|_{L^2}\nonumber\\
	&\leq C(\sigma+\tau) \big[\|\partial_{x}\mathfrak{u}^{\lambda}\|_{L^2}^2+\|\partial_{z}\mathfrak{u}^{\lambda}\|_{L^2}^2+\int_{\mathbb{T}}|\mathfrak{u}^{\lambda}(x,0)|^2\,{\rm d}x\Big]+C_{\tau}\|\langle z\rangle \mathfrak{f}\|_{L^2}^2.
\end{align}
Then, combining \eqref{4-5}, \eqref{4-6} and \eqref{4-8} and taking $\sigma$ and $\tau$ small enough, there holds that
\begin{align}\label{4-9}
	\lambda\|\partial_{x}\mathfrak{u}^{\lambda}\|_{L^2}^2+\|\partial_{z}\mathfrak{u}^{\lambda}\|_{L^2}^2+\int_{\mathbb{T}}|\mathfrak{u}^{\lambda}(x,0)|^2\,{\rm d}x\leq C\sigma\|\partial_{x}\mathfrak{u}^{\lambda}\|_{L^2}^2+C\|\langle z\rangle \mathfrak{f}\|_{L^2}^2.
\end{align}

It still remains to estimate $\|\partial_{x}\mathfrak{u}^{\lambda}\|_{L^2}^2$. To this end, we multiply $\eqref{4-4}_{1}$ by $\partial_{x}\mathfrak{u}^{\lambda}$ and integrating over $(x,z)\in D$, we have
\begin{align*}
&\int_{\mathbb{T}}\int_{0}^{\infty}\Big[u\partial_{x}\mathfrak{u}^{\lambda}+v\partial_{z}\mathfrak{u}^{\lambda}-\Big(\int_{0}^{z}\partial_{x}\mathfrak{u}^{\lambda}(x,s)\,{\rm d}s\Big)\partial_{z}u+\mathfrak{u}^{\lambda}\partial_{x}u-\partial_{zz}\mathfrak{u}^{\lambda}-\lambda\partial_{xx}\mathfrak{u}^{\lambda}\Big]\partial_{x}\mathfrak{u}^{\lambda}\,{\rm d}z{\rm d}x\nonumber\\
&=\int_{\mathbb{T}}\int_{0}^{\infty}\mathfrak{f}\partial_{x}\mathfrak{u}^{\lambda}\,{\rm d}z{\rm d}x.
\end{align*}
It is clear that
\begin{align}\label{4.11}
\int_{\mathbb{T}}\int_{0}^{\infty}u|\partial_{x}\mathfrak{u}^{\lambda}|^2\,{\rm d}z{\rm d}x\geq (\mathfrak{B}-C\sigma)\|\partial_{x}\mathfrak{u}^{\lambda}\|_{L^2}^2.	
\end{align}
By integrating by parts, the dissipation term can be computed as follows:
\begin{align}\label{4-12}
&\int_{\mathbb{T}}\int_{0}^{\infty}\big[-\partial_{zz}\mathfrak{u}^{\lambda}-\lambda \partial_{xx}\mathfrak{u}^{\lambda}\big]\partial_{x}\mathfrak{u}^{\lambda}\,{\rm d}z{\rm d}x\nonumber\\
&=\int_{\mathbb{T}}\int_{0}^{\infty}\Big(\frac{1}{2}\partial_{x}(\partial_{z}\mathfrak{u}^{\lambda})^2-\frac{\lambda}{2}\partial_{x}\big(\partial_{x}\mathfrak{u}^{\lambda}\big)^2\Big)\,{\rm d}z{\rm d}x+\int_{\mathbb{T}}\partial_{x}\mathfrak{u}^{\lambda}(x,0)\mathfrak{u}^{\lambda}(x,0)\,{\rm d}x=0.	
\end{align}
Moreover, the remaining terms
\begin{align}\label{4-13}
&\Big\vert\int_{\mathbb{T}}\int_{0}^{\infty}\Big[v\partial_{z}\mathfrak{u}^{\lambda}-\Big(\int_{0}^{z}\partial_{x}\mathfrak{u}^{\lambda}(x,s)\,{\rm d}s\Big)\partial_{z}u+\mathfrak{u}^{\lambda}\partial_{x}u\Big]\partial_{x}\mathfrak{u}^{\lambda}\,{\rm d}z{\rm d}x\Big\vert\nonumber\\
&\quad +\Big\vert\int_{\mathbb{T}}\int_{0}^{\infty}\mathfrak{f}\partial_{x}\mathfrak{u}^{\lambda}\,{\rm d}z{\rm d}x\Big\vert\nonumber\\
&\leq \|v\|_{\infty}\|\partial_{z}\mathfrak{u}^{\lambda}\|_{L^2}\|\partial_{x}\mathfrak{u}^{\lambda}\|_{L^2}+\|z\partial_{z}u\|_{\infty}\Big\|\frac{\int_{0}^{z}\partial_{x}\mathfrak{u}^{\lambda}(x,s)\,{\rm d}s}{z}\Big\|_{L^2}\|\partial_{x}\mathfrak{u}^{\lambda}\|_{L^2}\nonumber\\
&\quad +\|\langle z\rangle\partial_{x}u\|_{\infty}\Big\|\frac{\mathfrak{u}^{\lambda}}{\langle z\rangle}\Big\|_{L^2}\|\partial_{x}\mathfrak{u}^{\lambda}\|_{L^2}+\|\mathfrak{f}\|_{L^2}\|\partial_{x}\mathfrak{u}^{\lambda}\|_{L^2}\nonumber\\
&\leq C\sigma[\|\partial_{x}\mathfrak{u}^{\lambda}\|_{L^2}^2+\|\partial_{z}\mathfrak{u}^{\lambda}\|_{L^2}^2+\int_{\mathbb{T}}|\mathfrak{u}^{\lambda}(x,0)|^2\,{\rm d}x]+\|\mathfrak{f}\|_{L^2}\|\partial_{x}\mathfrak{u}^{\lambda}\|_{L^2}.
\end{align}
Thus, taking $\sigma$ small enough, we obtain from \eqref{4-9}, \eqref{4.11}--\eqref{4-13} that
\begin{align}\label{4-13-0}
\|\partial_{x}\mathfrak{u}^{\lambda}\|_{L^2}^2+\|\partial_{z}\mathfrak{u}^{\lambda}\|_{L^2}^2+\int_{\mathbb{T}}|\mathfrak{u}^{\lambda}(x,0)|^2\,{\rm d}x\leq C\|\langle z\rangle \mathfrak{f}\|_{L^2}^2,
\end{align}
where $C$ is independent of $\lambda$.

\smallskip

Now, for second order derivatives, it follows directly from $\eqref{4-4}_{1}$, the boundary conditions $\eqref{4-4}_{3}$ and \eqref{4-13-0} that
\begin{align*}
	&\|\partial_{zz}\mathfrak{u}^{\lambda}+\lambda \partial_{xx}\mathfrak{u}^{\lambda}\|_{L^2}^2\nonumber\\
    &=\|\partial_{zz}\mathfrak{u}^{\lambda}\|_{L^2}^2+2\lambda\|\partial_{xz}\mathfrak{u}^{\lambda}\|_{L^2}^2+\lambda^2\|\partial_{xx}\mathfrak{u}^{\lambda}\|_{L^2}^2+2\lambda\int_{\mathbb{T}}|\partial_{x}\mathfrak{u}^{\lambda}(x,0)|^2\,{\rm d}x\nonumber\\
	&=\|u\partial_{x}\mathfrak{u}^{\lambda}+v\partial_{z}\mathfrak{u}^{\lambda}-\Big(\int_{0}^z\partial_{x}\mathfrak{u}^{\lambda}(x,s)\,{\rm d}s\Big)\partial_{z}u+\mathfrak{u}^{\lambda}\partial_{x}u-\mathfrak{f}\|_{L^2}^2\nonumber\\
	&\leq C\big[\|\mathfrak{f}\|_{L^2}^2+\|\partial_{x}\mathfrak{u}^{\lambda}\|_{L^2}^2+\|\partial_{z}\mathfrak{u}^{\lambda}\|_{L^2}^2+\int_{\mathbb{T}}|\mathfrak{u}^{\lambda}(x,0)|^2\,{\rm d}x\big]\leq C\|\langle z\rangle \mathfrak{f}\|_{L^2}^2.
\end{align*}
Thus, we obtain
\begin{align}\label{4-13-1}
	\|\partial_{x}\mathfrak{u}^{\lambda}\|_{L^2}^2+\|\partial_{z}\mathfrak{u}^{\lambda}\|_{L^2}^2+\|\partial_{zz}\mathfrak{u}^{\lambda}\|_{L^2}^2\leq C\|\langle z\rangle \mathfrak{f}\|_{L^2}^2.
\end{align}

The above inequality shows the uniqueness of solution to the system \eqref{4-4} for any fixed $\lambda>0$. Moreover, the solution is smooth if $\mathfrak{f}$ is smooth enough. Set
\begin{align}\label{uv}
\mathfrak{u}=\lim\limits_{\lambda\to 0}\mathfrak{u}^{\lambda},\quad \mathfrak{v}=-\int_{0}^{z}\partial_{x}\mathfrak{u}(x,s)\,{\rm d}s,
\end{align}
then, taking the limit $\lambda\to 0$ in $\eqref{4-13-1}$ to get
\begin{align}\label{4-13-2}
\|\partial_{x}\mathfrak{u}\|_{L^2}^2+\|\partial_{z}\mathfrak{u}\|_{L^2}^2+\|\partial_{zz}\mathfrak{u}\|_{L^2}^2\leq C\|\langle z\rangle \mathfrak{f}\|_{L^2}^2.
\end{align}
Moreover, $(\mathfrak{u},\mathfrak{v})$ solves the system \eqref{4-3}.

\textbf{Step 4}. Finally, we show that the derivatives of $\mathfrak{u}$ and $\mathfrak{v}$ decay rapidly as $z\to \infty$. To this end, we introduce the stream function $\psi=\int_{0}^{z}u(x,s)\,{\rm d}s$, and denote
\begin{align*}
	\tilde{\mathfrak{u}}(x,\psi)=\mathfrak{u}(x,z(x,\psi)),\quad \tilde{\mathfrak{v}}(x,\psi)=\mathfrak{v}(x,z(x,\psi)),\quad \mathfrak{F}(x,\psi)=\frac{\mathfrak{f}(x,z(x,\psi))}{u(x,z(x,\psi))}.
\end{align*}
Then, it follows from \eqref{4-3} that $(\tilde{\mathfrak{u}},\tilde{\mathfrak{v}})$ satisfies
\begin{align}\label{4-14}
	\left\{
	\begin{aligned}
	&\partial_{x}\tilde{\mathfrak{u}}-\partial_{\psi}(a(x,\psi)\partial_{\psi}\tilde{\mathfrak{u}})+b(x,\psi)\tilde{\mathfrak{u}}+c(x,\psi)\tilde{\mathfrak{v}}=\mathfrak{F}(x,\psi),\\
	&(\tilde{\mathfrak{u}}-a(x,\psi)\partial_{\psi}\tilde{\mathfrak{u}})\vert_{\psi=0}=0,
	\end{aligned}
	\right.
\end{align}
where
$$
a(x,\psi)=u(x,z(x,\psi)),\quad b(x,\psi)=\frac{\partial_{x}u(x,z(x,\psi))}{u(x,z(x,\psi))},\quad c(x,\psi)=\frac{\partial_{z}u(x,z(x,\psi))}{u(x,z(x,\psi))}.
$$
Noticing from Corollary \ref{cor2} that
$$
\frac{\mathfrak{B}}{2}\leq u(x,z)\leq \frac{3\mathfrak{B}}{2},
$$
we deduce that $\frac{\mathfrak{B}}{2}\leq \frac{\psi}{z}\leq \frac{3\mathfrak{B}}{2}$.

We claim that for any $\ell\in \mathbb{N}$, there holds
\begin{align}\label{4-16}
	\|\partial_{x}\tilde{\mathfrak{u}}\psi^{\ell}\|_{L^2}+\|\partial_{\psi}\tilde{\mathfrak{u}}_{\neq}\psi^{\ell}\|_{L^2}+\|\partial_{\psi\psi}\tilde{\mathfrak{u}}\psi^{\ell}\|_{L^2}\leq C_{\ell}\|\mathfrak{F}\langle\psi\rangle^{\ell+1}\|_{L^2}.
\end{align}
where
$
\tilde{\mathfrak{u}}_{\neq}=\tilde{\mathfrak{u}}-\overline{\tilde{\mathfrak{u}}}(\psi).
$
In fact, from \eqref{4-13-2}, it holds that
\begin{align}\label{4-16-0}
\|\partial_{x}\tilde{\mathfrak{u}}\|_{L^2}+\|\partial_{\psi}\tilde{\mathfrak{u}}\|_{L^2}+\|\partial_{\psi\psi}\tilde{\mathfrak{u}}\|_{L^2}\leq C\|\mathfrak{F}\langle \psi\rangle\|_{L^2},
\end{align}
which implies \eqref{4-16} holds for $\ell=0$.

For any $\ell\geq 1$, multiplying $\tilde{\mathfrak{u}}_{\neq}\psi^{2\ell}$ and integrating $(x,\psi)\in \mathbb{T}\times (0,\infty)$, we obtain
\begin{align}\label{4-16-1}
	&\int_{\mathbb{T}}\int_{0}^{\infty}\partial_{x}\tilde{\mathfrak{u}}\tilde{\mathfrak{u}}_{\neq}\psi^{2\ell}\,{\rm d}\psi{\rm d}x-\int_{\mathbb{T}}\int_{0}^{\infty}\partial_{\psi}(a(x,\psi)\partial_{\psi}\tilde{\mathfrak{u}})\tilde{\mathfrak{u}}_{\neq}\psi^{2\ell}\,{\rm d}\psi{\rm d}x\nonumber\\
	&=\int_{\mathbb{T}}\int_{0}^{\infty}[\mathfrak{F}(x,\psi)-b(x,\psi)\tilde{\mathfrak{u}}-c(x,\psi)\tilde{\mathfrak{v}}]\tilde{\mathfrak{u}}_{\neq}\psi^{2\ell}\,{\rm d}\psi{\rm d}x.
\end{align}
For the first term on the left hand side of \eqref{4-16-1}, it is clear that
$$
\int_{\mathbb{T}}\int_{0}^{\infty}\partial_{x}\tilde{\mathfrak{u}}\tilde{\mathfrak{u}}_{\neq}\psi^{2\ell}\,{\rm d}\psi{\rm d}x=0.
$$
For the second term on the left hand side of \eqref{4-16-1}, using the integration by parts, one has
\begin{align}\label{4-17}
	&-\int_{\mathbb{T}}\int_{0}^{\infty}\partial_{\psi}(a(x,\psi)\partial_{\psi}\tilde{\mathfrak{u}})\tilde{\mathfrak{u}}_{\neq}\psi^{2\ell}\,{\rm d}\psi{\rm d}x\nonumber\\
    &=\int_{\mathbb{T}}\int_{0}^{\infty}a(x,\psi)\partial_{\psi}\tilde{\mathfrak{u}}\partial_{\psi}\tilde{\mathfrak{u}}_{\neq}\psi^{2\ell}\,{\rm d}\psi{\rm d}x+2\ell\int_{\mathbb{T}}\int_{0}^{\infty}a(x,\psi)\partial_{\psi}\tilde{\mathfrak{u}}\tilde{\mathfrak{u}}_{\neq}\psi^{2\ell-1}\,{\rm d}\psi{\rm d}x=:I_{1}+I_{2}.
\end{align}
For $I_{1}$, noting that $a(x,\psi)=\mathfrak{B}+u_{B}^{0}(x,\psi)$, we obtain from $|u_{B}^{0}\psi^{2\ell}|\leq C\sigma$ that
\begin{align}\label{4-17-1}
	I_{1}\geq \frac{\mathfrak{B}}{2}\|\partial_{\psi}\tilde{\mathfrak{u}}_{\neq}\psi^{\ell}\|_{L^2}^2-C\sigma\|\partial_{\psi}\tilde{\mathfrak{u}}\|_{L^2}^2.
\end{align}
For $I_{2}$, using the fast decay of $u_{B}^{0}(x,\psi)$, the H\"{o}lder and Poincar\'{e} inequalities, there holds
\begin{align}\label{4-17-2}
|I_{2}|\leq C_{\ell}\|\partial_{\psi}\tilde{\mathfrak{u}}_{\neq}\psi^{\ell}\|_{L^2}\|\partial_{x}\tilde{\mathfrak{u}}\psi^{\ell-1}\|_{L^2}+C_{\ell}\sigma\|\partial_{\psi}\tilde{\mathfrak{u}}\|_{L^2}\|\partial_{x}\tilde{\mathfrak{u}}\|_{L^2}.
\end{align}
For the terms on the right hand side of \eqref{4-16-1}, due to the fast decay of $b(x,\psi)$, $c(x,\psi)$ as $\psi\to \infty$, it can be bounded as
\begin{align}\label{4-18}
	&C\|\mathfrak{F}\langle \psi\rangle^{\ell+1}\|_{L^2}\|\partial_{x}\tilde{\mathfrak{u}}_{\neq}\psi^{\ell-1}\|_{L^2}+C\sigma(\|\partial_{x}\tilde{\mathfrak{u}}\|_{L^2}^2+\|\partial_{\psi}\tilde{\mathfrak{u}}\|_{L^2}^2+\int_{\mathbb{T}}|\mathfrak{u}(x,0)|^2\,{\rm d}x)\nonumber\\
	&\leq C\|\mathfrak{F}\langle \psi\rangle^{\ell+1}\|_{L^2}^2+C\|\partial_{x}\tilde{\mathfrak{u}}\psi^{\ell-1}\|_{L^2}^2.
\end{align}
Thus, substituting \eqref{4-17}--\eqref{4-18} into \eqref{4-16-1} gives
\begin{align}\label{4-19}
\|\partial_{\psi}\tilde{\mathfrak{u}}_{\neq}\psi^{\ell}\|_{L^2}\leq C\|\partial_{x}\tilde{\mathfrak{u}}\psi^{\ell-1}\|_{L^2}+C\|\mathfrak{F}\langle \psi\rangle^{\ell+1}\|_{L^2}.
\end{align}

Moreover, we rewrite $\eqref{4-14}_{1}$ as
\begin{align*}
\partial_{x}\tilde{\mathfrak{u}}	-a(x,\psi)\partial_{\psi\psi}\tilde{\mathfrak{u}}=\mathfrak{F}(x,\psi)+\partial_{\psi}a(x,\psi)\partial_{\psi}\tilde{\mathfrak{u}}-b(x,\psi)\tilde{\mathfrak{u}}-c(x,\psi)\tilde{\mathfrak{v}}.
\end{align*}
Thus, there holds
\begin{align*}
	\|[\partial_{x}\tilde{\mathfrak{u}}-a(x,\psi)\partial_{\psi\psi}\tilde{\mathfrak{u}}]\psi^{\ell}\|_{L^2}^2=\|[\mathfrak{F}(x,\psi)+\partial_{\psi}a(x,\psi)\partial_{\psi}\tilde{\mathfrak{u}}-b(x,\psi)\tilde{\mathfrak{u}}-c(x,\psi)\tilde{\mathfrak{v}}]\psi^{\ell}\|_{L^2}^2.
\end{align*}
Due to the fast decay of $\partial_{\psi}a(x,\psi)$, $b(x,\psi)$ and $c(x,\psi)$, the terms on the right hand side above can be controlled by
\begin{align*}
	C\|\mathfrak{F}\langle \psi\rangle^{\ell}\|_{L^2}^2+C\big(\|\partial_{x}\tilde{\mathfrak{u}}\|_{L^2}^2+\|\partial_{\psi}\tilde{\mathfrak{u}}\|_{L^2}^2+\int_{\mathbb{T}}|\tilde{\mathfrak{u}}(x,0)|^2\,{\rm d}x\big)\leq C\|\mathfrak{F}\langle \psi\rangle^{\ell}\|_{L^2}^2,
\end{align*}
where we have used \eqref{4-16-0}. Moreover, a direct calculation shows that
\begin{align*}
&\|[\partial_{x}\tilde{\mathfrak{u}}-a(x,\psi)\partial_{\psi\psi}\tilde{\mathfrak{u}}]\psi^{\ell}\|_{L^2}^2\\
&=\|\partial_{x}\tilde{\mathfrak{u}} \psi^{\ell}\|_{L^2}^2+\|a(x,\psi)\partial_{\psi\psi}\tilde{\mathfrak{u}}\psi^{\ell}\|_{L^2}^2-2\int_{\mathbb{T}}\int_{0}^{\infty}a(x,\psi)\partial_{x}\tilde{\mathfrak{u}}\partial_{\psi\psi}\tilde{\mathfrak{u}}\psi^{2\ell}\,{\rm d}\psi{\rm d}x.
\end{align*}
For the last term in the above equation, it follows from the integration by parts that
\begin{align*}
	-2\int_{\mathbb{T}}\int_{0}^{\infty}a(x,\psi)\partial_{x}\tilde{\mathfrak{u}}\partial_{\psi\psi}\tilde{\mathfrak{u}}\psi^{2\ell}\,{\rm d}\psi{\rm d}x&=2\int_{\mathbb{T}}\int_{0}^{\infty}\partial_{\psi}a(x,\psi)\partial_{x}\tilde{\mathfrak{u}}\partial_{\psi}\tilde{\mathfrak{u}}\psi^{2\ell}\,{\rm d}\psi{\rm d}x\nonumber\\
	&\quad +4\ell\int_{\mathbb{T}}\int_{0}^{\infty}a(x,\psi)\partial_{x}\tilde{\mathfrak{u}}\partial_{\psi}\tilde{\mathfrak{u}}\psi^{2\ell-1}\,{\rm d}\psi{\rm d}x\nonumber\\
	&\quad +2\int_{\mathbb{T}}\int_{0}^{\infty}a(x,\psi)\partial_{x\psi}\tilde{\mathfrak{u}}\partial_{\psi}\tilde{\mathfrak{u}}\psi^{2\ell}\,{\rm d}\psi{\rm d}x\nonumber\\
    &=J_{1}+J_{2}+J_{3}.
\end{align*}
For $J_{1}$, using $|\partial_{\psi} a(x,\psi)|\leq C\sigma$ and $\partial_{\psi}a$ has fast decay as $\psi\to \infty$, one has
$$
|J_{1}|\leq C\sigma(\|\partial_{x}\tilde{\mathfrak{u}}\|_{L^2}^2+\|\partial_{\psi}\tilde{\mathfrak{u}}\|_{L^2}^2),
$$
For $J_{3}$, using integration by parts, and the facts $|\partial_{x} a(x,\psi)|\leq C\sigma$ and $\partial_{x}a$ has fast decay as $\psi\to \infty$, we have
$$
|J_{3}|\leq C\sigma(\|\partial_{x}\tilde{\mathfrak{u}}\|_{L^2}^2+\|\partial_{\psi}\tilde{\mathfrak{u}}\|_{L^2}^2).
$$
Using the fact $a(x,\psi)-\mathfrak{B}=O(\sigma)$ has fast decay as $\psi\to \infty$, one has
\begin{align*}
&\|a(x,\psi)\partial_{\psi\psi}\tilde{\mathfrak{u}}\psi^{\ell}\|_{L^2}^2\geq \frac{|\mathfrak{B}|^2}{4}\|\partial_{\psi\psi}\tilde{\mathfrak{u}}\psi^{\ell}\|_{L^2}^2,\\
&|J_{2}|
	\leq C_{\ell}\|\partial_{x}\tilde{\mathfrak{u}}\psi^{\ell-1}\|_{L^2}\|\partial_{\psi}\tilde{\mathfrak{u}}_{\neq}\psi^{\ell}\|_{L^2}+C_{\ell}\sigma(\|\partial_{x}\tilde{\mathfrak{u}}\|_{L^2}^2+\|\partial_{\psi}\tilde{\mathfrak{u}}\|_{L^2}^2).
\end{align*}
Thus, we obtain from \eqref{4-19} that
\begin{align*}
	\|\partial_{x}\tilde{\mathfrak{u}}\psi^{\ell}\|_{L^2}+\|\partial_{\psi\psi}\tilde{\mathfrak{u}}\psi^{\ell}\|_{L^2}\leq C_{\ell}\|\partial_{x}\tilde{\mathfrak{u}}\psi^{\ell-1}\|_{L^2}+C_{\ell}\|\mathfrak{F}\langle \psi\rangle^{\ell+1}\|_{L^2},
\end{align*}
which, together with \eqref{4-16-0} and \eqref{4-19}, yields
\begin{align*}
	\|\partial_{x}\tilde{\mathfrak{u}}\psi^{\ell}\|_{L^2}+\|\partial_{\psi}\tilde{\mathfrak{u}}_{\neq}\psi^{\ell}\|_{L^2}+\|\partial_{\psi\psi}\tilde{\mathfrak{u}}\psi^{\ell}\|_{L^2}\leq C_{\ell}\|\partial_{x}\tilde{\mathfrak{u}}\psi^{\ell-1}\|_{L^2}+C_{\ell}\|\mathfrak{F}\langle \psi\rangle^{\ell+1}\|_{L^2}.
\end{align*}
Hence, by induction, it holds that
\begin{align*}
	\|\partial_{x}\tilde{\mathfrak{u}}\psi^{\ell}\|_{L^2}+\|\partial_{\psi}\tilde{\mathfrak{u}}_{\neq}\psi^{\ell}\|_{L^2}+\|\partial_{\psi\psi}\tilde{\mathfrak{u}}\psi^{\ell}\|_{L^2}\leq C_{\ell}\|\mathfrak{F}\langle\psi\rangle^{\ell+1}\|_{L^2},
\end{align*}
which concludes \eqref{4-16}.

Returning to the original variables $(x,z)$, we obtain that for any $\ell\in \mathbb{N}$,
\begin{align*}
	\|\langle z\rangle^{\ell}\partial_{zz}\mathfrak{u}\|_{L^2}^2+\|\langle z\rangle^{\ell-1}\partial_{x}\mathfrak{u}\|_{L^2}^2+\|\langle z\rangle^{\ell-1}\partial_{z}\mathfrak{u}\|_{L^2}^2\leq C_{\ell}\|\langle z\rangle^{\ell+1}\mathfrak{f}\|_{L^2}^2.
\end{align*}
Furthermore, by using the induction in $(x,\psi)$ coordinates and then returning to $(x,z)$ coordinates, we obtain that for any integers $m,\ell\geq 0$,
\begin{align}\label{4.25}
	&\sum\limits_{j+k\leq m}\Big(\|\langle z\rangle^{\ell}\partial_{x}^{j}\partial_{z}^{k}\partial_{zz}\mathfrak{u}\|_{L^2}^2+\|\langle z\rangle^{\ell-1}\partial_{x}^{j}\partial_{z}^{k}\partial_{x}\mathfrak{u}\|_{L^2}^2+\|\langle z\rangle^{\ell-1}\partial_{x}^{j}\partial_{z}^{k}\partial_{z}\mathfrak{u}\|_{L^2}^2\Big)\nonumber\\
	&\leq C(m,\ell)\sum\limits_{j+k\leq m}\|\langle z\rangle^{\ell+1}\partial_{x}^{j}\partial_{z}^{k}\mathfrak{f}\|_{L^2}^2.
\end{align}

Noting that $\lim\limits_{z\to \infty}(\mathfrak{u}_{x},\mathfrak{u}_{z})=0$ and $\mathfrak{A}_{\infty}=\lim\limits_{z\to \infty}\mathfrak{u}(x,z)$ is a constant independent of $x$, then by the Hardy inequality, we have for any $\ell\geq 2$,
\begin{align*}
	\|z^{\ell-2}(\mathfrak{u}-\mathfrak{A}_{\infty})\|_{L^2}\leq C_{\ell}\|z^{\ell-1}\partial_{z}\mathfrak{u}\|_{L^2}.
\end{align*}
The weighted estimates on $\mathfrak{v}-\mathfrak{v}_{\infty}:=\mathfrak{v}+\int_{0}^{\infty}\partial_{x}\mathfrak{u}(x,s)\,{\rm d}s$ follow directly from \eqref{4.25} and \eqref{uv}. Finally, noting that
$$
u_{k}=\mathfrak{u}+\mathfrak{R}(x)\kappa(z),\quad v_{k+\tilde{\alpha}}=\mathfrak{v}-\mathfrak{R}'(x)\int_{0}^{z}\kappa(s)\,{\rm d}s+\int_{0}^{\infty}\mathfrak{u}_{x}(x,s)\,{\rm d}s,
$$
the weighted estimates for $(u_{k},v_{k+\tilde{\alpha}})$ \eqref{4-30} follows from \eqref{4.25} and the estimates on $\mathfrak{v}-\mathfrak{v}_{\infty}$. Therefore, the proof of Lemma \ref{lem3} is complete. $\hfill\square$

\medskip

\subsubsection{Subcritical regime: $0\leq \alpha<1$} For the system \eqref{lp0} with \eqref{b1-1}, recall the Lemma \ref{lem2}, the leading boundary-layer profile $(u,v)\equiv (\mathfrak{B},0)$, and \eqref{lp0} reduces to
\begin{align}\label{lp1}
\left\{\begin{aligned}
&\mathfrak{B}\partial_{x}u_{k}=\partial_{zz}u_{k}+\mathfrak{f}_{k},\\
&\partial_{x}u_{k}+\partial_{z}v_{k+\tilde{\alpha}}=0,\quad \partial_{z}u_{k}\vert_{z=0}=\mathfrak{R}_{3}(x).
\end{aligned}
\right.
\end{align}
Noting that $u_{k}$ enters the equation \eqref{lp1} as its derivatives, so a change in the far-field constant of $u_{k}$ does not affect \eqref{lp1}. Therefore, we can impose the following far-field condition
\begin{align}\label{f2}
\lim\limits_{z\to \infty}(u_{k},v_{k+\tilde{\alpha}})=(0,0).
\end{align}
It follows from $\eqref{lp1}$ and \eqref{f2} that
$$
v_{k+\tilde{\alpha}}(x,z)=\int_{z}^{\infty}\partial_{x}u_{k}(x,s)\,{\rm d}s.
$$
Before showing the well-posedness of system \eqref{lp1} with \eqref{f2}, we first consider its solvability condition. Actually, integrating \eqref{lp1} with respect to $(x,z)\in \mathbb{T}\times (0,\infty)$, we have
\begin{align}\label{slp1}
\int_{\mathbb{T}}\int_{0}^{\infty}\mathfrak{f}_{k}\,{\rm d}z{\rm d}x=\int_{\mathbb{T}}\partial_{z}u_{k}(x,0)\,{\rm d}x=\int_{\mathbb{T}}\mathfrak{R}_{3}(x)\,{\rm d}x.
\end{align}
Applying $\partial_{z}$ to \eqref{lp1} denoting $\mathfrak{q}_{k}=\partial_{z}u_{k}$, one has
\begin{align*}
\left\{
\begin{aligned}   &\mathfrak{B}\partial_{x}\mathfrak{q}_{k}=\partial_{zz}\mathfrak{q}_{k}+\partial_{z}\mathfrak{f}_{k},\\
&\mathfrak{q}_{k}\vert_{z=0}=\mathfrak{R}_{3}(x),\quad \lim\limits_{z\to \infty}\mathfrak{q}_{k}=0.
\end{aligned}
\right.
\end{align*}
It reduces to the case of \eqref{lp0} with \eqref{b3-1} and $(u,v)=(\mathfrak{B},0)$.
Then, under the solvability condition \eqref{slp1}, using similar and simpler calculations as in Lemma \ref{lem3}, we have following lemma, see also \cite{PZ-2025}.

\begin{lemma}\label{lem4}
Let $\mathfrak{B}>0$. Let $\mathfrak{R}_{3}$ and $\mathfrak{f}_{k}$ be smooth, periodic in $x$ satisfying
\begin{align*}
\sum\limits_{i+j\leq m}\big[\|\langle z\rangle^{\ell}\partial_{x}^{i}\partial_{z}^{j}\mathfrak{f}_{k}\|_{L^2(D)}+\|\partial_{x}^{i}\mathfrak{R}_{3}\|_{L^2(\mathbb{T})}\big]<\infty,\quad \text{for any integers $m,\ell\geq 0$}.
\end{align*}
Assume \eqref{slp1}, there exists a unique smooth solution $(u_{k},v_{k+\tilde{\alpha}})$ to \eqref{lp1} with \eqref{f2} satisfying,
\begin{align*}
&\sum\limits_{i+j\leq m}\int_{\mathbb{T}}\int_{0}^{\infty}\big[|\partial_{x}^{i}\partial_{z}^{j}u_{k}|^2+|\partial_{x}^{i}\partial_{z}^{j}v_{k+\tilde{\alpha}}|^2\big]\langle z\rangle^{2\ell}\,{\rm d}z{\rm d}x\nonumber\\
&\leq C(m,\ell)\Big[\sum\limits_{i+j\leq m}\|\langle z\rangle^{\ell+3}\partial_{x}^{i}\partial_{z}^{j}\mathfrak{f}_{k}\|_{L^2(D)}^2+\sum\limits_{i\leq m+1}|\partial_{x}^{i}\mathfrak{R}_{3}|_{L^2(\mathbb{T})}^2\Big]
\end{align*}
 for any integers $m,\ell\geq 0$.
\end{lemma}

\section{Constructions of Approximate Solutions}

In this section, we rigorously construct the approximate solution via the Prandtl expansions, and show the interplay between the boundary-layer and interior profiles. Let $\chi(s)$ be a smooth and monotonously decreasing cut-off function satisfying \eqref{cut-off}
and $\tilde{\chi}(y)=1-\chi(y)$. Then we construct our approximate solutions as follows:
\begin{align}
	u^{a}&=\sum\limits_{i=0}^{k_{0}}\varepsilon^{\frac{i}{\tilde{\alpha}}}u_{I}^{i}(x,y)+\sum\limits_{i=0}^{k_{0}}\varepsilon^{\frac{i}{\tilde{\alpha}}}\big(\chi(y)u_{B}^{i}(x,z)+\tilde{\chi}(y)\hat{u}_{B}^{i}(x,\zeta)\big)+\varepsilon^{\frac{k_{0}}{\tilde{\alpha}}}h_{u}(x,y)\nonumber\\
    &=:u_{I}^{a}+u_{B}^{a}+\hat{u}_{B}^{a}+\varepsilon^{\frac{k_{0}}{\tilde{\alpha}}}h_{u},\label{3.1}\\
	v^{a}&=\sum\limits_{i=0}^{k_{0}}\varepsilon^{\frac{i}{\tilde{\alpha}}}v_{I}^{i}(x,y)+\sum\limits_{i=0}^{k_{0}+\tilde{\alpha}}\varepsilon^{\frac{i}{\tilde{\alpha}}}\big(\chi(y)v_{B}^{i}(x,z)+\tilde{\chi}(y)\hat{v}_{B}^{i}(x,\zeta)\big)+\varepsilon^{\frac{k_{0}}{\tilde{\alpha}}}h_{v}(x,y)\nonumber\\
    &=:v_{I}^{a}+v_{B}^{a}+\hat{v}_{B}^{a}+\varepsilon^{\frac{k_{0}}{\tilde{\alpha}}}h_{v},\label{3.2}\\
	p^{a}&=\sum\limits_{i=0}^{k_{0}}\varepsilon^{\frac{i}{\tilde{\alpha}}}p_{I}^{i}(x,y)+\sum\limits_{i=0}^{k_{0}+\tilde{\alpha}}\varepsilon^{\frac{i}{\tilde{\alpha}}}\big(\chi(y)p_{B}^{i}(x,z)+\tilde{\chi}(y)\hat{p}_{B}^{i}(x,\zeta)\big),\label{3.3}
\end{align}
where $z=\frac{y}{\varepsilon}$, $\zeta=\frac{y-1}{\varepsilon}$ and $k_{0}$ is an integer truncation index. Here, $(u_{I}^{a},v_{I}^{a})$, $(u_{B}^{a},v_{B}^{a})$ and $(\hat{u}_{B}^{a},\hat{v}_{B}^{a})$ denote the interior, lower and upper boundary-layer parts in approximate solutions $(u^{a},v^{a})$, respectively. $k_{0}$ is a constant which will be chosen later, $h_{u}(x,y)$ and $h_{v}(x,y)$ are two artificial auxiliary functions to ensure that the approximate solution $(u^{a},v^{a})$ will satisfy
\begin{align}\label{3.4}
	\begin{aligned}
	&\partial_{x}u^{a}+\partial_{y}v^{a}=0,\,\, v^{a}(x,0)=v^{a}(x,1)=0,\\ &(u^{a}-\varepsilon^{\alpha}\partial_{y}u^{a})\vert_{y=0}=f_{0}(x),\,\, (u^{a}+\varepsilon^{\alpha}\partial_{y}u^{a})\vert_{y=1}=f_{1}(x).
	\end{aligned}
\end{align}
A direct calculation shows that $(h_{u},h_{v})$ satisfies

\begin{align}\label{3.5-0}
\left\{\begin{aligned}
&\partial_{x}h_{u}+\partial_{y}h_{v}=-\varepsilon^{-\frac{k_{0}}{\tilde{\alpha}}}\chi'(y)\sum\limits_{i=\tilde{\alpha}}^{\tilde{\alpha}+k_{0}}\varepsilon^{\frac{i}{\tilde{\alpha}}}(v_{B}^{i}-\hat{v}_{B}^{i}),\\ 
&h_{v}(x,y)\vert_{y=0}=-\varepsilon^{-\frac{k_{0}}{\tilde{\alpha}}}\sum\limits_{i=k_{0}+1}^{\tilde{\alpha}+k_{0}}\varepsilon^{\frac{i}{\tilde{\alpha}}}v_{B}^{i}(x,0),\,\,h_{v}(x,y)\vert_{y=1}=-\varepsilon^{-\frac{k_{0}}{\tilde{\alpha}}}\sum\limits_{i=k_{0}+1}^{\tilde{\alpha}+k_{0}}\varepsilon^{\frac{i}{\tilde{\alpha}}}\hat{v}_{B}^{i}(x,0),
\end{aligned}
\right.
\end{align}
and the following boundary conditions on $h_{u}$:
\begin{align}\label{3.5-2}
	\left\{
	\begin{aligned}
&\big(h_{u}-\varepsilon^{\alpha}\partial_{y}h_{u}\big)\vert_{y=0}=\varepsilon\partial_{y}u_{I}^{k_{0}}(x,0),\\
&\big(h_{u}+\varepsilon^{\alpha}\partial_{y}h_{u}\big)\vert_{y=1}=-\varepsilon\partial_{y}u_{I}^{k_{0}}(x,1),
\end{aligned}
\right.\qquad\text{for $\alpha=1$},
\end{align}
and
\begin{flalign}\label{3.5-1}
\left\{\begin{aligned}
\big(h_{u}-\varepsilon^{\alpha}\partial_{y}h_{u}\big)\vert_{y=0}=&-u_{I}^{k_{0}}(x,0)-u_{B}^{k_{0}}(x,0)\\
&+\varepsilon^{\alpha-\frac{k_{0}}{\tilde{\alpha}}}\sum\limits_{i=k_{0}-\tilde{\alpha}+1}^{k_{0}}\varepsilon^{\frac{i}{\tilde{\alpha}}}\partial_{y}u_{I}^{i}(x,0),\\
\big(h_{u}+\varepsilon^{\alpha}\partial_{y}h_{u}\big)\vert_{y=1}=&-u_{I}^{k_{0}}(x,1)-\hat{u}_{B}^{k_{0}}(x,0)\\
& -\varepsilon^{\alpha-\frac{k_{0}}{\tilde{\alpha}}}\sum\limits_{i=k_{0}-\tilde{\alpha}+1}^{k_{0}}\varepsilon^{\frac{i}{\tilde{\alpha}}}\partial_{y}u_{I}^{i}(x,1),
\end{aligned}
\right.\quad \text{for $0\leq \alpha<1$},&&
\end{flalign}
and
\begin{flalign}\label{3.5-3}
\left\{
\begin{aligned}
    \big(h_{u}-\varepsilon^{\alpha}\partial_{y}h_{u}\big)\vert_{y=0}=&\varepsilon^{\alpha-\frac{k_{0}}{\tilde{\alpha}}}\sum\limits_{i=k_{0}-(\alpha-1)\tilde{\alpha}+1}^{k_{0}}\varepsilon^{\frac{i}{\tilde{\alpha}}}\partial_{y}\big[\chi(y)u_{B}^{i}(x,z)+\tilde{\chi}(y)\hat{u}_{B}^{i}(x,\zeta)\big]\big\vert_{y=0}\\
    &+\varepsilon^{\alpha-\frac{k_{0}}{\tilde{\alpha}}}\sum\limits_{i=k_{0}-\alpha\tilde{\alpha}+1}^{k_{0}}\varepsilon^{\frac{i}{\tilde{\alpha}}}\partial_{y}u_{I}^{i}(x,0),\\
    \big(h_{u}+\varepsilon^{\alpha}\partial_{y}h_{u}\big)\vert_{y=1}=&-\varepsilon^{\alpha-\frac{k_{0}}{\tilde{\alpha}}}\sum\limits_{i=k_{0}-(\alpha-1)\tilde{\alpha}+1}^{k_{0}}\varepsilon^{\frac{i}{\tilde{\alpha}}}\partial_{y}\big[\chi(y)u_{B}^{i}(x,z)+\tilde{\chi}(y)\hat{u}_{B}^{i}(x,\zeta)\big]\big\vert_{y=1}\\
    &-\varepsilon^{\alpha-\frac{k_{0}}{\tilde{\alpha}}}\sum\limits_{i=k_{0}-\alpha\tilde{\alpha}+1}^{k_{0}}\varepsilon^{\frac{i}{\tilde{\alpha}}}\partial_{y}u_{I}^{i}(x,1),
\end{aligned}
\right.&&
\end{flalign}
for $\alpha>1$. 

\begin{proposition}\label{prop3}
    There exists an approximate solution $(u^{a},v^{a},p^{a})$ in the form of \eqref{3.1}--\eqref{3.3} such that the following holds: 
       
    \noindent (\rmnum{1}) $(u_{I}^{0},v_{I}^{0},p_{I}^{0})=(Ay+B,0,0)$, where the constants $A$ and $B$ are chosen as follows:

    (a) for $\alpha=0$,
    \begin{align}\label{M00}
    A+B=\frac{2}{3}\bar{f}_{1}+\frac{1}{3}\bar{f}_{0},\quad B=\frac{2}{3}\bar{f}_{0}+\frac{1}{3}\bar{f}_{1}.
    \end{align}
 
    (b) for $0<\alpha<1$,
          \begin{align}\label{M0}
          A+B=\bar{f}_{1},\quad B=\bar{f}_{0}.
          \end{align}
          
    (c) for $\alpha=1$, $B>0$, $(A+B)>0$ and
          \begin{align}\label{M2}
          \begin{aligned}
          \bar{f}_{1}-(A+B)=-\frac{\sigma^2}{4\pi(A+B)}\int_{\mathbb{T}}\hat{\phi}^2(\xi,0)\,{\rm d}\xi,\quad \bar{f}_{0}-B=-\frac{\sigma^2}{4\pi B}\int_{\mathbb{T}}\phi^2(\xi,0)\,{\rm d}\xi,
         \end{aligned}
          \end{align}
          where $\phi(\xi,\eta)$ and $\hat{\phi}(\xi,\eta)$ are the solutions of following systems respectively
          \begin{align}\label{lowerboundary}
\left\{
\begin{aligned}
&\partial_{\xi}\phi=B\partial_{\eta\eta}\phi+\frac{\sigma}{2}\partial_{\eta\eta}(\phi^2),\\
&\Pi_{\perp}[\big(\phi-B\partial_{\eta}\phi-\frac{\sigma}{2}\partial_{\eta}(\phi^2)\big)\vert_{\eta=0}]=g_{0}(\xi),\\
&\lim\limits_{\eta\to \infty}\phi(\xi,\eta)=0,
\end{aligned}
\right.\qquad\text{for $(\xi,\eta)\in \mathbb{T}\times (0,\infty)$},
          \end{align}
and
\begin{align}\label{upperboundary}
\left\{
\begin{aligned}
&\partial_{\xi}\hat{\phi}=(A+B)\partial_{\eta\eta}\hat{\phi}+\frac{\sigma}{2}\partial_{\eta\eta}(\hat{\phi}^2)\\
&\Pi_{\perp}[\big(\hat{\phi}+(A+B)\partial_{\eta}\hat{\phi}+\frac{\sigma}{2}\partial_{\eta}(\hat{\phi}^2))\vert_{\eta=0}]=g_{1}(\xi),\\
&\lim\limits_{\eta\to -\infty}\hat{\phi}(\xi,\eta)=0.
\end{aligned}
\right.\qquad\text{for $(\xi,\eta)\in \mathbb{T}\times (-\infty,0)$}.
          \end{align}
Here $\Pi_{\perp}$ is the projection operator defined by
$$
\Pi_{\perp}f:=f-\bar{f}.
$$

(d) for $\alpha>1$, $B>0$, $(A+B)>0$ and
\begin{align}\label{M1}
(A+B)^2=\overline{f_{1}^2},\quad B^2=\overline{f_{0}^2}.
\end{align}
(\rmnum{2}) the higher-order interior profiles  $(u_{I}^{i},v_{I}^{i},p_{I}^{i})$ with $1\leq i\leq k_{0}$ satisfy the linearized Euler equations \eqref{1.9} with boundary condition \eqref{1.10} around the Couette flow $(u_{I}^{0},v_{I}^{0})=(Ay+B,0)$. Moreover, it holds that
\begin{align*}
\sum\limits_{j+l\leq m}\|\partial_{x}^{j}\partial_{y}^{l}(u_{I}^{i},v_{I}^{i},p_{I}^{i})\|_{L^2}\leq C(m)\quad\text{for any $1\leq i\leq k_{0}$}.
\end{align*}
(\rmnum{3}) For boundary-layer profiles, it holds that
$$
v_{B}^{i}=p_{B}^{i}\equiv 0\quad \text{for $0\leq i\leq \tilde{\alpha}-1$}.
$$
For $0\leq \alpha<1$, the zeroth-order lower boundary-layer profile $(u_{B}^{0},v_{B}^{\tilde{\alpha}},p_{B}^{\tilde{\alpha}})\equiv 0$. And in fact, 
$$
(u_{B}^{i},v_{B}^{i+\tilde{\alpha}},p_{B}^{i+\tilde{\alpha}})=(\hat{u}_{B}^{i},\hat{v}_{B}^{i+\tilde{\alpha}},\hat{p}_{B}^{i+\tilde{\alpha}})\equiv 0\quad \text{for $0\leq i\leq (1-\alpha)\tilde{\alpha}-1$},
$$
and the higher-order lower boundary-layer profile $(u_{B}^{i},v_{B}^{i+\tilde{\alpha}},p_{B}^{i+\tilde{\alpha}})$ with $i\geq (1-\alpha)\tilde{\alpha}$ is the solution of the linearized Prandtl equation \eqref{1.15} with boundary conditions \eqref{b00}--\eqref{b}.

For $\alpha=1$, the zeroth-order lower boundary-layer profile $(u_{B}^{0},v_{B}^{\tilde{\alpha}},p_{B}^{\tilde{\alpha}})$ is the solution of the Prandtl equation \eqref{1.13} with boundary conditions \eqref{b2}. The higher-order lower boundary-layer profile $(u_{B}^{i},v_{B}^{i+\tilde{\alpha}},p_{B}^{i+\tilde{\alpha}})$ with $1\leq i\leq k_{0}$ is the solution of the linearized Prandtl equations \eqref{1.15} with boundary conditions \eqref{b2}.

For $\alpha>1$, the zeroth-order lower boundary-layer profile $(u_{B}^{0},v_{B}^{\tilde{\alpha}},p_{B}^{\tilde{\alpha}})$ is the solution of the Prandtl equation \eqref{1.13} with boundary conditions \eqref{b3}. The higher-order lower boundary-layer profile $(u_{B}^{i},v_{B}^{i+\tilde{\alpha}},p_{B}^{i+\tilde{\alpha}})$ with $1\leq i\leq k_{0}$ is the solution of the linearized Prandtl equations \eqref{1.15} with boundary conditions \eqref{b3}.

(\rmnum{4}) For any $m\in \mathbb{N}_{+}$ and $\ell\gg 2$, the lower boundary-layer profile $(u_{B}^{i},v_{B}^{i+\tilde{\alpha}},p_{B}^{i+\tilde{\alpha}})$ with $i\geq 0$ satisfies
\begin{align}\label{4.1-1}
\sum\limits_{j+k\leq m}\int_{\mathbb{T}}\int_{0}^{\infty}\langle z\rangle^{\ell}|\partial_{x}^{j}\partial_{z}^{k}(u_{B}^{i},v_{B}^{i+\tilde{\alpha}},p_{B}^{i+\tilde{\alpha}})|^2\,{\rm d}z{\rm d}x\leq C(m,\ell)
\end{align}
provided $\alpha<1$, and 
\begin{align}\label{4.1}
\sum\limits_{j+k\leq m}\int_{\mathbb{T}}\int_{0}^{\infty}\langle z\rangle^{\ell}|\partial_{x}^{j}\partial_{z}^{k}(u_{B}^{i},v_{B}^{i+\tilde{\alpha}},p_{B}^{i+\tilde{\alpha}})|^2\,{\rm d}z{\rm d}x\leq C(m,\ell)(\sigma^2{\bf 1}_{i=0}+{\bf 1}_{i\geq 1}),
\end{align}
provided $\alpha\geq 1$.

(\rmnum{5}) The upper boundary profile  $(\hat{u}_{B}^{i},\hat{v}_{B}^{i+\tilde{\alpha}},\hat{p}_{B}^{i+\tilde{\alpha}})$ satisfies the similar boundary-layer equation as the corresponding lower boundary-layer profile with variable $z$ replaced by $\zeta$ and the traces of interior profiles at the lower boundary replaced by the ones at the upper boundary. Moreover, the estimates \eqref{4.1-1} and \eqref{4.1} also hold for $(\hat{u}_{B}^{i},\hat{v}_{B}^{i+\tilde{\alpha}},\hat{p}_{B}^{i+\tilde{\alpha}})$.
\end{proposition}

\noindent\textbf{Proof}. We divide the proof into four steps.

\textbf{Step 1}. The leading Euler flow  $(u_{I}^{0},v_{I}^{0},p_{I}^{0})=(Ay+B,0,0)$ is determined by Prandtl--Batchelor theory (cf. \cite{B-1956,FGLT-2023,FGLT-2024,GXX,PZ-2025}). The values $A$ and $B$ will be determined by the solvability of the corresponding boundary-layer problem, which will be presented in Step 2 below. The well-posedness of the higher-order interior profile $(u_{I}^{k},v_{I}^{k},p_{I}^{k})$ follows directly from Proposition \ref{interior}. Please note that the constants $\mathfrak{C}_{k}$ and $\tilde{\mathfrak{C}}_{k}$ will be determined by the corresponding boundary-layer problem, which will be explained in Step 3 below.

\textbf{Step 2}. For the leading-order boundary-layer profile $(u_{B}^{0},v_{B}^{\tilde{\alpha}},p_{B}^{\tilde{\alpha}})$, if $0\leq \alpha<1$, it follows from Lemma \ref{lem2} that $(u_{B}^{0},v_{B}^{\tilde{\alpha}},p_{B}^{\tilde{\alpha}})\equiv 0$. Moreover, noting \eqref{1.15} and $\eqref{b}_{1}$, one has $u_{B}^{i}$ ($0\leq i\leq (1-\alpha)\tilde{\alpha}-1$) satisfies \eqref{lp1} with homogeneous boundary condition $\eqref{b}_{1}$, so, using a similar argument as in the proof of Lemma \ref{lem2} and induction method, we obtain
\begin{align}\label{vanish}
(u_{B}^{i},v_{B}^{i+\tilde{\alpha}},p_{B}^{i+\tilde{\alpha}})\equiv (0,0,0)\quad \text{for }0\leq i\leq (1-\alpha)\tilde{\alpha}-1.
\end{align}
For $u_{B}^{(1-\alpha)\tilde{\alpha}}$, if $\alpha=0$, it satisfies
\begin{align*}
B\partial_{x}u_{B}^{(1-\alpha)\tilde{\alpha}}=\partial_{zz}u_{B}^{(1-\alpha)\tilde{\alpha}},
\end{align*}
with the boundary condition
\begin{align*}
\partial_{z}u_{B}^{(1-\alpha)\tilde{\alpha}}(x,0)=-f_{0}(x)+u_{I}^{0}(x,0)-\partial_{y}u_{I}^{0}(x,0)=-f_{0}(x)+B-A.
\end{align*}
Then the compatibility condition \eqref{slp1} reads
\begin{align}\label{B-A}
B-A=\bar{f}_{0}.
\end{align}
Moreover, the corresponding upper boundary profile $\hat{u}_{B}^{(1-\alpha)\tilde{\alpha}}$ satisfies
\begin{align*}
(A+B)\partial_{x}\hat{u}_{B}^{(1-\alpha)\tilde{\alpha}}=\partial_{\zeta\zeta}\hat{u}_{B}^{(1-\alpha)\tilde{\alpha}}
\end{align*}
with the boundary condition
\begin{align*}
\partial_{\zeta}\hat{u}_{B}^{(1-\alpha)\tilde{\alpha}}(x,0)=f_{1}(x)-u_{I}^{0}(x,1)-\partial_{y}u_{I}^{0}(x,1)=f_{1}(x)-(A+B)-A.
\end{align*}
Then the compatibility condition \eqref{slp1} reads
\begin{align}\label{B+A}
2A+B=\bar{f}_{1}.
\end{align}
Combining \eqref{B-A} and \eqref{B+A} gives \eqref{M00}

If $0<\alpha<1$, the first potentially nonzero boundary-layer profile $u_{B}^{(1-\alpha)\tilde{\alpha}}$ satisfies
\begin{align*}
B\partial_{x}u_{B}^{(1-\alpha)\tilde{\alpha}}=\partial_{zz}u_{B}^{(1-\alpha)\tilde{\alpha}},\quad \text{with }\partial_{z}u_{B}^{(1-\alpha)\tilde{\alpha}}(x,0)=-f_{0}(x)+u_{I}^{0}(x,0)=-f_{0}(x)+B.
\end{align*}
Then the compatibility condition \eqref{slp1} reads
$
B=\bar{f}_{0}.
$
Similarly, the choice of $(A+B)$ in \eqref{M0} comes from the compatibility condition of the equation of $\hat{u}_{B}^{(1-\alpha)\tilde{\alpha}}$.

If $\alpha>1$, it follows from \eqref{1.13} and the Batchelor--Wood formula \eqref{bw} that $A$ and $B$ should be determined by \eqref{M1}.
If $\alpha=1$, it follows from \eqref{1.13} and Lemma \ref{lem1} that $A$ and $B$ should be chosen to satisfy \eqref{M2}.

\textbf{Step 3}. For $0\leq \alpha<1$, to determine $\mathfrak{C}_{k}$ and $\tilde{\mathfrak{C}}_{k}$ in the construction of $u_{I}^{k}$, we consider the equation of $u_{B}^{(1-\alpha)\tilde{\alpha}+k}$. It follows from \eqref{1.15} and $\eqref{b00}$--$\eqref{b}$ that $u_{B}^{(1-\alpha)\tilde{\alpha}+k}$ satisfies
\begin{align*}
\left\{
\begin{aligned}
    &B\partial_{x}u_{B}^{(1-\alpha)\tilde{\alpha}+k}=\partial_{zz}u_{B}^{(1-\alpha)\tilde{\alpha}+k}+F_{B}^{k+(1-\alpha)\tilde{\alpha}},\\
&\partial_{z}u_{B}^{(1-\alpha)\tilde{\alpha}+k}(x,0)=u_{I}^{k}(x,0)+u_{B}^{k}(x,0)-\partial_{y}u_{I}^{k-\alpha\tilde{\alpha}}(x,0).
\end{aligned}
\right.
\end{align*}
Then the compatibility condition \eqref{slp1} gives
\begin{align}\label{Ck}
\int_{\mathbb{T}}u_{I}^{k}(x,0)\,{\rm d}x=\int_{\mathbb{T}}\int_{0}^{\infty}F_{B}^{k+(1-\alpha)\tilde{\alpha}}\,{\rm d}z{\rm d}x-\int_{\mathbb{T}}u_{B}^{k}(x,0)\,{\rm d}x+\int_{\mathbb{T}}\partial_{y}u_{I}^{k-\alpha\tilde{\alpha}}(x,0)\,{\rm d}x.
\end{align}
Noting \eqref{vanish}, it is easy to check that there is no interior profile of order greater than $k$ entering the integral $\int_{\mathbb{T}}\int_{0}^{\infty}F_{B}^{k+(1-\alpha)\tilde{\alpha}}\,{\rm d}z{\rm d}x$. In addition, for terms involving $u_{I}^{k}(x,0)$ in the integral $\int_{\mathbb{T}}\int_{0}^{\infty}F_{B}^{k+(1-\alpha)\tilde{\alpha}}\,{\rm d}z{\rm d}x$, using \eqref{vanish}, a direct calculation shows that
\begin{align*}
-\int_{\mathbb{T}}\int_{0}^{\infty}\big[u_{I}^{k}(x,0)\partial_{x}u_{B}^{(1-\alpha)\tilde{\alpha}}+u_{B}^{(1-\alpha)\tilde{\alpha}}\partial_{x}u_{I}^{k}(x,0)+v_{B}^{(1-\alpha)\tilde{\alpha}}\partial_{y}u_{I}^{k}(x,0)\big]\,{\rm d}z{\rm d}x=0.
\end{align*}
So, the terms on the right hand side of \eqref{Ck} only consist of $v_{I}^{k}(x,0)$ and all lower-order interior and boundary-layer profiles prior to $u_{I}^{k}$ when $0<\alpha<1$. And when $\alpha=0$, we can move the term $\int_{\mathbb{T}}\partial_{y}u_{I}^{k}(x,0)\,{\rm d}x$ in \eqref{Ck} to the left hand side. Similarly, from the compatibility condition of the upper boundary-layer profile $\hat{u}_{B}^{(1-\alpha)\tilde{\alpha}+k}$, one has 
\begin{align}\label{hCk}
\int_{\mathbb{T}}u_{I}^{k}(x,1)\,{\rm d}x=\int_{\mathbb{T}}\int_{-\infty}^{0}\hat{F}_{B}^{(1-\alpha)\tilde{\alpha}+k}\,{\rm d}\zeta{\rm d}x-\int_{\mathbb{T}}\hat{u}_{B}^{k}(x,0)\,{\rm d}x-\int_{\mathbb{T}}\partial_{y}u_{I}^{k-\alpha\tilde{\alpha}}(x,1)\,{\rm d}x.
\end{align}
Then the coefficients $\tilde{\mathfrak{C}}_{k}$ and $\mathfrak{C}_{k}$ in constructions of $u_{I}^{k}$ will be uniquely determined through \eqref{Ck} and \eqref{hCk}.

For $\alpha\geq 1$, 
we first choose $\mathfrak{C}_{k}=\tilde{\mathfrak{C}}_{k}=0$ in the construction of $u_{I}^{k}$. Then we substitute it into the equation of $u_{B}^{k}$ and $\hat{u}_{B}^{k}$. Then we obtain from Lemma \ref{lem3} that there exist two constants $\mathfrak{A}_{k,\infty}$ and $\hat{\mathfrak{A}}_{k,\infty}$ such that
$$
\lim\limits_{z\to \infty}u_{B}^{k}(x,z)=\mathfrak{A}_{k,\infty},\quad \lim\limits_{\zeta\to -\infty}\hat{u}_{B}^{k}(x,\zeta)=\hat{\mathfrak{A}}_{k,\infty}.
$$
Then we choose $\mathfrak{C}_{k}=\hat{\mathfrak{A}}_{k,\infty}-\mathfrak{A}_{k,\infty}$, $\tilde{\mathfrak{C}}_{k}=\mathfrak{A}_{k,\infty}$, and modify $u_{I}^{k}$, $u_{B}^{k}$ and $\hat{u}_{B}^{k}$ as
$$
\tilde{u}_{I}^{k}(x,y)=u_{I}^{k}(x,y)+\mathfrak{C}_{k}y+\tilde{\mathfrak{C}}_{k}=u_{I}^{k}(x,y)+(1-y)\mathfrak{A}_{k,\infty}+y\hat{\mathfrak{A}}_{k,\infty}
$$
and
$$
\tilde{u}_{B}^{k}(x,z)=u_{B}^{k}(x,z)-\mathfrak{A}_{k,\infty},\quad \tilde{\hat{u}}_{B}^{k}(x,\zeta)=\hat{u}_{B}^{k}(x,\zeta)-\hat{\mathfrak{A}}_{k,\infty}.
$$
Then, $\lim\limits_{z\to \infty}\tilde{u}_{B}^{k}=\lim\limits_{\zeta\to -\infty}\tilde{\hat{u}}_{B}^{k}=0$. Moreover, the $k$th order Prandtl profiles
\begin{align*}
&\tilde{u}_{p}^{k}(x,z):=\tilde{u}_{B}^{k}(x,z)+\tilde{u}_{I}^{k}(x,0)=u_{B}^{k}(x,z)+u_{I}^{k}(x,0),\\
&\tilde{\hat{u}}_{p}^{k}(x,\zeta)=\tilde{\hat{u}}_{B}^{k}(x,\zeta)+\tilde{u}_{I}^{k}(x,1)=\hat{u}_{B}^{k}(x,\zeta)+u_{I}^{k}(x,1),
\end{align*}
remain unchanged and preserve the corresponding boundary conditions.
Then, for subsequent interior and boundary-layer equations, we substitute modified interior and boundary-layer profiles, and repeat the above procedure to pick up the subsequent coefficients $\tilde{\mathfrak{C}}_{k}$ and $\mathfrak{C}_{k}$. It is worth pointing out that this modification does not affect well-posedness of the equation of subsequent interior and boundary-layer profiles as it appears in the source terms. For simplicity of notations, we still denote the modified interior and boundary-layers profiles by $u_{I}^{k}$, $u_{B}^{k}$ and $\hat{u}_{B}^{k}$ respectively in our expansions.

\textbf{Step 4}. The estimates \eqref{4.1-1} and \eqref{4.1} follow from the induction method and Corollaries \ref{cor1}, \ref{cor2} and Lemmas \ref{lem2}, \ref{lem3} and \ref{lem4}. Therefore, the proof of Proposition \ref{prop3} is complete. $\hfill\square$

By Proposition \ref{prop3}, we summarize the properties of the approximate solutions $(u^{a},v^{a})$.

\begin{lemma}\label{lem5.1}
	The approximate solution $(u^{a},v^{a})$ satisfies \eqref{3.4} and for all integers $\ell\gg 2$ and $k$, we have the following
	\begin{align*}
	\left\{\begin{aligned}
    &v_{I}^{j}\equiv 0,\quad u_{I}^{j}=\tilde{C}_{j}y+C_{j},\quad \text{for $0\leq j\leq \tilde{\alpha}-1$},\\
	&\|\nabla_{x,y}^{k}(u_{I}^{a}-(Ay+B),v_{I}^{a})\|_{L^{\infty}}+\|\langle z\rangle^{\ell}\nabla_{x,z}^{k}(u_{B}^{a},v_{B}^{a})\|_{L^{\infty}}+\|\langle \zeta\rangle^{\ell}\nabla_{x,\zeta}^{k}(\hat{u}_{B}^{a},\hat{v}_{B}^{a})\|_{L^{\infty}}\leq C(\varepsilon^{\frac{1}{\tilde{\alpha}}}+\sigma),
	\end{aligned}
	\right.
	\end{align*}
    for $\alpha\geq 1$, and
    \begin{align*}
    	\left\{\begin{aligned}
        &v_{I}^{j}\equiv 0,\quad u_{I}^{j}=\tilde{C}_{j}y+C_{j},\quad \text{for $0\leq j\leq (2-\alpha)\tilde{\alpha}-1$},\\
    	&\|\nabla_{x,y}^{k}(u_{I}^{a}-(Ay+B),v_{I}^{a})\|_{L^{\infty}}+\|\langle z\rangle^{\ell}\nabla_{x,z}^{k}(u_{B}^{a},v_{B}^{a})\|_{L^{\infty}}+\|\langle \zeta\rangle^{\ell}\nabla_{x,\zeta}^{k}(\hat{u}_{B}^{a},\hat{v}_{B}^{a})\|_{L^{\infty}}\leq C\varepsilon^{1-\alpha},
    	\end{aligned}
    	\right.
    \end{align*}
for $0\leq \alpha<1$.
\end{lemma}

\noindent\textbf{Proof}. The estimates follow directly from \eqref{4.1-1}--\eqref{4.1} and the Sobolev embedding. Moreover, for $\alpha\geq 1$, noting that
$$
v_{B}^{i}=\hat{v}_{B}^{i}\equiv 0\quad \text{for $1\leq i\leq \tilde{\alpha}-1$},
$$
we obtain from \eqref{v} that $v_{I}^{i}\equiv 0$ for $1\leq i\leq \tilde{\alpha}-1$. Then it follows from \eqref{u} that
$$
u_{I}^{i}=\tilde{C}_{i}y+C_{i}\quad \text{for $1\leq i\leq \tilde{\alpha}-1$}.
$$
For $\alpha<1$, since $u_{B}^{i}=\hat{u}_{B}^{i}\equiv 0$ for $0\leq i\leq (1-\alpha)\tilde{\alpha}-1$, one obtains from the incompressible condition that
$$
v_{B}^{i}=\hat{v}_{B}^{i}\equiv 0\quad \text{for $0\leq i\leq (2-\alpha)\tilde{\alpha}-1$}.
$$
Then, we obtain from \eqref{v}--\eqref{u} that
$$
v_{I}^{i}\equiv 0,\quad u_{I}^{i}=\tilde{C}_{i}y+C_{i},\quad \text{for $1\leq i\leq (2-\alpha)\tilde{\alpha}-1$}.
$$
This completes the proof of Lemma \ref{lem5.1} $\hfill\square$

For the auxiliary functions $(h_{u},h_{v})$ in \eqref{3.5-0} with \eqref{3.5-2}--\eqref{3.5-3}, noting from the incompressible condition that 
$$
\int_{\mathbb{T}}v_{B}(x,z)\,{\rm d}x=\int_{\mathbb{T}}\hat{v}_{B}(x,\zeta)\,{\rm d}x=0,
$$
it is direct to check that the condition \eqref{B1-1} holds for \eqref{3.5-0}. The constructions of $(h_{u},h_{v})$ are presented in Appendix \ref{AppendixB}. Moreover, it follows from such constructions and Proposition \ref{prop3}, Lemma \ref{lem5.1} that
$$
\sum\limits_{j+l\leq m}\|\partial_{x}^{j}\partial_{y}^{l}(h_{u},h_{v})\|_{L^2}\leq C(m)\quad \text{for any $m\geq 0$}.
$$

\section{Remainder Estimates and Proof of Main Theorems}

In this section, we establish the uniform estimates and existence of the remainder equations. Let
$$
u=u^{\varepsilon}-u^{a},\quad v=v^{\varepsilon}-v^{a},\quad p=p^{\varepsilon}-p^{a}.
$$
Then, the remainder $(u,v,p)$ satisfies
\begin{align}\label{3.8}
	\left\{
	\begin{aligned}
		&u^{a}\partial_{x}u+v^{a}\partial_{y}u+u\partial_{x}u^{a}+v\partial_{y}u^{a}+\partial_{x}p=\varepsilon^{2}\Delta u+R_{u},\\
		&u^{a}\partial_{x}v+v^{a}\partial_{y}v+u\partial_{x}v^{a}+v\partial_{y}v^{a}+\partial_{y}p=\varepsilon^2\Delta v+R_{v},\\
		&\partial_{x}u+\partial_{y}v=0,\\
		&v(x,0)=v(x,1)=0,\\
		&(u-\varepsilon^{\alpha}\partial_{y}u)(x,0)=(u+\varepsilon^{\alpha}\partial_{y}u)(x,1)=0,
	\end{aligned}
	\right.
\end{align}
where $(R_{u},R_{v})$ are given as follows:
\begin{align*}
	R_{u}&=-(u^{a}\partial_{x}+v^{a}\partial_{y})u^{a}-\partial_{x}p^{a}+\varepsilon^2\Delta u^{a}-(u\partial_{x}+v\partial_{y})u=:R_{u}^{a}-(u\partial_{x}+v\partial_{y})u,\\
	R_{v}&=-(u^{a}\partial_{x}+v^{a}\partial_{y})v^{a}-\partial_{y}p^{a}+\varepsilon^{2}\Delta v^{a}-(u\partial_{x}+v\partial_{y})v=:R_{v}^{a}-(u\partial_{x}+v\partial_{y})v.
\end{align*}
By a direct and careful computation, we have
\begin{align}\label{3.9}
    \|(R_{u}^{a},R_{v}^{a})\|_{L^{\infty}}\leq C\varepsilon^{\frac{k_{0}}{\tilde{\alpha}}},
\end{align}
provided that $\mathfrak{b}\geq \frac{k_{0}}{\tilde{\alpha}}$ in \eqref{1.11-1}--\eqref{1.11-2}. Moreover, it follows from $\eqref{3.8}_{3}$--$\eqref{3.8}_{4}$ that
\begin{align}\label{3.8-1}
	\int_{\mathbb{T}}v(x,y)\,{\rm d}x=0.
\end{align}

For later use, we recall the classical Hardy inequality.
\begin{lemma}\label{Hardy}
Let $f\in H^{1}(\Omega)$ with $f(x,0)=0$,  then 
\begin{align*}
\int_{\mathbb{T}}\int_{0}^{1}\frac{f^2}{y^2}\,{\rm d}y{\rm d}x\leq 4\int_{\mathbb{T}}\int_{0}^{1}|f_{y}|^2\,{\rm d}y{\rm d}x.
\end{align*}
\end{lemma}

\noindent\textbf{Proof}. Using $f(x,0)=0$, one has
$$
f(x,y)=\int_0^yf_y(x,s)\,{\rm d}s,
$$
which, together with the H\"{o}lder inequality, yields
\begin{align*}
f^2(x,y)\leq y\int_0^yf_y^2(x,s)\,{\rm d}s.
\end{align*}
In particular, for any $\tau>0$, 
\begin{align}\label{tau}
\int_{\mathbb{T}}\frac{f^2(x,\tau)}{\tau}\,{\rm d}x\leq \int_{\mathbb{T}}\int_0^\tau f_y^2(x,s)\,{\rm d}s{\rm d}x\to 0,\quad \text{as $\tau\to 0$}.
\end{align}
Then,
\begin{align*}
\int_{\mathbb{T}}\int_\tau^1\frac{f^2}{y^2}\,{\rm d}y{\rm d}x&=-\int_{\mathbb{T}}f^2(x,1)\,{\rm d}x+\int_{\mathbb{T}}\frac{f^2(x,\tau)}{\tau}\,{\rm d}x+2\int_{\mathbb{T}}\int_\tau^1\frac{ff_y}{y}\,{\rm d}y{\rm d}x\nonumber\\
&\leq \int_{\mathbb{T}}\frac{f^2(x,\tau)}{\tau}\,{\rm d}x+2\Big(\int_{\mathbb{T}}\int_\tau^1\frac{f^2}{y^2}\,{\rm d}y{\rm d}x\Big)^{\frac{1}{2}}\Big(\int_{\mathbb{T}}\int_\tau^1|f_y|^2\,{\rm d}y{\rm d}x\Big)^{\frac{1}{2}}.
\end{align*}
Using the fact $2ab\leq \frac{1}{2}a^2+2b^2$, we obtain
\begin{align*}
\int_{\mathbb{T}}\int_\tau^1\frac{f^2}{y^2}\,{\rm d}y{\rm d}x\leq 2\int_{\mathbb{T}}\frac{f^2(x,\tau)}{\tau}\,{\rm d}x+4\int_{\mathbb{T}}\int_\tau^1|f_y|^2\,{\rm d}y{\rm d}x.
\end{align*}
Taking $\tau\to 0$ in the above inequality and using \eqref{tau} completes the proof of Lemma \ref{Hardy}. $\hfill\square$

By similar arguments, we also have following lemma.
\begin{lemma}\label{Hardy3}
Let $f\in H^{1}(\Omega)$ with $f(x,1)=0$, then
\begin{align*}
\int_{\mathbb{T}}\int_{0}^{1}\frac{f^2}{(1-y)^2}\,{\rm d}y{\rm d}x\leq 4\int_{\mathbb{T}}\int_{0}^{1}|f_{y}|^2\,{\rm d}y{\rm d}x.
\end{align*}
\end{lemma}

\begin{lemma}\label{So1}
For any $u\in H^{1}(\Omega)$, it holds that
\begin{align*}
\|u\|_{L^2}^2\leq C\Big(\|\partial_{y}u\|_{L^2}^2+\|\partial_{x}u\|_{L^2}^2+\int_{\mathbb{T}}|u(x,0)|^2\,{\rm d}x\Big).
\end{align*}
\end{lemma}

\noindent\textbf{Proof}. We decompose $u$ as
\begin{align}\label{dec1}
u(x,y)=\bar{u}(y)+\tilde{u}(x,y).
\end{align}
It follows directly from the Poincar\'{e} inequality that
\begin{align}\label{tildeu1}
\|\tilde{u}\|_{L^2}\leq C\|\partial_{x}u\|_{L^2}.
\end{align}
For $\bar{u}(y)$, a direct calculation shows that
\begin{align}\label{tildeu}
\bar{u}(y)=\bar{u}(0)+\int_{0}^{y}\bar{u}'(s)\,{\rm d}s.
\end{align}
Integrating the above equation over $\mathbb{T}\times [0,1]$ gives
\begin{align}\label{baru}
\|\bar{u}(y)\|_{L^2}^2\leq C\Big(\|\partial_{y}u\|_{L^2}^2+\int_{\mathbb{T}}u^2(x,0)\,{\rm d}x\Big).
\end{align}
Combining \eqref{tildeu1}--\eqref{baru} completes the proof of Lemma \ref{So1}. $\hfill\square$

\subsection{Basic energy estimates} Multiplying $\eqref{3.8}_{1}$ by $u$, $\eqref{3.8}_{2}$ by $v$, and adding the resulting equations together, we have
\begin{align}\label{3.11}
	&\varepsilon^{2}\int_{\mathbb{T}}\int_{0}^{1}|\nabla u|^2+|\nabla v|^2\,{\rm d}y{\rm d}x+\varepsilon^{2-\alpha}\int_{\mathbb{T}}|u(x,0)|^2+|u(x,1)|^2\,{\rm d}x\nonumber\\
	&\leq -\int_{\mathbb{T}}\int_{0}^{1}\big(u^{a}u\partial_{x}u+v^{a}u\partial_{y}u+u^2\partial_{x}u^{a}+uv\partial_{y}u^{a}\big)\,{\rm d}y{\rm d}x\nonumber\\
	&\quad -\int_{\mathbb{T}}\int_{0}^{1}\big(u^{a}v\partial_{x}v+v^{a}v\partial_{y}v+uv\partial_{x}v^{a}+v^{2}\partial_{y}v^{a}\big)\,{\rm d}y{\rm d}x\nonumber\\
	&\quad -\int_{\mathbb{T}}\int_{0}^{1}\big(u\partial_{x}p+v\partial_{y}p\big)\,{\rm d}y{\rm d}x+\int_{\mathbb{T}}\int_{0}^{1}\big(uR_{u}+vR_{v}\big)\,{\rm d}y{\rm d}x\nonumber\\
	&=:\sum\limits_{i=1}^{4}S_{i}.
\end{align}
It follows directly from the divergence-free condition and integration by parts that
\begin{align}\label{3.12}
	S_{3}=\int_{\mathbb{T}}\int_{0}^{1}p(\partial_{x}u+\partial_{y}v)\,{\rm d}y{\rm d}x=0.
\end{align}
For $S_{4}$, using the Poincar\'{e} inequality and Lemma \ref{So1}, a direct calculation shows that
\begin{align}\label{3.12-1}
|S_{4}|&\leq C\varepsilon^{2}\tau(\|u\|_{L^2}^2+\|v\|_{L^2}^2)+C_{\tau}\varepsilon^{-2}\|(R_{u},R_{v})\|_{L^2}^2\nonumber\\
&\leq C\varepsilon^2\tau\Big(\|\partial_{y}u\|_{L^2}^2+\|\partial_{x}u\|_{L^2}^2+\|\partial_{x}v\|_{L^2}^2+\int_{\mathbb{T}}u^{2}(x,0)\,{\rm d}x\Big)+C_{\tau}\varepsilon^{-2}\|(R_{u},R_{v})\|_{L^2}^2.
\end{align}
For the convection terms $S_{1}$ and $S_{2}$, one has
\begin{align}\label{3.13}
	S_{1}+S_{2}&=-\int_{\mathbb{T}}\int_{0}^{1}\big(u^{a}u\partial_{x}u+v^{a}u\partial_{y}u\big)\,{\rm d}y{\rm d}x-\int_{\mathbb{T}}\int_{0}^{1}\big(u^{a}v\partial_{x}v+v^{a}v\partial_{y}v\big)\,{\rm d}y{\rm d}x\nonumber\\
	&\quad -\int_{\mathbb{T}}\int_{0}^{1}\big(uv\partial_{x}v^{a}+v^2\partial_{y}v^{a}\big)\,{\rm d}y{\rm d}x
	-\int_{\mathbb{T}}\int_{0}^{1}\big(u^{2}\partial_{x}u^{a}+uv\partial_{y}u^{a}\big)\,{\rm d}y{\rm d}x\nonumber\\
	&=:\sum\limits_{i=1}^{4}H_{i}.
\end{align}
For $H_{1}$ and $H_{2}$, it follows from integration by parts and the fact $\partial_{x}u^{a}+\partial_{y}v^{a}=0$ that
\begin{align}\label{3.14}
	H_{1}+H_{2}=\int_{\mathbb{T}}\int_{0}^{1}\frac{1}{2}u^{2}(\partial_{x}u^{a}+\partial_{y}v^{a})\,{\rm d}y{\rm d}x+\int_{\mathbb{T}}\int_{0}^{1}\frac{1}{2}v^2\big(\partial_{x}u^{a}+\partial_{y}v^{a}\big)\,{\rm d}y{\rm d}x=0.
\end{align}
For $H_{3}$, using \eqref{3.8-1}, the Poincar\'{e} inequality, Proposition \ref{prop3} and Lemmas \ref{lem5.1}, \ref{So1}, we obtain
\begin{align}\label{3.16}
|H_{3}|&\leq C\varepsilon^{1-\alpha}\|v\|_{L^2}^2+C\varepsilon^{2-\alpha}\|u\|_{L^2}\|v\|_{L^2}\leq C\|\partial_{x}v\|_{L^2}^2+C\varepsilon^{4-2\alpha}\|u\|_{L^2}^2\nonumber\\
&\leq C\|\partial_{x}v\|_{L^2}^2+C\varepsilon^{4-2\alpha}\Big(\|\partial_{x}u\|_{L^2}^2+\|\partial_{y}u\|_{L^2}^2+\int_{\mathbb{T}}|u(x,0)|^2\,{\rm d}x\Big)
\end{align}
provided that $0\leq \alpha<1$ and $k_{0}\geq (2-\alpha)\tilde{\alpha}$. And for $\alpha\geq 1$, 
\begin{align}\label{3.16-1}
|H_{3}|&\leq C\|v\|_{L^2}^2+C\varepsilon(\varepsilon^{\frac{1}{\tilde{\alpha}}}+\sigma)\|u\|_{L^2}\|v\|_{L^2}\leq C\|\partial_{x}v\|_{L^2}^2+C\varepsilon^2(\varepsilon^{\frac{1}{\tilde{\alpha}}}+\sigma)\|u\|_{L^2}^2\nonumber\\
&\leq C\|\partial_{x}v\|_{L^2}^2+C\varepsilon^2(\varepsilon^{\frac{1}{\tilde{\alpha}}}+\sigma)\Big(\|\partial_{x}u\|_{L^2}^2+\|\partial_{y}u\|_{L^2}^2+\int_{\mathbb{T}}|u(x,0)|^2\,{\rm d}x\Big).
\end{align}
For $H_{4}$, the terms involving the interior approximate solution $u_{I}^{a}$ are
\begin{align}\label{3.26}
	&\int_{\mathbb{T}}\int_{0}^{1}\big(u^2\partial_{x}u_{I}^{a}+uv\partial_{y}u_{I}^{a}\big)\,{\rm d}y{\rm d}x\nonumber\\
	&=\int_{\mathbb{T}}\int_{0}^{1}\big(\bar{u}^2+2\bar{u}\tilde{u}+\tilde{u}^2\big)\partial_{x}u_{I}^{a}\,{\rm d}y{\rm d}x+\int_{\mathbb{T}}\int_{0}^{1}(\bar{u}v+\tilde{u}v)\partial_{y}u_{I}^{a}\,{\rm d}y{\rm d}x=:H_{41}+H_{42},
\end{align}
where we have used the decomposition \eqref{dec1}. For $H_{41}$, it is straightforward to see
\begin{align}\label{3.27}
	\int_{\mathbb{T}}\int_{0}^1\bar{u}^2\partial_{x}u_{I}^{a}\,{\rm d}y{\rm d}x=0.
\end{align}
Using Lemma \ref{lem5.1}, the Poincar\'{e} and H\"{o}lder inequalities, we have that for $0\leq \alpha<1$,
\begin{align}\label{3.28}
	\int_{\mathbb{T}}\int_{0}^{1}\tilde{u}^2\partial_{x}u_{I}^{a}\,{\rm d}y{\rm d}x&\leq C\|\partial_{x}u_{I}^{a}\|_{L^{\infty}}\|\partial_{x}u\|_{L^2}^2\leq C\|\partial_{x}u\|_{L^2}^2,
\end{align}
and
\begin{align}\label{3.29}
	&\int_{\mathbb{T}}\int_{0}^{1}\bar{u}\tilde{u}\partial_{x}u_{I}^{a}\,{\rm d}y{\rm d}x
	\nonumber\\&\leq 
    \left\{
\begin{aligned}
& C\|\partial_{x}u\|_{L^2}^2+C\varepsilon^{4-2\alpha}\Big(\|\partial_{y}u\|_{L^2}^2+\int_{\mathbb{T}}|u(x,0)|^2\,{\rm d}x\Big),\,\qquad\qquad\text{when $0\leq \alpha<1$},\\
 &C\|\partial_{x}u\|_{L^2}^2+C\varepsilon^2(\varepsilon^{\frac{1}{\tilde{\alpha}}}+\sigma)\Big(\|\partial_{y}u\|_{L^2}^2+\int_{\mathbb{T}}|u(x,0)|^2\,{\rm d}x\Big),\qquad \text{when $\alpha\geq 1$}.
\end{aligned}
    \right.
\end{align}
Therefore, $H_{41}$ can be bounded as
\begin{align}\label{3.30}
	|H_{41}|\leq  
     \left\{
\begin{aligned}
& C\|\partial_{x}u\|_{L^2}^2+C\varepsilon^{4-2\alpha}\Big(\|\partial_{y}u\|_{L^2}^2+\int_{\mathbb{T}}|u(x,0)|^2\,{\rm d}x\Big),\,\qquad\qquad\text{when $0\leq \alpha<1$},\\
 &C\|\partial_{x}u\|_{L^2}^2+C\varepsilon^2(\varepsilon^{\frac{1}{\tilde{\alpha}}}+\sigma)\Big(\|\partial_{y}u\|_{L^2}^2+\int_{\mathbb{T}}|u(x,0)|^2\,{\rm d}x\Big),\qquad \text{when $\alpha\geq 1$}.
\end{aligned}
    \right.
\end{align}
For $H_{42}$, noting from Lemma \ref{lem5.1} that
$$
\partial_{y}u_{I}^{j}=\left\{
\begin{aligned}
&\tilde{C}_{j},\quad \text{for $0\leq j\leq (2-\alpha)\tilde{\alpha}-1$ when $0\leq \alpha<1$},\\
&\tilde{C}_{j},\quad \text{for $0\leq j\leq \tilde{\alpha}-1$ when $\alpha\geq 1$},
\end{aligned}
\right.
$$
one has from \eqref{3.8-1} and Proposition \ref{prop3} that
\begin{align*}
|H_{42}|
	\leq 
    \left\{
    \begin{aligned}
    &C(\|\partial_{x}v\|_{L^2}^2+\|\partial_{x}u\|_{L^2}^2)+C\varepsilon^{4-2\alpha}\Big(\|\partial_{y}u\|_{L^2}^2+\int_{\mathbb{T}}|u(x,0)|^2\,{\rm d}x\Big),\qquad\,\,\,\text{when $0\leq \alpha<1$},\\
&C(\|\partial_{x}v\|_{L^2}^2+\|\partial_{x}u\|_{L^2}^2)+C\varepsilon^2(\varepsilon^{\frac{1}{\tilde{\alpha}}}+\sigma)\Big(\|\partial_{y}u\|_{L^2}^2+\int_{\mathbb{T}}|u(x,0)|^2\,{\rm d}x\Big),\,\,\text{when $\alpha\geq 1$}.
    \end{aligned}
    \right.
\end{align*}
The terms involving $u_{B}^{a}+\hat{u}_{B}^{a}$ in $H_{4}$ are
\begin{align}\label{3.34-0}
&	\iintxy\big[u^2\partial_{x}(u_{B}^{a}+\hat{u}_{B}^{a})+uv\partial_{y}(u_{B}^{a}+\hat{u}_{B}^{a})\big]\,{\rm d}y{\rm d}x
=:H_{43}+H_{44}.
\end{align}
For $H_{43}$, when $0\leq \alpha<1$, using decompositions \ref{dec1}, one has
\begin{align*}
    u(x,y)=\bar{u}(y)+\tilde{u}(x,y).
\end{align*}
Then, 
\begin{align*}
H_{43}&=\iintxy(\bar{u}^2+2\bar{u}\tilde{u}+\tilde{u}^2)\partial_{x}(u_{B}^{a}+\hat{u}_{B}^{a})\,{\rm d}y{\rm d}x\nonumber\\
&=2\iintxy\bar{u}\tilde{u}\partial_{x}(u_{B}^{a}+\hat{u}_{B}^{a})\,{\rm d}y{\rm d}x+\iintxy\tilde{u}^2\partial_{x}(u_{B}^{a}+\hat{u}_{B}^{a})\,{\rm d}y{\rm d}x.
\end{align*}
Using the Poincar\'{e} inequality, one has
\begin{align*}
\iintxy\tilde{u}^2\partial_{x}(u_{B}^{a}+\hat{u}_{B}^{a})\,{\rm d}y{\rm d}x
\leq C\|\partial_{x}u\|_{L^2}^2.
\end{align*}
For the remaining term, we decompose $\bar{u}(y)$ as
$$
\bar{u}(y)=\bar{u}(0)+(\bar{u}(y)-\bar{u}(0))=:\bar{u}(0)-\bar{u}_{0}(y),
$$
near the lower boundary, and
$$
\bar{u}(y)=\bar{u}(1)+(\bar{u}(y)-\bar{u}(1))=:\bar{u}(1)+\bar{u}_{1}(y),
$$
near the upper boundary. Then $\bar{u}_{0}(0)=\bar{u}_{1}(1)=0$, and we obtain from Lemmas \ref{Hardy}--\ref{Hardy3} that
\begin{align*}
&2\iintxy\bar{u}\tilde{u}\partial_{x}(u_{B}^{a}+\hat{u}_{B}^{a})\,{\rm d}y{\rm d}x\nonumber\\
&\leq C\|\partial_{x}u\|_{L^2}^2+C\varepsilon^{4-2\alpha}\iintxy\Big(\frac{|\bar{u}_{0}(y)|^2}{y^2}\chi(y)+\frac{|\bar{u}_{1}(y)|^2}{(y-1)^2}\tilde{\chi}(y)\Big)\,{\rm d}y{\rm d}x\nonumber\\
&\quad +C\varepsilon^{2-2\alpha+2\lambda}\iintxy\Big(\frac{|\bar{u}(0)|^2}{y^{2\lambda}}+\frac{|\bar{u}(1)|^2}{|y-1|^{2\lambda}}\Big)\,{\rm d}y{\rm d}x\nonumber\\
&\leq C\|\partial_{x}u\|_{L^2}^2+C\varepsilon^{4-2\alpha}\|\partial_{y}u\|_{L^2}^2\nonumber\\
&\quad +C\varepsilon^{2-2\alpha+2\lambda}\Big(\int_{\mathbb{T}}\big(|u(x,0)|^2+|u(x,1)|^2\big)\,{\rm d}x\Big),
\end{align*}
provided that $2\lambda<1$. Therefore, when $0\leq \alpha<1$, $H_{43}$ can be bounded as
\begin{align*}
|H_{43}|\leq C\|\partial_{x}u\|_{L^2}^2+C\varepsilon^{4-2\alpha}\|\partial_{y}u\|_{L^2}^2+C\varepsilon^{2-2\alpha+2\lambda}\Big(\int_{\mathbb{T}}\big(|u(x,0)|^2+|u(x,1)|^2\big)\,{\rm d}x\Big).
\end{align*}
For $H_{44}$, by using Proposition \ref{prop3}, Lemmas \ref{Hardy}--\ref{Hardy3}, we have
\begin{align*}
&\iintxy\tilde{u}v\partial_{y}(u_{B}^{a}+\hat{u}_{B}^{a})\,{\rm d}y{\rm d}x\nonumber\\
&\leq C\|\partial_{x}u\|_{L^2}^2+C\varepsilon^{2-2\alpha}\iintxy \big(\frac{|v|^2}{y^2}\chi(y)+\frac{|v|^2}{(y-1)^2}\tilde{\chi}(y)\big)\,{\rm d}y{\rm d}x\nonumber\\
&\leq C\|\partial_{x}u\|_{L^2}^2+C\varepsilon^{2-2\alpha}\|\partial_{y}v\|_{L^2}^2\leq C\|\partial_{x}u\|_{L^2}^2,
\end{align*}
and
\begin{align*}
\iintxy\bar{u}(y)v\partial_{y}(u_{B}^{a}+\hat{u}_{B}^{a})\,{\rm d}y{\rm d}x&\leq C\iintxy \big(\frac{|v|^2}{y^2}\chi(y)+\frac{|v|^2}{(y-1)^2}\tilde{\chi}(y)\big)\,{\rm d}y{\rm d}x\nonumber\\
&\quad +C\varepsilon^{4-2\alpha}\iintxy\Big(\frac{|\bar{u}_{0}(y)|^2}{y^2}\chi(y)+\frac{|\bar{u}_{1}(y)|^2}{(y-1)^2}\tilde{\chi}(y)\Big)\,{\rm d}y{\rm d}x\nonumber\\
&\quad +C\varepsilon^{2-2\alpha+2\lambda}\iintxy\Big(\frac{|\bar{u}(0)|^2}{y^{2\lambda}}+\frac{|\bar{u}(1)|^2}{|y-1|^{2\lambda}}\Big)\,{\rm d}y{\rm d}x\nonumber\\
&\leq C\|\partial_{x}u\|_{L^2}^2+C\varepsilon^{4-2\alpha}\|\partial_{y}u\|_{L^2}^2\nonumber\\
&\quad +C\varepsilon^{2-2\alpha+2\lambda}\Big(\int_{\mathbb{T}}\big(|u(x,0)|^2+|u(x,1)|^2\big)\,{\rm d}x\Big).
\end{align*}
Therefore, when $0\leq \alpha<1$, $H_{44}$ can be bounded as
\begin{align*}
    |H_{44}|\leq C\|\partial_{x}u\|_{L^2}^2+C\varepsilon^{4-2\alpha}\|\partial_{y}u\|_{L^2}^2+C\varepsilon^{2-2\alpha+2\lambda}\Big(\int_{\mathbb{T}}\big(|u(x,0)|^2+|u(x,1)|^2\big)\,{\rm d}x\Big).
\end{align*}

In view of \eqref{3.11}, we need to choose $\lambda$ to satisfy
$$
2-2\alpha+2\lambda> 2-\alpha\Rightarrow 2\lambda>\alpha.
$$
So, the above arguments will not work when $\alpha\geq 1$. To overcome this issue, for $\alpha\geq 1$, we decompose $u(x,y)$ as
$$
u(x,y)=\big(u(x,y)-\varepsilon^{\alpha}\partial_{y}u\big)+\varepsilon^{\alpha}\partial_{y}u
$$
near the lower boundary, and
$$
u(x,y)=\big(u(x,y)+\varepsilon^{\alpha}\partial_{y}u\big)-\varepsilon^{\alpha}\partial_{y}u
$$
near the upper boundary. Noting the boundary conditions in \eqref{3.8}, we obtain from Proposition \ref{prop3} and Lemmas \ref{Hardy}--\ref{Hardy3} that for $\alpha\geq 1$,
\begin{align}\label{3.34-1}
	|H_{43}|
	&\leq C(\varepsilon^{\frac{1}{\tilde{\alpha}}}+\sigma)\Big[\varepsilon^2\iintxy\Big[\frac{|u-\varepsilon^{\alpha}\partial_{y}u|^2}{y^2}\chi(y)+\frac{|u+\varepsilon^{\alpha}\partial_{y}u|^2}{(y-1)^2}\tilde{\chi}(y)\Big]\,{\rm d}y{\rm d}x+\varepsilon^{2\alpha}\|\partial_{y}u\|_{L^2}^2\Big]\nonumber\\
	&\leq C \varepsilon^2(\varepsilon^{\frac{1}{\tilde{\alpha}}}+\sigma)\|\partial_{y}u\|_{L^2}^2+C \varepsilon^{2+2\alpha}(\varepsilon^{\frac{1}{\tilde{\alpha}}}+\sigma)\|\partial_{yy}u\|_{L^2}^2,
\end{align}
and
\begin{align}\label{3.34-2}
|H_{44}|&\leq C\varepsilon(\varepsilon^{\frac{1}{\tilde{\alpha}}}+\sigma)\Big[\iintxy\Big\vert\frac{u-\varepsilon^{\alpha}\partial_{y}u}{y}\frac{v}{y}\chi(y)\Big\vert\,{\rm d}y{\rm d}x+\iintxy\Big\vert\frac{u+\varepsilon^{\alpha}\partial_{y}u}{y-1}\frac{v}{y-1}\tilde{\chi}(y)\Big\vert\,{\rm d}y{\rm d}x\Big]\nonumber\\
&\quad +C\varepsilon^{\alpha}(\varepsilon^{\frac{1}{\tilde{\alpha}}}+\sigma)\Big[\iintxy \Big\vert\frac{v\partial_{y}u}{y}\Big\vert\,{\rm d}y{\rm d}x+\iintxy \Big\vert\frac{v\partial_{y}u}{1-y}\Big\vert\,{\rm d}y{\rm d}x\Big]\nonumber\\
&\leq C\varepsilon^2(\varepsilon^{\frac{1}{\tilde{\alpha}}}+\sigma)\iintxy\Big\vert\frac{u-\varepsilon^{\alpha}\partial_{y}u}{y}\Big\vert^2\,{\rm d}y{\rm d}x+C\Big[\iintxy \Big\vert\frac{v}{y}\Big\vert^2\,{\rm d}y{\rm d}x+\iintxy \Big\vert\frac{v}{y-1}\Big\vert^2\,{\rm d}y{\rm d}x\Big]\nonumber\\
&\quad +C\varepsilon^{2\alpha}(\varepsilon^{\frac{1}{\tilde{\alpha}}}+\sigma)\|\partial_{y}u\|_{L^2}^2+C\varepsilon^{2}(\varepsilon^{\frac{1}{\tilde{\alpha}}}+\sigma)\iintxy\Big\vert\frac{u+\varepsilon^{\alpha}\partial_{y}u}{y-1}\Big\vert^2\,{\rm d}y{\rm d}x\nonumber\\
&\leq C\varepsilon^2(\varepsilon^{\frac{1}{\tilde{\alpha}}}+\sigma)\|\partial_{y}u\|_{L^2}^2+C\varepsilon^{2+2\alpha}(\varepsilon^{\frac{1}{\tilde{\alpha}}}+\sigma)\|\partial_{yy}u\|_{L^2}^2+C\|u_{x}\|_{L^2}^2.
\end{align}
For the terms involving $h_{u}$ in $H_{4}$, a direct calculation shows that
\begin{align}\label{3.44-4}
	\varepsilon^{\frac{k_{0}}{\tilde{\alpha}}}\int_{\mathbb{T}}\int_{0}^{1}(u^2\partial_{x}h_{u}+uv\partial_{y}h_{u})\,{\rm d}y{\rm d}x&\leq C\big(\|\partial_{x}u\|_{L^2}^2+\|\partial_{x}v\|_{L^2}^2\big)\nonumber\\
    &\quad +C\varepsilon^{\frac{2k_{0}}{\tilde{\alpha}}}\Big(\|\partial_{y}u\|_{L^2}^2+\int_{\mathbb{T}}|u(x,0)|^2\,{\rm d}x\Big).
\end{align}
In conclusion, combining \eqref{3.26}--\eqref{3.44-4} together, we have
\begin{align}\label{3.52}
|H_{4}|&\leq C\big(\|\partial_{x}u\|_{L^2}^2+\|\partial_{x}v\|_{L^2}^2\big)+C\varepsilon^{4-2\alpha}\|\partial_{y}u\|_{L^2}^2\nonumber\\
&\quad +C\varepsilon^{2-2\alpha+2\lambda}\Big(\int_{\mathbb{T}}\big(|u(x,0)|^2+|u(x,1)|^2\big)\,{\rm d}x\Big)
\end{align}
when $0\leq \alpha<1$, and
\begin{align}\label{3.52-1}
|H_{4}|&\leq C\big(\|\partial_{x}u\|_{L^2}^2+\|\partial_{x}v\|_{L^2}^2\big)+C\varepsilon^{2}(\varepsilon^{\frac{1}{\tilde{\alpha}}}+\sigma)\|\partial_{y}u\|_{L^2}^2+C\varepsilon^{2+2\alpha}(\varepsilon^{\frac{1}{\tilde{\alpha}}}+\sigma)\|\partial_{yy}u\|_{L^2}^2\nonumber\\
&\quad +C\varepsilon^{2}(\varepsilon^{\frac{1}{\tilde{\alpha}}}+\sigma)\int_{\mathbb{T}}|u(x,0)|^2\,{\rm d}x
\end{align}
when $\alpha\geq 1$. Here, $\lambda\in (\frac{\alpha}{2},\frac{1}{2})$ in \eqref{3.52}, and we have chosen $k_{0}\geq (2-\alpha)\tilde{\alpha}$.

Substituting \eqref{3.12}, \eqref{3.12-1}, \eqref{3.13}, \eqref{3.14}, \eqref{3.16}--\eqref{3.16-1}, and \eqref{3.52}--\eqref{3.52-1} into \eqref{3.11}, and taking $\varepsilon$, $\tau$ and $\sigma$ small enough (when $\alpha\geq 1$), we obtain
\begin{align}\label{3.11-2}
&\varepsilon^2\int_{\mathbb{T}}\int_{0}^{1}|\nabla u|^2+|\nabla v|^2\,{\rm d}y{\rm d}x+\varepsilon^{2-\alpha}\int_{\mathbb{T}}|u(x,0)|^2+|u(x,1)|^2\,{\rm d}x\nonumber\\
	&\leq C\big(\|\partial_{x}u\|_{L^2}^2+\|\partial_{x}v\|_{L^2}^2\big)+C\varepsilon^{-2}\|(R_{u},R_{v})\|_{L^2}^2
\end{align}
when $0\leq \alpha<1$, and
\begin{align}\label{3.11-1}
	&\varepsilon^2\int_{\mathbb{T}}\int_{0}^{1}|\nabla u|^2+|\nabla v|^2\,{\rm d}y{\rm d}x+\varepsilon^{2-\alpha}\int_{\mathbb{T}}|u(x,0)|^2+|u(x,1)|^2\,{\rm d}x\nonumber\\
	&\leq C\big(\|\partial_{x}u\|_{L^2}^2+\|\partial_{x}v\|_{L^2}^2\big)+C\varepsilon^{2+2\alpha}(\varepsilon^{\frac{1}{\tilde{\alpha}}}+\sigma)\|\partial_{yy}u\|_{L^2}^2+C\varepsilon^{-2}\|(R_{u},R_{v})\|_{L^2}^2
\end{align}
when $\alpha\geq 1$.

\subsection{Positivity estimates} To control $\|(\partial_{x}u,\partial_{x}v)\|_{L^2}^2$, we multiply $\eqref{3.8}_{1}$ by $\partial_{x}u$ and $\eqref{3.8}_{2}$ by $\partial_{x}v$ and integrate the resultant equation over $\mathbb{T}\times (0,1)$ to get
\begin{align}\label{3.55}
	&\int_{\mathbb{T}}\int_{0}^{1}u^{a}(|\partial_{x}u|^{2}+|\partial_{x}v|^2)\,{\rm d}y{\rm d}x\nonumber\\
	&=-\int_{\mathbb{T}}\int_{0}^{1}\big(v^{a}\partial_{y}u\partial_{x}u+u\partial_{x}u\partial_{x}u^{a}+v\partial_{x}u\partial_{y}u^{a}\big)\,{\rm d}y{\rm d}x\nonumber\\
	&\quad-\int_{\mathbb{T}}\int_{0}^{1}\big(v^{a}\partial_{y}v\partial_{x}v+u\partial_{x}v^{a}\partial_{x}v+v\partial_{y}v^{a}\partial_{x}v\big)\,{\rm d}y{\rm d}x\nonumber\\
	&\quad -\int_{\mathbb{T}}\int_{0}^{1}\big(\partial_{x}p\partial_{x}u+\partial_{y}p\partial_{x}v\big)\,{\rm d}y{\rm d}x+\int_{\mathbb{T}}\int_{0}^{1}\big(\varepsilon^2\Delta u\partial_{x}u+\varepsilon^2\Delta v\partial_{x}v\big)\,{\rm d}y{\rm d}x\nonumber\\
	&\quad +\int_{\mathbb{T}}\int_{0}^{1}\big(R_{u}\partial_{x}u+R_{v}\partial_{x}v\big)\,{\rm d}y{\rm d}x=\sum\limits_{i=1}^{5}Q_{i}.
\end{align}

Recalling \eqref{bd}, without loss of generality, we assume $\bar{f}_{0}$ and $\bar{f}_{1}$ are both positive. By the choices of $A$ and $B$ in \eqref{M00}, \eqref{M0}, \eqref{M2} and \eqref{M1}, there exists a positive constant $c_{0}$ depending on $\bar{f}_{0}$ and $\bar{f}_{1}$ such that
$$
u_{I}^{0}=Ay+B\geq c_{0}>0.
$$
Noting $u_{B}^{0}=0$ for $0\leq \alpha<1$ and $u_{B}^{0}=O(\sigma)$ for $\alpha\geq 1$, we have
$$
u^{a}\geqslant \frac{1}{2}c_{0}>0,
$$
provided that $\varepsilon$ and $\sigma$ are sufficiently small. Then
\begin{align}\label{3.56}
	\int_{\mathbb{T}}\int_{0}^{1}u^{a}\big(|\partial_{x}u|^2+|\partial_{x}v|^2\big)\,{\rm d}y{\rm d}x\geq \frac{1}{2}c_{0}\int_{\mathbb{T}}\int_{0}^{1}\big(|\partial_{x}u|^2+|\partial_{x}v|^2\big)\,{\rm d}y{\rm d}x.
\end{align}
For $Q_{3}$, it follows from integration by parts and the divergence free condition that
\begin{align}\label{3.57}
	Q_{3}=&-\int_{\mathbb{T}}\int_{0}^{1}\big(\partial_{x}p\partial_{x}u+\partial_{y}p\partial_{x}v\big)\,{\rm d}y{\rm d}x=\int_{\mathbb{T}}\int_{0}^{1}p\partial_{x}\big(\partial_{x}u+\partial_{y}v\big)\,{\rm d}y{\rm d}x=0.
\end{align}
For $Q_{4}$, a direct calculation shows that
\begin{align}\label{3.58}
	Q_{4}&=\varepsilon^2\int_{\mathbb{T}}\int_{0}^{1}\big(\Delta u\partial_{x}u+\Delta v\partial_{x}v\big)\,{\rm d}y{\rm d}x\nonumber\\
	&=\varepsilon^2\int_{\mathbb{T}}\int_{0}^{1}\big(\partial_{yy}u\partial_{x}u+\partial_{yy}v\partial_{x}v\big)\,{\rm d}y{\rm d}x\nonumber\\
	&=\varepsilon^2\int_{\mathbb{T}}\big(\partial_{x}u\partial_{y}u\big)(x,1)\,{\rm d}x-\varepsilon^2\int_{\mathbb{T}}\big(\partial_{x}u\partial_{y}u\big)(x,0)\,{\rm d}x\nonumber\\
	&\quad-\varepsilon^2\int_{\mathbb{T}}\int_{0}^{1}\frac{1}{2}\partial_{x}(\partial_{y}u)^2\,{\rm d}y{\rm d}x-\varepsilon^2\int_{\mathbb{T}}\int_{0}^{1}\frac{1}{2}\partial_{x}(\partial_{y}v)^2\,{\rm d}y{\rm d}x\nonumber\\
	&=-\varepsilon^{2-\alpha}\int_{\mathbb{T}}\big(u\partial_{x}u\big)(x,1)\,{\rm d}x-\varepsilon^{2-\alpha}\int_{\mathbb{T}}\big(u\partial_{x}u\big)(x,0)\,{\rm d}x=0.
\end{align}
For $Q_{5}$, by the H\"{o}lder inequality, one has
\begin{align}\label{3.57-1}
Q_{5}\leq \tau(\|\partial_{x}u\|_{L^2}^2+\|\partial_{x}v\|_{L^2}^2)+C_{\tau}\big(\|R_{u}\|_{L^2}^2+\|R_{v}\|_{L^2}^2\big).
\end{align}

For $Q_{1}$ and $Q_{2}$, it follows from Proposition \ref{prop3} and Lemma \ref{lem5.1} that
$$
|v^{a}|+|\partial_{x}v^{a}|\leq \left\{
\begin{aligned}
&C\varepsilon^{(2-\alpha)},\qquad\,\,\,\text{for $0\leq \alpha<1$},\\
&C\varepsilon(\varepsilon^{\frac{1}{\tilde{\alpha}}}+\sigma),\quad\text{for $\alpha\geq 1$},
\end{aligned}
\right.\quad |\partial_{y}v^{a}|\leq \left\{
\begin{aligned}
&C\varepsilon^{(1-\alpha)},\qquad\text{for $0\leq \alpha<1$},\\
&C(\varepsilon^{\frac{1}{\tilde{\alpha}}}+\sigma),\quad\text{for $\alpha\geq 1$}.
\end{aligned}
\right.
$$
Then we have
\begin{align}\label{3.59}
	&-\int_{\mathbb{T}}\int_{0}^{1}\big(v^{a}\partial_{y}u\partial_{x}u+v^{a}\partial_{y}v\partial_{x}v\big)\,{\rm d}y{\rm d}x\nonumber\\
    &\leq \left\{
\begin{aligned}
    &\tau\big(\|\partial_{x}u\|_{L^2}^2+\|\partial_{x}v\|_{L^2}^2\big)+C_{\tau}\varepsilon^{4-2\alpha}\big(\|\partial_{y}u\|_{L^2}^2+\|\partial_{y}v\|_{L^2}^2\big),\qquad\qquad\,\,\text{for $0\leq \alpha<1$},\\
    &\tau\big(\|\partial_{x}u\|_{L^2}^2+\|\partial_{x}v\|_{L^2}^2\big)+C_{\tau}\varepsilon^{2}(\varepsilon^{\frac{1}{\tilde{\alpha}}}+\sigma)\big(\|\partial_{y}u\|_{L^2}^2+\|\partial_{y}v\|_{L^2}^2\big),\qquad \text{for $\alpha\geq 1$},
\end{aligned}
    \right.
\end{align}
\begin{align}\label{3.60}
	-\int_{\mathbb{T}}\int_{0}^{1}v\partial_{y}v^{a}\partial_{x}v\,{\rm d}y{\rm d}x\leq \left\{
\begin{aligned}
    &C\varepsilon^{1-\alpha}\|\partial_{x}v\|_{L^2}^2,\qquad\qquad \text{for $0\leq \alpha<1$},\\
    &C(\varepsilon^{\frac{1}{\tilde{\alpha}}}+\sigma)\|\partial_{x}v\|_{L^2}^2,\qquad\,\,\text{for $\alpha\geq 1$},
\end{aligned}
    \right.
\end{align}
and
\begin{align}\label{3.61}
	&-\int_{\mathbb{T}}\int_{0}^{1}u\partial_{x}v^{a}\partial_{x}v\,{\rm d}y{\rm d}x\nonumber\\
    &\leq 
    \left\{
    \begin{aligned}
        &\tau\|\partial_{x}v\|_{L^2}^2+C_{\tau}\varepsilon^{4-2\alpha}\Big(\|\partial_{y}u\|_{L^2}^2+\|\partial_{x}u\|_{L^2}^2+\int_{\mathbb{T}}|u(x,0)|^2\,{\rm d}x\Big),\qquad\quad\,\text{for $0\leq \alpha<1$},\\
        &\tau\|\partial_{x}v\|_{L^2}^2+C_{\tau}\varepsilon^{2}(\varepsilon^{\frac{1}{\tilde{\alpha}}}+\sigma)\Big(\|\partial_{y}u\|_{L^2}^2+\|\partial_{x}u\|_{L^2}^2+\int_{\mathbb{T}}|u(x,0)|^2\,{\rm d}x\Big),\quad \text{for $\alpha\geq 1$},
    \end{aligned}
    \right.
\end{align}
where we have used Lemma \ref{So1} in \eqref{3.61}.

For $\int_{\mathbb{T}}\int_{0}^{1}u\partial_{x}u\partial_{x}u^{a}\,{\rm d}y{\rm d}x$, we split it into three parts:
\begin{align}\label{9.1}
	\iintxy u\partial_{x}u\partial_{x}u^{a}\,{\rm d}y{\rm d}x&=\iintxy u\partial_{x}u\partial_{x}u_{I}^{a}\,{\rm d}y{\rm d}x+\iintxy u\partial_{x}u\partial_{x}(u_{B}^{a}+\hat{u}_{B}^{a})\,{\rm d}y{\rm d}x\quad \nonumber\\
    &\quad +\varepsilon^{\frac{k_{0}}{\tilde{\alpha}}}\iintxy u\partial_{x}u\partial_{x}h_{u}\,{\rm d}y{\rm d}x\nonumber\\
	&=I_{1}+I_{2}+I_{3}.
\end{align}
For $I_{1}$ and $I_{3}$, it follows from Lemma \ref{lem5.1} that
\begin{align*}
	|\partial_{x}u_{I}^{a}|\leq
    \left\{
\begin{aligned}
&C\varepsilon^{2-\alpha},\qquad\quad\,\,\text{for $0\leq \alpha<1$},\\
&C\varepsilon(\varepsilon^{\frac{1}{\tilde{\alpha}}}+\sigma),\quad \text{for $\alpha\geq 1$},
\end{aligned}
    \right.
\end{align*}
then
\begin{align}\label{9.2}
	&|I_{1}|+|I_{3}|\nonumber\\
    &\leq \left\{
\begin{aligned}
    &\tau\|\partial_{x}u\|_{L^2}^2+C_{\tau}\varepsilon^{4-2\alpha}\Big(\|\partial_{x}u\|_{L^2}^2+\|\partial_{y}u\|_{L^2}^2+\int_{\mathbb{T}}|u(x,0)|^2\,{\rm d}x\Big),\qquad\quad\,\text{for $0\leq \alpha<1$},\\
    &\tau\|\partial_{x}u\|_{L^2}^2+C_{\tau}\varepsilon^{2}(\varepsilon^{\frac{1}{\tilde{\alpha}}}+\sigma)\Big(\|\partial_{x}u\|_{L^2}^2+\|\partial_{y}u\|_{L^2}^2+\int_{\mathbb{T}}|u(x,0)|^2\,{\rm d}x\Big),\quad \text{for $\alpha\geq 1$},
\end{aligned}
\right.
\end{align}
provided that $k_{0}\geq (2-\alpha)\tilde{\alpha}$ for $0\leq \alpha<1$, and $k_{0}\geq \tilde{\alpha}$ for $\alpha\geq 1$.

For $I_{2}$, when $0\leq \alpha<1$, using similar arguments as in the estimation of $H_{43}$ in \eqref{3.34-0}, we have
\begin{align*}
|I_{2}|\leq \tau\|\partial_{x}u\|_{L^2}^2+C_{\tau}\varepsilon^{4-2\alpha}\|\partial_{y}u\|_{L^2}^2+C_{\tau}\varepsilon^{2-2\alpha+2\lambda}\Big(\int_{\mathbb{T}}|u(x,0)|^2+|u(x,1)|^2\,{\rm d}x\Big).
\end{align*}
When $\alpha\geq 1$, by using similar decomposition and arguments as in \eqref{3.34-1}, we have
\begin{align}\label{9.6}
	|I_{2}|\leq \tau\|\partial_{x}u\|_{L^2}^2+C_{\tau}\big(\varepsilon^2(\varepsilon^{\frac{1}{\tilde{\alpha}}}+\sigma)\|\partial_{y}u\|_{L^2}^2+\varepsilon^{2+2\alpha}(\varepsilon^{\frac{1}{\tilde{\alpha}}}+\sigma)\|\partial_{yy}u\|_{L^2}^2\big).
\end{align}
Therefore, substituting \eqref{9.2}--\eqref{9.6} into \eqref{9.1}, one has
\begin{align*}
\iintxy u\partial_{x}u\partial_{x}u^{a}\,{\rm d}y{\rm d}x\leq \tau\|\partial_{x}u\|_{L^2}^2+C_{\tau}\Big(\varepsilon^{4-2\alpha}\|\partial_{y}u\|_{L^2}^2+\varepsilon^{2-2\alpha+2\lambda}\int_{\mathbb{T}}|u(x,0)|^2+|u(x,1)|^2\,{\rm d}x\Big),
\end{align*}
provided that $0\leq \alpha<1$, and
\begin{align*}
\iintxy u\partial_{x}u\partial_{x}u^{a}\,{\rm d}y{\rm d}x\leq \tau\|\partial_{x}u\|_{L^2}^2+C_{\tau}\varepsilon^2(\varepsilon^{\frac{1}{\tilde{\alpha}}}+\sigma)\Big(\|\partial_{y}u\|_{L^2}^2+\int_{\mathbb{T}}|u(x,0)|^2\,{\rm d}x+\varepsilon^{2\alpha}\|\partial_{yy}u\|_{L^2}^2\Big),
\end{align*}
provided that $\alpha\geq 1$.

For the remaining term $\iintxy v\partial_{x}u\partial_{y}u^{a}\,{\rm d}y{\rm d}x$, using $\partial_{x}u=-\partial_{y}v$, we also split it into three parts
\begin{align}\label{9.7}
	\iintxy v\partial_{x}u\partial_{y}u^{a}\,{\rm d}y{\rm d}x&=-\iintxy v\partial_{y}v\partial_{y}u_{I}^{a}\,{\rm d}y{\rm d}x-\iintxy v\partial_{y}v\partial_{y}(u_{B}^{a}+\hat{u}_{B}^{a})\,{\rm d}y{\rm d}x\nonumber\\
    &\quad -\varepsilon^{\frac{k_{0}}{\tilde{\alpha}}}\iintxy v\partial_{y}v\partial_{y}h_{u}\,{\rm d}y{\rm d}x\nonumber\\
	&=J_{1}+J_{2}+J_{3}.
\end{align}
For $J_{1}$ and $J_{3}$, it follows from Lemma \ref{lem5.1} that
$$
\partial_{y}u_{I}^{i}=\tilde{C}_{i},\quad \text{for $0\leq i\leq \tilde{\alpha}-1$},
$$
we have
\begin{align}\label{9.8}
|J_{1}|+|J_{3}|\leq C\varepsilon\|v\|_{L^2}\|\partial_{y}v\|_{L^2}\leq C\varepsilon\|\partial_{x}v\|_{L^2}^2+C\varepsilon\|\partial_{x}u\|_{L^2}^2,
\end{align}
provided that $k_{0}\geq \tilde{\alpha}$.
For $J_{2}$, it follows from Proposition \ref{prop3} and Lemma \ref{Hardy} that
\begin{align}\label{9.10}
	|J_{2}|\leq \left\{
    \begin{aligned}
    &C\varepsilon^{1-\alpha}\|\partial_{y}v\|_{L^2}^2= C\varepsilon^{1-\alpha}\|\partial_{x}u\|_{L^2}^2,\qquad\qquad\quad\,\,\,\text{for $0\leq \alpha<1$},\\
    &C(\varepsilon^{\frac{1}{\tilde{\alpha}}}+\sigma)\|\partial_{y}v\|_{L^2}^2=C(\varepsilon^{\frac{1}{\tilde{\alpha}}}+\sigma)\|\partial_{x}u\|_{L^2}^2,\qquad\text{for $\alpha\geq 1$}.
    \end{aligned}
    \right.
\end{align}
Substituting \eqref{9.8}--\eqref{9.10} into \eqref{9.7}, one has
\begin{align}\label{9.10-1}
\iintxy v\partial_{x}u\partial_{y}u^{a}\,{\rm d}y{\rm d}x\leq \left\{
\begin{aligned}
&C\varepsilon^{1-\alpha}(\|\partial_{x}u\|_{L^2}^2+\|\partial_{x}v\|_{L^2}^2),\qquad\qquad\,\text{for $0\leq \alpha<1$},\\
&C\varepsilon\|\partial_{x}v\|_{L^2}^2+C(\varepsilon^{\frac{1}{\tilde{\alpha}}}+\sigma)\|\partial_{x}u\|_{L^2}^2,\qquad \text{for $\alpha\geq 1$}.
\end{aligned}
\right.
\end{align}
Therefore, combining \eqref{3.59}--\eqref{9.10-1} together, we obtain
\begin{align}\label{3.67-1}
	Q_{1}+Q_{2}&\leq (5\tau+C\varepsilon^{1-\alpha})(\|\partial_{x}u\|_{L^2}^2+\|\partial_{x}v\|_{L^2}^2)+C_{\tau}\varepsilon^{4-2\alpha}\|\partial_{y}u\|_{L^2}^2\nonumber\\
    &\quad +C_{\tau}\varepsilon^{2-2\alpha+2\lambda}\Big(\int_{\mathbb{T}}|u(x,0)|^2+|u(x,1)|^2\,{\rm d}x\Big)
\end{align}
provided that $0\leq \alpha<1$, and
\begin{align}\label{3.67-2}
	Q_{1}+Q_{2}&\leq (5\tau+C(\varepsilon^{\frac{1}{\tilde{\alpha}}}+\sigma))(\|\partial_{x}u\|_{L^2}^2+\|\partial_{x}v\|_{L^2}^2)\nonumber\\
    &\quad +C_{\tau}\varepsilon^2(\varepsilon^{\frac{1}{\tilde{\alpha}}}+\sigma)\Big(\|\partial_{y}u\|_{L^2}^2+\int_{\mathbb{T}}|u(x,0)|^2\,{\rm d}x+\varepsilon^{2\alpha}\|\partial_{yy}u\|_{L^2}^2\Big)
\end{align}
provided that $\alpha\geq 1$. Substituting \eqref{3.57}, \eqref{3.58}, \eqref{3.57-1}, and \eqref{3.67-1}--\eqref{3.67-2} into \eqref{3.55}, and using \eqref{3.11-2}--\eqref{3.11-1}, we obtain
\begin{align}\label{3.68-1}
	&\|\partial_{x}u\|_{L^2}^2+\|\partial_{x}v\|_{L^2}^2+\varepsilon^2\|\partial_{y}u\|_{L^2}^2+\varepsilon^{2-\alpha}\int_{\mathbb{T}}|u(x,0)|^2+|u(x,1)|^2\,{\rm d}x\leq  C\varepsilon^{-2}\|(R_{u},R_{v})\|_{L^2}^2
\end{align}
provided that $0\leq \alpha<1$ and $\tau$, $\varepsilon$ are small enough, and 
\begin{align}\label{3.69-1}
&\|\partial_{x}u\|_{L^2}^2+\|\partial_{x}v\|_{L^2}^2+\varepsilon^2\|\partial_{y}u\|_{L^2}^2+\varepsilon^{2-\alpha}\int_{\mathbb{T}}|u(x,0)|^2+|u(x,1)|^2\,{\rm d}x\nonumber\\
&\leq C\varepsilon^{-2}\|(R_{u},R_{v})\|_{L^2}^2+C\varepsilon^{2+2\alpha}(\varepsilon^{\frac{1}{\tilde{\alpha}}}+\sigma)\|\partial_{yy}u\|_{L^2}^2
\end{align}
provided that $\alpha\geq 1$ and $\tau$, $\varepsilon$ and $\sigma$ are small enough.

\subsection{Stokes' estimates}
We still need to derive the second order derivatives estimates, especially for the estimates of $\|\partial_{yy}u\|_{L^2}^2$ when $\alpha\geq 1$. The proof is divided into four steps.

\textbf{Step 1. Tangential second derivative estimates}. We rewrite the remainder equation \eqref{3.8} as
\begin{align}\label{9.12}
	\left\{
	\begin{aligned}
	&-\varepsilon^2\Delta U+\nabla p=F,\\
	&\operatorname{div}U=0,	
	\end{aligned}
	\right.
\end{align}
where $U=(u,v)$, $U^{a}=(u^{a},v^{a})$ and $F=\mathfrak{R}-\mathfrak{N}(U,U^{a})$ with
$\mathfrak{R}=(R_{u},R_{v})$ and
\begin{align*}
	\mathfrak{N}(U,U^{a})=\left(
	\begin{matrix}
		&u^{a}u_{x}+v^{a}u_{y}+uu_{x}^{a}+vu_{y}^{a}\\
		&u^{a}v_{x}+v^{a}v_{y}+uv_{x}^{a}+vv_{y}^{a}
	\end{matrix}
	\right).
\end{align*}
Differentiating \eqref{9.12} with respect to $x$ to obtain
\begin{align}\label{9.13}
	\left\{
	\begin{aligned}
		&-\varepsilon^2\Delta U_{x}+\nabla p_{x}=F_{x},\\
		&\operatorname{div}U_{x}=0,	
	\end{aligned}
	\right.
\end{align}
supplemented with the boundary conditions
$$
v_{x}(x,0)=v_{x}(x,1)=0,\quad (u_{xy}-\varepsilon^{-\alpha}u_{x})(x,0)=0,\quad (u_{xy}+\varepsilon^{-\alpha}u_{x})(x,1)=0.
$$
Multiplying \eqref{9.13} by $U_{x}$ to see
\begin{align}\label{9.14}
	\varepsilon^2\|\nabla U_{x}\|_{L^2}^2+\varepsilon^{2-\alpha}\int_{\mathbb{T}}(|u_{x}(x,0)|^2+|u_{x}(x,1)|^2)\,{\rm d}x=\iintxy F_{x}\cdot U_{x}\,{\rm d}y{\rm d}x.
\end{align}
Integrating by parts to get
\begin{align*}
	\Big\vert\iintxy F_{x}\cdot U_{x}\,{\rm d}y{\rm d}x\Big\vert=\Big\vert-\iintxy F\cdot U_{xx}\,{\rm d}y{\rm d}x\Big\vert\leq \|F\|_{L^2}\|U_{xx}\|_{L^2}\leq \|F\|_{L^2}\|\nabla U_{x}\|_{L^2}.
\end{align*}
Then we obtain from \eqref{9.14} that
\begin{align}\label{9.14-1}
\varepsilon^4\|\nabla U_{x}\|_{L^2}^2+\varepsilon^{4-\alpha}\int_{\mathbb{T}}(|u_{x}(x,0)|^2+|u_{x}(x,1)|^2)\,{\rm d}x\leq C\|F\|_{L^2}^2.
\end{align}

\textbf{Step 2. Estimate of $p_{x}$}. Since $p$ is periodic in $x$, we have
\begin{align*}
	\iintxy p_{x}\,{\rm d}y{\rm d}x=0.
\end{align*}
Then there exists a vector function $\Phi=(\phi,\psi)\in H_{0}^{1}(\mathbb{T}\times (0,1))$ such that
\begin{align*}
	\operatorname{div}\Phi=p_{x},\quad \|\Phi\|_{H^{1}(\mathbb{T}\times (0,1))}\leq C\|p_{x}\|_{L^{2}(\mathbb{T}\times (0,1))}.
\end{align*}
Now, we use $\Phi$ as a test function to test $\partial_{x}\eqref{9.12}$ to obtain
\begin{align*}
	\varepsilon^2\iintxy \nabla U_{x}:\nabla \Phi\,{\rm d}y{\rm d}x-\iintxy p_{x}\operatorname{div}\Phi\,{\rm d}y{\rm d}x=\langle F_{x},\Phi\rangle.
\end{align*}
Noting that $\operatorname{div}\Phi=p_{x}$, one has from the above equality that
\begin{align}\label{9.15}
	\|p_{x}\|_{L^2}^2\leq \varepsilon^2\Big\vert\iintxy \nabla U_{x}:\nabla \Phi\,{\rm d}y{\rm d}x\Big\vert+|\langle F_{x},\Phi\rangle|.
\end{align}
Using \eqref{9.14-1}, a direct calculation shows that
\begin{align}\label{9.16}
	\varepsilon^2\Big\vert\iintxy \nabla U_{x}:\nabla \Phi\,{\rm d}y{\rm d}x\Big\vert\leq \varepsilon^2\|\nabla U_{x}\|_{L^2}\|\nabla \Phi\|_{L^2}\leq C\|F\|_{L^2}\|p_{x}\|_{L^2},
\end{align}
and
\begin{align}\label{9.17}
	|\langle F_{x},\Phi\rangle|=\Big\vert-\iintxy F\cdot \Phi_{x}\,{\rm d}y{\rm d}x\Big\vert\leq C\|F\|_{L^2}\|p_{x}\|_{L^2}.
\end{align}
Therefore, it holds that
\begin{align}\label{9.18}
	\|p_{x}\|_{L^2}\leq C\|F\|_{L^2}.
\end{align}

\textbf{Step 3. Estimate of $p_{y}$}. From the second momentum equation of \eqref{9.12}, one has
$$
-\varepsilon^2\Delta v+p_{y}=F_{2},
$$
which implies
$$
p_{y}=F_{2}+\varepsilon^2\Delta v=F_{2}+\varepsilon^2(v_{xx}+v_{yy})=F_{2}+\varepsilon^2(v_{xx}-u_{xy}).
$$
Using \eqref{9.14-1}, a direct calculation shows that
$$
\|p_{y}\|_{L^2}\leq C\|F\|_{L^2}+C\varepsilon^2\|\nabla U_{x}\|_{L^2}\leq C\|F\|_{L^2},
$$
which, together with \eqref{9.18}, yields that
\begin{align}\label{9.18-1}
\|\nabla p\|_{L^2}\leq C\|F\|_{L^2}.
\end{align}

\textbf{Step 4. Estimates of $u_{yy}$}. From the first component of \eqref{9.12}, one has
\begin{align}\label{9.19}
	\varepsilon^2u_{yy}=p_{x}-F_{1}-\varepsilon^2u_{xx},
\end{align}
which, together with \eqref{9.14-1} and \eqref{9.18}, implies that
\begin{align}\label{9.19-1}
\varepsilon^4\|u_{yy}\|_{L^2}^2\leq c\|F_{1}\|_{L^2}^2+C\|p_{x}\|_{L^2}^2+C\varepsilon^4\|u_{xx}\|_{L^2}^2\leq C\|F\|_{L^2}^2.
\end{align}
Combining \eqref{9.14-1}, \eqref{9.18-1} and \eqref{9.19-1} and using $v_{yy}=-u_{xy}$, we obtain the following Stokes' estimates
\begin{align}\label{9.20}
	\varepsilon^4\|\nabla^2 U\|_{L^2}^2+\|\nabla p\|_{L^2}^2\leq C\|F\|_{L^2}^2.
\end{align}

Using \eqref{3.68-1}--\eqref{3.69-1}, Proposition \ref{prop3}, Lemma \ref{lem5.1} and the H\"{o}lder and Hardy inequalities, we have 
\begin{align}\label{9.20-1}
\|F\|_{L^2}^2&\leq C\|(R_{u},R_{v})\|_{L^2}^2+C\big(\|u^{a}(u_{x},v_{x})\|_{L^2}^2+\|v^{a}(u_{y},v_{y})\|_{L^2}^2+\|u(u_{x}^{a},v_{x}^{a})\|_{L^2}^2+\|v(u_{y}^{a},v_{y}^{a})\|_{L^2}^2\big)\nonumber\\
&\leq C\|(R_{u},R_{v})\|_{L^2}^2+C\|(u_{x},v_{x})\|_{L^2}^2+C\varepsilon^{4-2\alpha}\|(u_{y},v_{y})\|_{L^2}^2+C\|v\|_{L^2}^2\nonumber\\
&\quad +C\varepsilon^{2-2\alpha+2\lambda}(\overline{u^2(x,0)}+\overline{u^2(x,1)})+C\varepsilon^{2-2\alpha}\|\partial_{y}v\|_{L^2}^2\nonumber\\
&\leq C\varepsilon^{-2}\|(R_{u},R_{v})\|_{L^2}^2,
\end{align}
for any $\frac{\alpha}{2}<\lambda<\frac{1}{2}$ provided that $0\leq \alpha<1$, and
\begin{align}\label{9.20-2}
\|F\|_{L^2}^2&\leq C\|(R_{u},R_{v})\|_{L^2}^2+C\big(\|u^{a}(u_{x},v_{x})\|_{L^2}^2+\|v^{a}(u_{y},v_{y})\|_{L^2}^2+\|u(u_{x}^{a},v_{x}^{a})\|_{L^2}^2+\|v(u_{y}^{a},v_{y}^{a})\|_{L^2}^2\big)\nonumber\\
&\leq C\|(R_{u},R_{v})\|_{L^2}^2+C\|(u_{x},v_{x})\|_{L^2}^2+C\varepsilon^2(\varepsilon^{\frac{1}{\tilde{\alpha}}}+\sigma)\|(u_{y},v_{y})\|_{L^2}^2+C\|v\|_{L^2}^2\nonumber\\
&\quad +C\varepsilon^{2+2\alpha}(\varepsilon^{\frac{1}{\tilde{\alpha}}}+\sigma)\|\partial_{yy}u\|_{L^2}^2+C(\varepsilon^{\frac{1}{\tilde{\alpha}}}+\sigma)\|\partial_{y}v\|_{L^2}^2\nonumber\\
&\leq C\varepsilon^{-2}\|(R_{u},R_{v})\|_{L^2}^2+C\varepsilon^{2+2\alpha}(\varepsilon^{\frac{1}{\tilde{\alpha}}}+\sigma)\|\partial_{yy}u\|_{L^2}^2,
\end{align}
provided that $\alpha\geq 1$. Substituting \eqref{9.20-1}--\eqref{9.20-2} into \eqref{9.20}, and combining \eqref{3.68-1}--\eqref{3.69-1}, we obtain that
\begin{align}\label{9.22}
	&\|(u_{x},v_{x},v_{y})\|_{L^2}^2+\varepsilon^2\|u_{y}\|_{L^2}^2+\varepsilon^4\|\nabla^2(u,v)\|_{L^2}^2+\varepsilon^{2-\alpha}\int_{\mathbb{T}}(|u(x,0)|^2+|u(x,1)|^2)\,{\rm d}x\nonumber\\
    &\leq C\varepsilon^{-2}\|(R_{u},R_{v})\|_{L^2}^2.
\end{align}
provided that $\sigma$ and $\varepsilon$ are sufficiently small when $\alpha\geq 1$.

\subsection{Existence of remainder equation}. With the uniform {\it a priori} estimates \eqref{9.22}, we shall establish the existence of the remainder equation \eqref{3.8} through the contraction mapping theorem. 

\begin{proposition}\label{prop4}
    Let $k_{0}\geq 7\tilde{\alpha}$ be an integer. There exists $\varepsilon_{0}>0$ and $\sigma_{0}$ ${\rm(}$required only when $\alpha\geq 1{\rm)}$ such that for any $\varepsilon\in (0,\varepsilon_{0})$, and $\sigma\in [0,\sigma_{0})$, the remainder equation \eqref{3.8} has a unique solution $(u,v)$ satisfying
$$
\varepsilon\|\nabla(u,v)\|_{L^2}+\varepsilon^{3}\|\nabla^2(u,v)\|_{L^2}+\varepsilon^{1-\frac{\alpha}{2}}\Big(\int_{\mathbb{T}}(|u(x,0)|^2+|u(x,1)|^2)\,{\rm d}x\Big)^{\frac{1}{2}}\leq C\varepsilon^{6}.
$$
    In particular, it holds that
    \begin{align*}
        \|(u,v)\|_{L^{\infty}}\leq C\varepsilon^3.
    \end{align*}
\end{proposition}

\noindent\textbf{Proof}. For each smooth function $(\tilde{u},\tilde{v})\in H^2(\Omega)\times H^2(\Omega)$ which satisfies
\begin{align}\label{9.23}
\left\{
\begin{aligned}
&\partial_{x}\tilde{u}+\partial_{y}\tilde{v}=0,\\
&\tilde{v}(x,0)=\tilde{v}(x,1)=0,\\
&(\tilde{u}-\varepsilon^{\alpha}\partial_{y}\tilde{u})(x,0)=(\tilde{u}+\varepsilon^{\alpha}\partial_{y}\tilde{u})(x,1)=0,
\end{aligned}
\right.
\end{align}
we consider the following linear system
\begin{align}\label{9.24}
\left\{
\begin{aligned}
&u^{a}\partial_{x}u+v^{a}\partial_{y}u+u\partial_{x}u^{a}+v\partial_{y}u^{a}+\partial_{x}p=\varepsilon^{2}\Delta u+R_{u}^{a}-(\tilde{u}\partial_{x}\tilde{u}+\tilde{v}\partial_{y}\tilde{u}),\\
		&u^{a}\partial_{x}v+v^{a}\partial_{y}v+u\partial_{x}v^{a}+v\partial_{y}v^{a}+\partial_{y}p=\varepsilon^2\Delta v+R_{v}^{a}-(\tilde{u}\partial_{x}\tilde{v}+\tilde{v}\partial_{y}\tilde{v}),\\
		&\partial_{x}u+\partial_{y}v=0,\\
		&v(x,0)=v(x,1)=0,\\
		&(u-\varepsilon^{\alpha}\partial_{y}u)(x,0)=(u+\varepsilon^{\alpha}\partial_{y}u)(x,1)=0.
\end{aligned}
\right.
\end{align}
The solvability of the linear system \eqref{9.24} follows from the standard Lax-Milgram theorem combined with the Garding inequality, Fredholm alternative theorem and the {\it a priori } estimate \eqref{9.22}. Moreover, by \eqref{9.22}, we deduce that there exist $\varepsilon_{0}>0$ and $\sigma_{0}>0$ (required only when $\alpha\geq 1$) such that for $\varepsilon<\varepsilon_{0}$ and $\sigma\in [0,\sigma_{0})$,
\begin{align}\label{9.25}
&\|(u_{x},v_{x},v_{y})\|_{L^2}^2+\varepsilon^2\|u_{y}\|_{L^2}^2+\varepsilon^{4}\|\nabla^2(u,v)\|_{L^2}^2+\varepsilon^{2-\alpha}\int_{\mathbb{T}}\big(|u(x,0)|^2+|u(x,1)|^2)\,{\rm d}x\nonumber\\
&\leq C\varepsilon^{-2}\|(R_{u}^{a},R_{v}^{a})\|_{L^2}^2+C\varepsilon^{-2}\|(\tilde{u},\tilde{v})\cdot \nabla (\tilde{u},\tilde{v})\|_{L^2}^2\nonumber\\
&\leq C\varepsilon^{-2}\|(R_{u}^{a},R_{v}^{a})\|_{L^2}^2+C\varepsilon^{-2}\|(\tilde{u},\tilde{v})\|_{L^{\infty}}^2\|\nabla (\tilde{u},\tilde{v})\|_{L^2}^2.
\end{align}
Let $\mathbb{X}$ be the function space
\begin{align*}
\mathbb{X}=\{(u,v)\in H^{2}(\Omega)\times H^2(\Omega)\,\vert\,(u,v)\text{ satisfies }\eqref{9.23}\}
\end{align*}
equipped with the norm
$$
\|(u,v)\|_{\mathbb{X}}:=\varepsilon\|\nabla(u,v)\|_{L^2}+\varepsilon^{3}\|\nabla^2(u,v)\|_{L^2}+\varepsilon^{1-\frac{\alpha}{2}}\Big(\int_{\mathbb{T}}(|u(x,0)|^2+|u(x,1)|^2)\,{\rm d}x\Big)^{\frac{1}{2}}.
$$
Using the anisotropic Sobolev inequality, one has
$$
\|(u,v)\|_{L^{\infty}}\leq C\Big(\|(u,v)\|_{L^2}+\|\nabla (u,v)\|_{L^2}+\|(u_{xy},v_{xy})\|_{L^2}\Big).
$$
Moreover, it follows from the Poincar\'{e} inequality and Lemma \ref{So1} that
\begin{align*}
\|v\|_{L^2}\leq C\|v_{x}\|_{L^2},\quad \|u\|_{L^2}\leq C\|\nabla u\|_{L^2}+C\Big(\int_{\mathbb{T}}|u(x,0)|^2\,{\rm d}x\Big)^{\frac{1}{2}}.
\end{align*}
Thus,
$$
\|(u,v)\|_{L^{\infty}}\leq C\varepsilon^{-3}\|(u,v)\|_{\mathbb{X}}.
$$
Recalling \eqref{3.9} and choosing $k_{0}\geq 7\tilde{\alpha}$, we obtain from \eqref{9.25} that
\begin{align*}
\|(u,v)\|_{\mathbb{X}}\leq C\varepsilon^{6}+C\varepsilon^{-5}\|(\tilde{u},\tilde{v})\|_{\mathbb{X}}^2.
\end{align*}
Then the mapping
$$
(\tilde{u},\tilde{v})\mapsto (u,v),
$$
maps the ball $\{(u,v):\|(u,v)\|_{\mathbb{X}}\leq 2C\varepsilon^{6}\}$ in $\mathbb{X}$ into itself provided that $\varepsilon$ is suitably small.

Furthermore, for every two pairs $(\tilde{u}_{1},\tilde{v}_{1})$ and $(\tilde{u}_{2},\tilde{v}_{2})$ in this ball, let $(u_{1},v_{1})$ and $(u_{2},v_{2})$ be the corresponding solutions of \eqref{9.24}, respectively. Then we have
\begin{align*}
&\|(u_{1}-u_{2},v_{1}-v_{2})\|_{\mathbb{X}}\nonumber\\
&\leq C\varepsilon^{-5}(\|(\tilde{u}_{1},\tilde{v}_{1})\|_{\mathbb{X}}+\|(\tilde{u}_{2},\tilde{v}_{2})\|_{\mathbb{X}})\|(\tilde{u}_{1}-\tilde{u}_{2},\tilde{v}_{1}-\tilde{v}_{2})\|_{\mathbb{X}}.
\end{align*}
Thus, for each fixed integer $k_{0}\geq 7\tilde{\alpha}$, there exist $\varepsilon_{0}$ and $\sigma_{0}$ (required only when $\alpha\geq 1$) such that for any $\varepsilon\in (0,\varepsilon_{0})$, $\sigma\in [0,\sigma_{0})$, the mapping 
$$
(\tilde{u},\tilde{v})\mapsto (u,v),
$$
maps the ball $\{(u,v):\|(u,v)\|_{\mathbb{X}}\leq 2C\varepsilon^{6}\}$ in $\mathbb{X}$ into itself and is a contraction mapping. Therefore, the remainder equation \eqref{3.8} has a unique solution $(u,v)$ satisfying
$$
\|(u,v)\|_{\mathbb{X}}\leq 2C\varepsilon^{6}.
$$
Hence, there holds $\|(u,v)\|_{L^{\infty}}\leq C\varepsilon^{3}$. This completes the proof of Proposition \ref{prop4}. $\hfill\square$

Now we conclude the proof of main theorems.

\textbf{Proof of Theorems \ref{thm1}--\ref{thm3}}. Theorems \ref{thm1}--\ref{thm2}, and the $H^{1}$ estimate in \eqref{vorticity} follow directly from Lemma \ref{lem5.1}, Propositions \ref{prop3} and \ref{prop4}. For the $L^{\infty}$ estimate in \eqref{vorticity}, by the Sobolev embedding $H^{1}(\mathbb{T})\hookrightarrow L^{\infty}(\mathbb{T})$, we can choose $y_{-}\in (d/2,d)$ and $y_{+}\in (1-d,1-d/2)$ such that
\begin{align}\label{vorticity2}
    \|(\omega^{\var}-A)(\cdot,y_{\pm})\|_{L^{\infty}(\mathbb{T})}\leq C_{d}\varepsilon^{\frac{1}{\tilde{\alpha}}}.
\end{align}
Let $\tilde{\Omega}=\mathbb{T}\times (y_{-},y_{+})$. Noting that
$$
\varepsilon^2\Delta(\omega^{\var}-A)=u^{\varepsilon}\partial_{x}(\omega^{\var}-A)+v^{\varepsilon}\partial_{y}(\omega^{\var}-A),
$$
we obtain the $L^{\infty}$ estimate in \eqref{vorticity} from \eqref{vorticity2} and the maximum principle on $\tilde{\Omega}$. This completes the proof of Theorem \ref{thm3}. $\hfill\square$

\appendix
\renewcommand{\appendixname}{Appendix~\Alph{section}}

\section{Constructions of Auxiliary Functions}\label{AppendixB}

In this section, we give explicit constructions of auxiliary functions $h_{u}$ and $h_{v}$ in \eqref{3.5-0} with \eqref{3.5-2}--\eqref{3.5-3}. Consider the following system:
\begin{align}\label{B1}	\left\{\begin{aligned}
		&\partial_{x}g_{1}+\partial_{y}g_{2}=F(x,y),\\
		&g_{2}(x,0)=G_{1}(x),\quad g_{2}(x,1)=G_{2}(x),\\
		&(g_{1}-\varepsilon^{\alpha}\partial_{y}g_{1})(x,0)=H_{1}(x),\quad (g_{1}+\varepsilon^{\alpha}\partial_{y}g_{1})(x,1)=H_{2}(x),
	\end{aligned}
	\right.
\end{align}
where $F(x,y)$, $G_{1}(x)$, $G_{2}(x)$, $H_{1}(x)$ and $H_{2}(x)$ are the known smooth periodic functions, satisfying the following condition:
\begin{align}\label{B1-1}
\begin{aligned}
&F(x,0)=F(x,1)=0,\,\,\partial_{y}F(x,0)=\partial_{y}F(x,1)=0,\\
&\int_{\mathbb{T}}\int_{0}^{1}F(x,y)\,{\rm d}y{\rm d}x=\int_{\mathbb{T}}\big(G_{2}-G_{1}\big)(x)\,{\rm d}x.
\end{aligned}
\end{align}

By using the Fourier transform with respect to $x$, for any smooth $2\pi-$periodic function $f(x,y)$, we have
$$
f(x,y)=\sum\limits_{k\in \mathbb{Z}}f_{k}(y)e^{ikx},\quad f_{k}(y)=\frac{1}{2\pi}\int_{0}^{2\pi}f(x,y)e^{-ikx}\,{\rm d}x.
$$
Then, applying the Fourier transform to $\eqref{B1}_{1}$, we have
\begin{align}\label{B4}
	ikg_{1,k}(y)+\partial_{y}g_{2,k}(y)=F_{k}(y).
\end{align}
For zero mode $k=0$, we can choose
$$
g_{2,0}(y)=\int_{0}^{y}F_{0}(s)\,{\rm d}s+\frac{1}{2\pi}\int_{0}^{2\pi}G_{1}(x)\,{\rm d}x,
$$
and $g_{1,0}(y)=k_{1}+k_{2}y$ to satisfy
$$
g_{1,0}(0)-\varepsilon^{\alpha}\partial_{y}g_{1,0}(0)=H_{1,0},
$$
and
$$
g_{1,0}(1)+\varepsilon^{\alpha}\partial_{y}g_{1,0}(1)=H_{2,0},
$$
that is,
$$
k_{1}=\frac{H_{1,0}+\varepsilon^{\alpha}(H_{1,0}+H_{2,0})}{1+2\varepsilon^{\alpha}},\quad k_{2}=\frac{H_{2,0}-H_{1,0}}{1+2\varepsilon^{\alpha}}.
$$

For nonzero mode $k\neq 0$, we have
\begin{align}\label{B5}
	g_{1,k}(y)=\frac{1}{ik}\Big(F_{k}(y)-\partial_{y}g_{2,k}(y)\Big),\quad k\neq 0,
\end{align}
and
\begin{align}\label{B5-1}
	\partial_{y}g_{1,k}(y)=\frac{1}{ik}\Big(\partial_{y}F_{k}(y)-\partial_{y}^2g_{2,k}(y)\Big),\quad k\neq 0,
\end{align}
which, together with the boundary conditions $\eqref{B1}_{2}$--$\eqref{B1}_{3}$, yields
\begin{align}\label{B6-0}
	g_{2,k}(0)=G_{1,k},\quad g_{2,k}(1)=G_{2,k}.
\end{align}
\begin{align}\label{B6}
	\varepsilon^{\alpha}\partial_{y}^2g_{2,k}(0)-\partial_{y}g_{2,k}(0)+\big(F_{k}(0)-\varepsilon^{\alpha}\partial_{y}F_{k}(0)\big)=ikH_{1,k},
\end{align}
and
\begin{align}\label{B6-1}
	-\varepsilon^{\alpha}\partial_{y}^2g_{2,k}(1)-\partial_{y}g_{2,k}(1)+\big(F_{k}(1)+\varepsilon^{\alpha}\partial_{y}F_{k}(1)\big)=ikH_{2,k}.
\end{align}
Using \eqref{B1-1}, one has $F_{k}(0)=F_{k}(1)=\partial_{y}F_{k}(0)=\partial_{y}F_{k}(1)=0$, so \eqref{B6}--\eqref{B6-1} can be rewritten as
\begin{align}\label{B6-2}
	\varepsilon^{\alpha}\partial_{y}^2g_{2,k}(0)-\partial_{y}g_{2,k}(0)=ikH_{1,k},\quad -\varepsilon^{\alpha}\partial_{y}^2g_{2,k}(1)-\partial_{y}g_{2,k}(1)=ikH_{2,k}.
\end{align}
To find explicit formula for $g_{2,k}(y)$, we assume
$$
g_{2,k}(y)=A_{k}+B_{k}y+C_{k}y^2+D_{k}y^3.
$$
It follows from \eqref{B6-0}, \eqref{B6-2} and direct calculations that
\begin{align*}
	&A_{k}=G_{1,k},\quad B_{k}=\frac{2\varepsilon^{\alpha}(G_{2,k}-G_{1,k})}{1+2\varepsilon^{\alpha}}-\frac{ikH_{1,k}}{1+2\varepsilon^{\alpha}}+\frac{2\varepsilon^{\alpha}\big[2(G_{2,k}-G_{1,k})+ik(H_{1,k}+H_{2,k})\big]}{(1+6\varepsilon^{\alpha})(1+2\varepsilon^{\alpha})},\\ &C_{k}=\frac{G_{2,k}-G_{1,k}}{1+2\varepsilon^{\alpha}}+\frac{ikH_{1,k}}{1+2\varepsilon^{\alpha}}+\frac{2(G_{2,k}-G_{1,k})+ik(H_{1,k}+H_{2,k})}{(1+6\varepsilon^{\alpha})(1+2\varepsilon^{\alpha})},\\ &D_{k}=-\frac{2(G_{2,k}-G_{1,k})+ik(H_{1,k}+H_{2,k})}{1+6\varepsilon^{\alpha}}.
\end{align*}
Once $g_{2,k}$ is determined, $g_{1,k}$ follows directly from \eqref{B5}.

\section{Remark on General Slip Length}
In this appendix, motivated by \cite[Section 4]{WWX-2011}, we briefly sketch the main idea for generalizing the above analysis to the case that $\alpha$ is any nonnegative rational number
in the subcritical and supercritical regimes.

Since the discussions in subcritical and supercritical cases are similar, we only take the supercritical case $\alpha>1$ for example. Let $r>0$ such that $a=\frac{\alpha}{r}$ and $b=\frac{1}{r}$ are co-prime integers. Then we construct the interior expansions as
\begin{align*}
	u_{I}(x,y)=\sum\limits_{j=0}^{k_{0}}\varepsilon^{rj}u_{I}^{j}(x,y),\quad v_{I}(x,y)=\sum\limits_{j=0}^{k_{0}}\varepsilon^{rj}v_{I}^{j}(x,y),\quad p_{I}(x,y)=\sum\limits_{j=0}^{k_{0}}\varepsilon^{rj}p_{I}^{j}(x,y),
\end{align*}
and the lower boundary expansions as
\begin{align*}
	u_{B}(x,z)=\sum\limits_{j=0}^{k_{0}}\varepsilon^{rj}u_{B}^{j}(x,z),\quad v_{B}(x,z)=\sum\limits_{j=0}^{k_{0}+b}\varepsilon^{rj}v_{B}^{j}(x,z),\quad p_{B}(x,z)=\sum\limits_{j=0}^{k_{0}+b}\varepsilon^{rj}p_{B}^{j}(x,z),
\end{align*}
and the upper boundary expansions as
 \begin{align*}
 	\hat{u}_{B}(x,\zeta)=\sum\limits_{j=0}^{k_{0}}\varepsilon^{rj}\hat{u}_{B}^{j}(x,\zeta),\quad \hat{v}_{B}(x,\zeta)=\sum\limits_{j=0}^{k_{0}+b}\varepsilon^{rj}\hat{v}_{B}^{j}(x,\zeta),\quad \hat{p}_{B}(x,\zeta)=\sum\limits_{j=0}^{k_{0}+b}\varepsilon^{rj}\hat{p}_{B}^{j}(x,\zeta).
 \end{align*}
Then, similarly, one has $(u_{I}^{0},v_{I}^{0},p_{I}^{0})$ satisfy the Euler equations with boundary conditions
$$
v_{I}^{0}(x,0)\equiv v_{I}^{0}(x,1)=0,
$$
and $(u_{I}^{j},v_{I}^{j},p_{I}^{j})$ ($j\geq 1$) satisfy the linearized Euler equation with the terms $\Delta (u_{I}^{j-2\tilde{\alpha}},v_{I}^{j-2\tilde{\alpha}})$ replaced by $\Delta (u_{I}^{j-2b},v_{I}^{j-2b})$ and the boundary conditions
$$
v_{I}^{j}(x,0)=-v_{B}^{j}(x,0),\quad v_{I}^{j}(x,1)=-\hat{v}_{B}^{j}(x,0).
$$

For the boundary-layer profiles $(u_{B}^{j},v_{B}^{j})$, we observe from the divergence free condition that
\begin{align*}
\partial_{z}v_{B}^{j}\equiv 0\,\,\forall j\leq b-1\Rightarrow v_{B}^{j}\equiv 0\,\,\forall j\leq b-1,
\end{align*}
which, together with \eqref{2.13}, yields that
\begin{align*}
	v_{I}^{j}\equiv 0\,\,\forall j\leq b-1.
\end{align*}
Moreover, we observe from the equation of vertical velocity that
\begin{align*}
	\partial_{z}p_{B}^{j}\equiv 0\,\,\forall j\leq b-1.
\end{align*}
For $(u_{B}^{0},v_{B}^{b},p_{B}^{b})$, it still satisfies the Prandtl equation with boundary condition
$$
u_{B}^{0}(x,0)=f_{0}(x)-u_{I}^{0}(x,0),
$$
and the far-field conditions
\begin{align*}
	\lim\limits_{z\to \infty}(u_{B}^{0},v_{B}^{b},p_{B}^{b})=(0,0,0).
\end{align*}
And $(u_{B}^{k},v_{B}^{k+b},p_{B}^{k+b})$ $(k\geq 1)$ satisfies the linearized Prandtl equation with boundary conditions
$$
u_{B}^{k}(x,0)=-u_{I}^{k}(x,0)+\partial_{z}u_{B}^{k+b-a}(x,0)+\partial_{y}u_{I}^{k-a}(x,0).
$$
The upper boundary profiles can similarly be treated as in the lower boundary profiles. By similar arguments as in Sections 2--5, we can establish the existence of the interior and boundary-layer profiles, and construct our approximate solution as in Proposition \ref{prop3}. Finally, using the arguments in the remainder equation, we can derive the existence and uniform estimates of remainder equations. Moreover, by the above expansion, we have following convergence rate in $L^{\infty}$:
\begin{align*}
    \Big\|u^{\varepsilon}(x,y)-(Ay+B)-u_{B}^{0}\Big(x,\frac{y}{\varepsilon}\Big)\chi(y)-\hat{u}_{B}^{0}\Big(x,\frac{y-1}{\varepsilon}\Big)\tilde{\chi}(y)\Big\|_{L^{\infty}}+\|v^{\varepsilon}(x,y)\|_{L^{\infty}}\leq C\varepsilon^{r}.
\end{align*}
The same convergence rate $O(\varepsilon^{r})$ holds in the subcritical regime for every rational exponent $0<\alpha<1$, in which $u_{B}^{0}=\hat{u}_{B}^{0}=0$. In particular, when $1/\alpha$ is a positive integer, the convergence rate will be $O(\varepsilon^{\alpha})$.

\begin{remark}
    By using similar construction as in \cite[Section 4.2]{WWX-2011}, we can also justify the validity of steady Prandtl expansion in the subcritical regime when $\alpha$ is a general nonnegative irrational number. However, when $\alpha<0$ or $\alpha$ is an irrational number in the supercritical case, some technical calculations are needed, so we leave the cases $\alpha<0$ and irrational $\alpha>1$ for future work.
\end{remark}

\bigskip

\noindent\textbf{Acknowledgments}: R. M. Chen was supported in part by the NSF grant DMS-2205910. D. Wang was supported in part by NSF grant DMS-2510532.

\smallskip





\noindent\textbf{AI Disclosure}: No generative AI tools were used in this paper.


\bigskip


\begin{thebibliography}{99}

\bibitem{AWXY} R. Alexandre, Y. Wang, C. Xu, and T. Yang, Well-posedness of the Prandtl equation in Sobolev spaces. J. Am. Math. Soc. 28 (2015), 745--784.
	
	\bibitem{B-1956} G. K. Batchelor, On steady laminar flow with closed streamlines at large Reynolds number, J. Fluid Mech. 1 (1956) 177--190.

    \bibitem{Beirao-Crispo} H. Beir\~ao da Veiga and F. Crispo, A survey on some vanishing viscosity limit results, Adv. Nonlinear Anal. 12 (2023) 28 pp.

    \bibitem{Berselli-Spirito} L. C. Berselli and S. Spirito, On the vanishing viscosity limit of 3D Navier-Stokes equations under slip boundary conditions in general domains, Commun. Math. Phys. 316 (2012) 171--198.

    \bibitem{Chemetov-Cirpriano} N. V. Chemetov and F. Cipriano, Inviscid limit for Navier-Stokes equations in domains with permeable boundaries, Appl. Math. Lett. 33 (2014) 6--11.

    \bibitem{CMR} T. Clopeau, A. Mikeli\'{c} and R. Robert, On the vanishing viscosity limit for the 2D incompressible Navier-Stokes equations with the friction type boundary conditions, Nonlinearity 11 (1998) 1625--1636.

    \bibitem{DIN} T. D. Drivas, S. Iyer and T. T. Nguyen, The Feynman--Lagerstrom criterion for boundary layers, Arch. Ration. Mech. Anal., 248 (2024), 41 pp.

    \bibitem{EHSX} T. M. Elgindi, Y. Huang, A. Said, C. Xie, A classification theorem for steady Euler flows. Duke Math. J. Advance Publication 1--35 (2026).
	
	\bibitem{FTZ-2018} M. Fei, T. Tao, and Z. Zhang, On the zero-viscosity limit of the Navier–Stokes equations in $\R_{+}^3$ without analyticity, J. Math. Pures Appl. 112 (2018) 170--229.
	
	\bibitem{FGLT-2023} M. Fei, C. Gao, Z. Lin, T. Tao, Prandtl--Batchelor flows on a disk, Commun. Math. Phys. 397 (3) (2023) 1103--1161.
	
	\bibitem{FGLT-2024} M. Fei, C. Gao, Z. Lin, T. Tao, Prandtl--Batchelor flows on an annulus. Adv. Math. 458 (2024), part B, 86 pp.
	
	\bibitem{FPZ-2025} M. Fei, X. Pan, and J. Zhao, Vanishing viscosity limit of the Navier-Stokes equations in a horizontally periodic strip, Sci China Math. 69 (2026), 2337--2372.

    \bibitem{FL-1956} R. P. Feynman and P. A. Lagerstrom, Remarks on high Reynolds number flows in finite domains, Proc. IX International Congress on Applied Mechanics, Brussels, 3 (1956), 342--343.

    \bibitem{GZ-2023} C. Gao and L. Zhang, On the steady Prandtl boundary layer expansions. Sci. China Math. 66 (2023), 1993--2020.

    \bibitem{GZ-2024} C. Gao and L. Zhang, The steady Prandtl boundary layer expansions for non-shear Euler flow. J. Math. Phys. 65 (2024), 27 pp.

    \bibitem{GL-2025} Y. Ge and C. Liu, Prandtl boundary layer expansions over a moving periodic plate. Commun. Pure Appl. Anal. 31 (2026), 86--121.

    \bibitem{G-D} D. G\'{e}rard-Varet and E. Dormy, On the ill-posedness of the Prandtl equation. J. Am. Math. Soc. 23 (2010), 591--609.

    \bibitem{GLNR} D. G\'{e}rard-Varet, C. Lacave, T.T. Nguyen and F. Rousset, The vanishing viscosity limit for 2D Navier-Stokes in a rough domain, J. Math. Pures Appl. 119 (2018) 45--84.

    \bibitem{GM-2015} D. G\'{e}rard-Varet and N. Masmoudi, Well-posedness for the Prandtl system without analyticity or monotonicity. Ann. Sci. \'{E}c. Norm. Sup\'{e}r. 48  (2015), 1273--1325.


    \bibitem{GM-2019} D. G\'{e}rard-Varet and Y. Maekawa, Sobolev stability of Prandtl expansions for the steady Navier-Stokes equations, Arch. Ration. Mech. Anal. 233 (2019), 1319--1382. 

    \bibitem{GGT} E. Grenier, Y. Guo and T. Nguyen, Spectral instability of Prandtl boundary layers: an overview. Analysis (Berlin) 35 (2015), 343--355.

    \bibitem{GT} E. Grenier and T. Nguyen, $L^{\infty}$ instability of Prandtl layers. Ann. PDE 5 (2019), 36 pp.

\bibitem{GXX} C. Gui, C. Xie, and H. Xu, A selection principle for 2D steady Euler flows via the vanishing viscosity limit. arXiv:2601.08647.

\bibitem{GXX-2026} C. Gui, C. Xie, and H. Xu, On a classification of steady solutions to two-dimensional Euler equations. Comm. Math. Phys. 407 (2026), 31 pp.

    \bibitem{GI-2021} Y. Guo and S. Iyer, Regularity and Expansion for Steady Prandtl Equations. Commun. Math. Phys. 382 (2021). 1403--1447.

    \bibitem{GI-2023} Y. Guo and S. Iyer, Validity of Steady Prandtl Layer Expansions. Comm. Pure Appl. Math. 76 (2023), 3150--3232.

    \bibitem{GT-2011} Y. Guo and T.T. Nguyen, A note on Prandtl boundary layer, Comm. Pure Appl. Math. 64 (2011), 1416--1438.
	
	\bibitem{GT-2017} Y. Guo and T.T. Nguyen, Prandtl boundary layer expansions of steady Navier–Stokes flows over a moving plate, Ann. PDE 3 (2017).

    \bibitem{GW-arXiv} Y. Guo and Y. Wang, Steady Prandtl Expansion for Subsonic Flows. arXiv:2501.13360.


    \bibitem{HN-2017} F. Hamel and N. Nadirashvili, Shear flows of an ideal fluid and elliptic equations in unbounded domains. Comm. Pure Appl. Math. 70 (2017), 590--608.

     \bibitem{HN-2019} F. Hamel and N. Nadirashvili, A Liouville theorem for the Euler equations in the plane.  Arch. Ration. Mech. Anal., 233 (2019), 599--642.

     \bibitem{HN-2023} F. Hamel and N. Nadirashvili, Circular flows for the Euler equations in two-dimensional annular domains, and related free boundary problems. 
J. Eur. Math. Soc. 25 (2023), 323--368.
	
	\bibitem{IP-2006} D. Iftimie and G. Planas, Inviscid limits for the Navier–Stokes equations with Navier friction boundary conditions, Nonlinearity 19 (2006) 899--918.
	
	\bibitem{IS-2011} D. Iftimie and F. Sueur, Viscous boundary layer for the Navier–Stokes equations with the Navier slip conditions, Arch. Ration. Mech. Anal. 199 (2011) 145--175.

    \bibitem{Iyer-2019} S. Iyer, Steady Prandtl Layers over a Moving Boundary: Non-Shear Euler Flows. SIAM J. Math. Anal., 51 (2019),1657--1695.

    \bibitem{Iyer-2017} S. Iyer, Steady Prandtl boundary layer expansions over a rotating disk. Arch. Ration.
Mech. Anal. 224 (2017), 421--469.

\bibitem{Iyer-PK} S. Iyer, Global steady Prandtl expansion over a moving boundary ({\uppercase\expandafter{\romannumeral1}}, {\uppercase\expandafter{\romannumeral2}}, {\uppercase\expandafter{\romannumeral3}}). Peking Math. J. 2 (2019), no. 2, 155--238; 2 (2019), no. 3--4, 353--437; 3 (2020), no. 1, 47--102.


\bibitem{LMY} W. Li, N. Masmoudi and T. Yang, Well-posedness in Gevrey function space for 3D Prandtl equations without structural assumption. Comm. Pure Appl. Math. 75 (2022), 1755--1797.


	
	
	
	\bibitem{FMT-2008} M. C. Lopes Filho, A. L. Mazzucato, and M.E. Taylor, Vanishing viscosity limits and boundary layers for circularly symmetric 2D flows, Bull. Braz. Math. Soc. 39 (2008) 471--513.
	
	\bibitem{LNP} M. C. Lopes Filho, H. J. Nussenzveig Lopes and G. Planas, On the inviscid limit for two dimensional incompressible flow with Navier friction condition, SIAM J. Math. Anal., 36, 1130--1141, 2005.
	
	
	\bibitem{JM-2001} W. J\"{a}ger and A. Mikeli\'{c}, On the roughness-induced effective boundary conditions for an incompressible viscous flow,
	J. Differ. Equations, 170 (2001), 96--122.
	
	
	\bibitem{M-2014} Y. Maekawa, On the inviscid limit problem of the vorticity equations for viscous incompressible flows in the half-plane. Comm. Pure Appl. Math. 67 (2014), 1045--1128.

    \bibitem{MM-2017} Y. Maekawa and A. L. Mazzucato, The inviscid limit and boundary layers for Navier-Stokes flows. In: Giga, Y., Novotny, A. (eds) Handbook of Mathematical Analysis in Mechanics of Viscous Fluids. Springer, Cham. 2016.

    \bibitem{MW-2015} N. Masmoudi and T. K. Wong, Local in time existence and uniqueness of solutions to the Prandtl equations by energy methods. Commun. Pure. Appl. Math. 68 (2015), 1683--1741.
    
	
	\bibitem{Maxwell} J. C. Maxwell, On stresses in rarefied gases arising from inequalities of temperature, Philos. Trans. R. Soc. London
	170 (1879), 231--256.
	

    \bibitem{N} C. L. M. H. Navier, M\'{e}moire sur les lois du mouvement des fluides, M\'{e}moires de l'Acad\'{e}mie des sciences, vol. 6 (1827), 389--440.

    \bibitem{NP} J. Neustupa and P. Penel, On the vanishing viscosity limit of strong solutions of the Navier-Stokes equation with slip boundary conditions, in: Mathematical Analysis on the Navier-Stokes Equations and Related Topics, Past and Future, in: GAKUTO Internat. Ser. Math. Sci. Appl., vol. 35, Gakk$\bar{o}$tosho, Tokyo, 2011, pp. 130--150.

    \bibitem{Nguyen-2020} T. Nguyen, The inviscid limit of Navier-Stokes with critical Navier-slip boundary conditions for analytic data, SIAM J. Math. Anal. 52 (2020) 4616--4637.

    \bibitem{Oleinik-1967} O. A. Oleinik, On the mathematical theory of boundary layer for an unsteady flow of incompressible fluid. J.
Appl. Math. Mech. 30 (1966), 951--974.

\bibitem{OS-1999} O. A. Oleinik and V. N. Samokhin, Mathematical models in boundary layer theory. Applied Mathematics and Mathematical Computation, 15. Chapman and Hall/CRC, Boca Raton, FL, 1999.

\bibitem{PZ-2021} M. Paicu and P. Zhang, Global existence and the decay of solutions to the Prandtl system with small analytic
data. Arch. Ration. Mech. Anal. 241 (2021), 403--446.

\bibitem{PZZ-2026} X. Pan, C. Zhao and J. Zhao, Vanishing viscosity limit of the stationary Navier-Stokes equations in an expanding annulus. arXiv:2610.17754.

\bibitem{PZ-2025} X. Pan and J. Zhao, Vanishing viscosity limit of the stationary Navier-Stokes equations with
the Navier-slip boundary and its application, J. Differential Equations., 443 (2025), 58 pp.

\bibitem{PZ-2025-2} X. Pan and J. Zhao, Vanishing viscosity limit of the 2D stationary Navier-Stokes equations outside a rotating disc and its application. arXiv:2511.12953.
	
	\bibitem{P-1904} L. Prandtl, \"{U}ber Fl\"{u}ssigkeitsbewegung bei sehr kleiner Reibung, in: Actes du 3me Congr\'{e}s international des Math\'{e}maticiens, Teubner Leipzig, 1905, 484--491.


    \bibitem{QDW} Y. Qin, X. Dong and X. Wang, Prandtl Equations and Related Boundary Layer Equations, Springer, Singapore, 2024.

	
	
	\bibitem{SC-1998-1} M. Sammartino and R. E. Caflisch, Zero viscosity limit for analytic solutions of the Navier-Stokes equation on a half-space. I. Existence for Euler and Prandtl equations, Comm. Math. Phys. 192 (1998) 433--461.
	
	\bibitem{SC-1998-2} M. Sammartino and R. E. Caflisch, Zero viscosity limit for analytic solutions of the Navier-Stokes equation on a half-space. II. Construction for the Navier–Stokes solution, Comm. Math. Phys. 192 (1998) 463--491.

    \bibitem{SG-2018} H. Schlichting and K. Gersten, Boundary-Layer Theory, ninth edition, Springer-Verlag, Berlin, 2017, With contributions from Egon Krause and Herbert Oertel Jr., Translated from the German by Katherine Mayes.

    \bibitem{ST-2027} J. Song and T. Tao, Vanishing viscosity limit of the forced steady Navier-Stokes Equations with slip boundary condition in a finite channel. Discrete Contin. Dyn. Syst. 59 (2027), 331--371.

    \bibitem{TWZ} T. Tao, W. Wang and Z. Zhang, Zero-viscosity limit of the Navier-Stokes equations with the Navier friction boundary condition, SIAM J. Math. Anal. 52 (2020) 1040--1095.

    \bibitem{TW-2002} R. Temam and X. Wang, Boundary layers associated with Incompressible Navier-Stokes equations: the noncharacteristic
boundary case. J. Differ. Equations. 179 (2002), 647--686.
	
	\bibitem{WWZ-2017} C. Wang, Y. Wang and Z. Zhang, Zero-viscosity limit of the Navier–Stokes equations in the analytic setting, Arch. Ration. Mech. Anal. 224 (2017) 555--595.
	
	
	\bibitem{WWX-2007} X. Wang, Y. Wang and Z. Xin, Boundary layers of incompressible fluids with Navier
	boundary conditions, Proceedings of ICCM 2007, Higher Edu. Press, 3, 145--159, 2007.
	
	\bibitem{WWX-2011} X. Wang, Y. Wang and Z. Xin, Boundary layers in incompressible Navier-Stokes equations with Navier boundary conditions for the vanishing viscosity limit, Commun. Math. Sci. 8 (2010), 965--998.

    \bibitem{WZ-2021} Y. Wang and Z. Zhang, Global $C^{\infty}$ regularity of the steady Prandtl equation with favorable pressure gradient. Ann. Inst. H. Poincar\'{e} C Anal. Non Lin\'{e}aire, 38 (2021), no. 6, 1989--2004.
	
	
	\bibitem{WYZ-2017} Y. Wang, J. Yin, and S. Zhu, Vanishing viscosity limit for incompressible Navier-Stokes equations with
	Navier boundary conditions for small slip length, J. Math. Phys. 58 (2017), 18 pp.

    \bibitem{WXZ} L. Wang, Z. Xin and A. Zang, Vanishing viscous limits for 3D Navier-Stokes equations with a Navier-slip boundary condition, J. Math. Fluid Mech. 14 (2012) 791--825.

    \bibitem{Wood} W. W. Wood, Boundary layers whose streamlines are closed. J. Fluid Mech., 2 (1957), 77--87. 
	
	\bibitem{XX-2007} Y. Xiao and Z. Xin, On the vanishing viscosity limit for the 3D Navier–Stokes equations with a slip boundary condition, Comm. Pure Appl. Math. 60 (2007) 1027--1055.

    \bibitem{XZ-2004} Z. Xin and L. Zhang, On the global existence of solutions to the Prandtl’s system. Adv. Math. 181 (2004), 88--133.

    \bibitem{Zhong} X. Zhong, Vanishing viscosity limits for the 3D Navier-Stokes equations with a slip boundary condition, Proc. Am. Math. Soc. 145 (2017) 1615--1628.
	
\end{thebibliography}
    \end{document}